\documentclass[english]{article}
\usepackage[T1]{fontenc}
\usepackage[utf8]{inputenc}
\usepackage{geometry}
\usepackage{color}
\usepackage{booktabs}
\usepackage{babel}
\usepackage{amsmath}
\usepackage{amssymb}
\usepackage{graphicx}
\usepackage{amsfonts}\usepackage{epsfig}\usepackage{graphics}\usepackage{oldgerm}\relax
\usepackage{color}
\usepackage{siunitx}
\usepackage{svg}
\usepackage{color}
\usepackage[unicode=true,pdfusetitle,bookmarks=true,bookmarksnumbered=false,bookmarksopen=false,breaklinks=false,pdfborder={0 0 1},backref=false,colorlinks=false]{hyperref}
\hypersetup{
  colorlinks   = true, %Colours links instead of ugly boxes
  urlcolor     = black, %Colour for external hyperlinks
  linkcolor    = black, %Colour of internal links
  citecolor   = black %Colour of citations
}

\makeatletter

\providecommand{\tabularnewline}{\\}
\newtheorem{theo}{Theorem}[section]
\newtheorem{prop}[theo]{Proposition}

\newtheorem{lem}[theo]{Lemma}

\newtheorem{rem}[theo]{Remark}

\newcommand{\Margo}[1]{\textcolor{red}{#1}}

\title{Augmented Skew-Symetric System\\
for Shallow-Water System\\
with Surface Tension Allowing\\
Large Gradient of Density}

\author{D. Bresch, N. Cellier, F. Couderc, M. Gisclon,\\
P. Noble, G.--L. Richard, C. Ruyer-Quil, J.--P. Vila}
\makeatother

\begin{document}
\maketitle
\begin{abstract}
    In this paper, we introduce a new extended version of the shallow-water equations with surface tension which may be decomposed into
    a hyperbolic part and a second order derivative part which is skew-symmetric with respect to the $L^{2}$ scalar product. This reformulation allows for large gradients of fluid height simulations using a splitting method. This result is a generalization of the results published by P. Noble and J.--P. Vila in {[}\textit{SIAM J. Num. Anal.} (2016){]} and by D.~Bresch, F. Couderc, P.~Noble and J.P. Vila in {[}\textit{C.R. Acad. Sciences Paris} (2016){]} which are restricted to quadratic
    forms of the capillary energy respectively in the one dimensional
    and two dimensional setting. This is also an improvement of the results
    by J.~Lallement, P. Villedieu \textit{et al.} published in {[}AIAA Aviation
    Forum 2018{]} where the augmented version is not skew-symetric with
    respect to the $L^{2}$ scalar product. Based on this new formulation,
    we propose a new numerical scheme and perform a nonlinear stability
    analysis. Various numerical simulations of the shallow water equations
    are presented to show differences between quadratic (w.r.t the gradient
    of the height) and general surface tension energy when high gradients
    of the fluid height occur.
\end{abstract}

\section{Introduction}
In this paper, we consider compressible Euler type equations with capillarity (such as the shallow-water system with surface tension),  in the two-dimensional setting, issued from Hamiltonian formulation in the spirit of P. Casal and H. Gouin (\cite{Casal-Gouin}) (see also D.~Serre~(\cite{Serre_1993}).
There exists a large body of literature on various numerical techniques for simulating the shallow-water equations without capillarity terms for a variety of applications, such as discontinuous Galerkin methods (e.g., Giraldo \cite{Gi}, Giraldo et al. \cite{Gi2}, Eskilsson and Sherwin \cite{Es}, Nair et al. \cite{Na}, Xing et al. \cite{Xi}, Blaise and St-Cyr \cite{Bl}), in addition to spectral methods (e.g., Giraldo and Warburton \cite{Gi3}, Giraldo \cite{Gi1}), and purely Lagrangian approaches (e.g., Frank and Reich \cite{Fr}, Capecelatro \cite{Ca}). When dispersive effect is included, everything change and numerous attempts have been conducted to try to get rid spurious currents (also known as parasitic currents) generated at the free surface due to the presence of the third-order term coming from the capillarity quantity. Augmented versions have been proposed to decrease the level of derivative in the system but with not enough properties to allow to design an efficient numerical method to compute for instance large gradient of density. This is the objective of our paper to propose an appropriate extended formulation which allows an appropriate splitting method.
To be more precise, let us define the internal energy $E$ as follows
\begin{equation}
    E \left(h,\boldsymbol{p}\right)=\Phi\left(h\right)+\sigma\left(h\right)\mathcal{E_{\mathrm{cap}}}\left(\|\boldsymbol{p}\|\right)\label{E}
\end{equation}
with $h$ the density of the fluid (or the fluid height if we consider the shallow-water system), $\boldsymbol{p}=\nabla h$ and $\Phi(h)$ the pressure contribution and $\sigma(h)\mathcal{E_{\mathrm{cap}}}(\|{\boldsymbol{p}}\|)$ the capillarity energy (Note that $s\mapsto \Phi(s)$, $s\mapsto \sigma(s)$ and  $s \mapsto \mathcal{E_{\mathrm{cap}}}(s)$ are three given positive scalar functions). We consider the
following system
\begin{equation}
    \left\{ \begin{array}{cc}
        \partial_{t}h+{\rm div}\left(h\mathbf{u}\right)=0                                                                                                                                                                          & (i)  \\
        \partial_{t}\left(h\mathbf{u}\right)+{\rm div}\left(h\mathbf{u}\otimes\mathbf{u}\right)+\nabla P=-{\rm div}\left(\nabla h\otimes\nabla_{\mathbf{p}}E\right)+\nabla\left(h{\rm div}\left(\nabla_{\mathbf{p}}E\right)\right) & (ii)
    \end{array}\right.\label{eq:Sys_SW_gen}
\end{equation}
which is obtained as the Euler-Lagrange equation related to the total energy (sum of the kinetic $h|\mathbf{u}|^2/2$ and internal energy $E(h,\mathbf{p})$ given by \eqref{E}) under the conservation of mass constraint  with ${\mathbf{u}}$ the fluid velocity vector field and $P$ the pressure law
given by
\begin{equation}
    P(h,\boldsymbol{p}):=h \, \partial_{h}E(h,\boldsymbol{p})-E(h,\boldsymbol{p})=\pi(h)-\left(\sigma(h)-h\sigma'(h)\right)\mathcal{E_{\mathrm{cap}}}(\|\boldsymbol{p}\|),\label{P}
\end{equation}
where
\begin{equation}
    {\displaystyle \frac{\pi(h)}{h^{2}}=
        \left(\frac{\Phi(s)}{s}\right)'\Bigr\vert_{s=h}.}\label{pi}
\end{equation}

\smallskip

In all the paper long, we will do the following hypothesis:
\begin{itemize}
    \item $s\mapsto\mathcal{E_{\mathrm{cap}}}(s)$,
          $s\mapsto\Phi(s)$ and $s\mapsto\sigma(s)$ are assumed to be positive,
    \item $\mathcal{E_{\mathrm{cap}}}$ invertible from $\mathbb{R}^{+}$
          to $\mathbb{R}^{+}$ with $\mathcal{E_{\mathrm{cap}}}(0)=0$,
    \item $\pi'(h)>0$ so that $\Phi$ is strictly convex as soon
          as $h>0$.
\end{itemize}
System \eqref{E}--\eqref{pi} is supplemented with the initial data
\begin{equation}
    h\vert_{t=0} = h_0, \qquad h\mathbf{u}\vert_{t=0} = m_0.\label{IC}
\end{equation}
In this context, System \eqref{E}--\eqref{IC} admits
an additional energy conservation law which reads
\begin{equation}
    \partial_{t}\left(\frac{1}{2}h\left\Vert \mathbf{u}\right\Vert ^{2}+E\right)+{\rm div}\left(\left(\frac{1}{2}h\left\Vert \mathbf{u}\right\Vert ^{2}+E+P\right)\mathbf{u}\right)-{\rm div}\left({\rm div}\left(\nabla_{\mathbf{p}}E\right)h\mathbf{u}\right)+{\rm div}\left({\rm div}\left(h\mathbf{u}\right)\nabla_{\mathbf{p}}E\right)=0.\label{eq:Energie_K}
\end{equation}

\bigskip

\noindent {\it Remark.}
For specific choices of the capillary energy, we note that the system
\eqref{eq:Sys_SW_gen} reduces to classical models of the fluid mechanics
literature like
\begin{itemize}
    \item The Euler-Korteweg isothermal system when :
          \[
              E\left(h,\nabla h\right)=\Phi\left(h\right)+\frac{1}{2}\sigma\left(h\right)\|\nabla h\|^{2}
          \]
          where $h$ is the density and $\sigma(h)$ is the capillary coefficient.
    \item The shallow-water type system for thin film flows both:
    \item In the quadratic capillary case
          \[
              E\left(h,\nabla h\right)=\frac{h^2}{2} +\frac{1}{2}\sigma\|\nabla h\|^{2}
              \hbox{ with } h \hbox{ the height of the fluid and }
              \sigma \hbox{ is constant }
          \]

    \item In the fully nonlinear capillary case :
          \[E\left(h,\nabla h\right)=\frac{h^2}{2} +\sigma \sqrt{1- \|\nabla h\|^{2}}
              \hbox{ with } h \hbox{ the height of the fluid and }
              \sigma \hbox{ is constant }
          \]
          Note that the fully nonlinear case admits the following  two asymptotics
          \[
              \begin{array}{ll}

                  {\displaystyle E(h,\nabla h)=
                  {h^2}/{2} +\sigma
                  {\|\nabla h\|^2}/{2}
                  +o_{\|\nabla h\|\to 0}\left(\|\nabla h\|\right)}
              \end{array}
          \]
          and
          \[
              \begin{array}{ll}
                  {\displaystyle E(h,\nabla h)=
                  {h^2}/{2}
                  +\sigma{\|\nabla h\|}+o_{\|\nabla h\|\to\infty}\left(\|\nabla h\|\right)}
              \end{array}
          \]
\end{itemize}

\bigskip

It is a hard problem to propose a discretization of System \eqref{E}--\eqref{IC} that is compatible with the energy
equation \eqref{eq:Energie_K} and this is the objective of our
paper. The main issue is that one cannot
adapt the proof of the energy estimate \eqref{eq:Energie_K} derived
from \eqref{eq:Sys_SW_gen} at a discretized level due to the presence
of high-order derivatives associated to the capillarity energy \eqref{E}. The readers interested in understanding the mathematical
and numerical difficulties are referred to \cite{Po} and important references cited therein.
The strategy first consists in performing a reduction of order in spatial derivatives and in introducing an alternative system (called augmented system) which contains lower order derivatives. It consists secondly in checking that the augmented system may be decomposed into two parts: A conservative hyperbolic part and a second order derivatives part which is skew-symmetric with respect to the $L^2$ scalar product. Such system is really adapted to discretization compatible with the energy: it is obtained by taking $L^2$ scalar products with respect to the new unknowns.
This strategy was applied successfully in the context of Euler-Korteweg isothermal system for numerical purposes when
the internal energy is quadratic with respect to $\nabla h$: see
\cite{NoVi} in the one dimensional case and \cite{BCNV} in the two
dimensional case. In both cases, the augmented version is obtained
by introducing an auxiliary velocity ${\bf v}$ which is proportional
to $\nabla h$ and admits an additional skew-symmetric structure
with respect to the $L^{2}$ scalar product which makes the proof
of energy estimates and the design of compatible numerical scheme
easier. However, this approach was not extended to more general
internal energy \eqref{E}. This is the objective of the paper to define the appropriate unknowns in order to get an appropriate augmented system for numerical purposes.

Note that there exists several interesting papers developing augmented systems such as \cite{GaGo1} and \cite{Go} for symmetric form for capillarity fluids with a capillarity energy $E(h,\nabla h)$ or multi-gradient fluids with a capillarity energy $E(h,\nabla h, \cdots,\nabla^n h)$. See also recently \cite{DaFaGa} for the defocusing Schr\"odinger equation which is linked to the quantum-Euler system ($E(h,\mathbf{p}) = \Phi(h) + \sigma \|\mathbf{p}\|^2/h$ where $\sigma$ is constant) through the Madelung transform and some numerical simulation.

It is interesting to note that the augmented system in \cite{GaGo1} and \cite{Go} is related to the unknowns ($h$, $\mathbf{u}$, $\nabla h$, $\cdots$, $\nabla^n h$). In \cite{La}, the authors developed a similar augmented version in order to deal with internal capillarity energies \eqref{E} for numerical purposes namely:
\begin{equation}
    \left\{ \begin{array}{cc}
        \partial_{t}{\displaystyle h+{\rm div}\left(h\mathbf{u}\right)=0}                                                                                                                                                                                           & (i)   \\
        \partial_{t}{\displaystyle \left(h\mathbf{u}\right)+{\rm div}\left(h\mathbf{u}\otimes\mathbf{u}\right)+\nabla P+{\rm div}\left(\boldsymbol{p}\otimes\nabla_{\mathbf{p}}E_{tot}\right)=\nabla\left(h{\rm div}\left(\nabla_{\mathbf{p}}E_{tot}\right)\right)} & (ii)  \\
        \partial_{t}{\displaystyle \mathbf{p}+\nabla\left(\mathbf{p}^t\mathbf{u}\right)=-\nabla\left(h\,{\rm div}\left(\mathbf{u}\right)\right)}                                                                                                                    & (iii)
    \end{array}\right.\label{eq:Form2D_Lall}
\end{equation}
where $E_{tot} = h |\mathbf{u}|^2/2 + E(h,\mathbf{p})$.
However, in the $2$-dimensional setting, the assumption ${\rm curl}\,\mathbf{p}=0$ has to be made to show the conservation of the total energy and therefore it has to be satisfied initially:
The interested reader is referred pages 166--168. This constraint is
not satisfied at the discretized level and it creates instabilities.

\noindent {\it Remark.} In order to avoid such a constraint which is hardly
guaranteed in the discrete case, one could use instead the following
modified formulation
\[
    \left\{ \begin{array}{cc}
        \partial_{t}{\displaystyle h+{\rm div}\left(h\mathbf{u}\right)=0}                                                                                                                                                                                                                                                                                                  & (i)   \\
        \partial_{t}{\displaystyle \left(h\mathbf{u}\right)+{\rm div}\left(h\mathbf{u}\otimes\mathbf{u}\right)+\nabla P+{\rm div}\left(\boldsymbol{p}\otimes\nabla_{\mathbf{p}}E_{tot}\right)-\left(\left(\nabla\mathbf{p}\right)^{t}-\left(\nabla\mathbf{p}\right)\right)\nabla_{\mathbf{p}}E_{tot}=\nabla\left(h{\rm div}\left(\nabla_{\mathbf{p}}E_{tot}\right)\right)} & (ii)  \\
        \partial_{t}{\displaystyle \mathbf{p}+\nabla\left(\mathbf{p}^t\mathbf{u}\right)=-\nabla\left(h{\rm div}\left(\mathbf{u}\right)\right)}                                                                                                                                                                                                                             & (iii)
    \end{array}\right.
\]
for which it is easy to prove the conservation of the total energy
\[
    \begin{array}{c}
        \partial_{t}\left(E_{tot}\right)+div\left(\boldsymbol{u}\left(E_{tot}+\pi\right)\right)=\left({\rm div}\left(h(\boldsymbol{u}^{t}\nabla)(\nabla_{\mathbf{p}}E_{tot})\right)-{\rm div}(h(\nabla_{\mathbf{p}}E_{tot}^{t}\nabla){\bf u})\right) \\
        -{\rm div}\left(\boldsymbol{u}(\boldsymbol{p}^{t}\nabla_{\mathbf{p}}E_{tot}-\left(\sigma-h\sigma'\right)\mathcal{E_{\mathrm{cap}}})\right)
    \end{array}
\]
for any smooth solution of the above system without assuming the curl free
assumption on $\mathbf{p}$. However, this formulation introduces non-conservative terms in the left-hand side of the momentum equation
and it is then hard to satisfy for conservation of momentum and energy at the discrete level.

\bigskip

In our paper, defining an appropriate velocity field $\mathbf{v}$ instead
of $\mathbf{p}$, we are able to design an appropriate augmented version which may be decomposed as the sum of a conservative hyperbolic part and
a skew-symmetric second order differential operator for the $L^2$ scalar product. The system is solved in the variables $(h,\mathbf{u},\mathbf{v})$
and if regularity occurs we recover the expression of $\mathbf{v}$ in terms of $h$, $\mathbf{p}$ and $\|\mathbf{p}\|$. The important property is
that the energy conservation law may be satisfied easily at the
discretized level using the particular structure of our augmented
system. The particular form allows also an efficient splitting method
allowing to simulate complex situations like large gradients of fluid height.

In the small gradient limit, this fomulation is equivalent to the
one derived by D. Bresch, F. Couderc, P. Noble and J.--P. Vila in
\cite{BCNV}. Our formulation is valid for any internal energy in
the form $E(h,{\bf p})=\Phi(h)+\sigma(h)\mathcal{E}_{cap}(\|{\bf p}\|)$.
When specified to $\mathcal{E}_{cap}(q)=\sqrt{1+q^{2}}-1$, we see
that in the high gradient limit, $\mathcal{E}_{cap}(q)\sim_{q\to\infty}|q|$
which is a capillary term found usually in two fluids systems. We
thus expect our approach to be useful in the context of bi-fluid flows.
Note also that our paper could be also of practical interest to deal
with generalization of Euler-Korteweg system: see \cite{Ko} and \cite{Jo}
for discussions on compressible Korteweg type systems.

We rely on the new augmented system to propose a numerical scheme
which is energetically stable and extends what was done in \cite{BCNV}
and \cite{NoVi}. Note that skew-symmetric augmented versions of the
capillary shallow water equations in the $L^2$ scalar product are also useful from a theoretical point of view: see e.g. \cite{BrGiLa} for the proof of existence of dissipative solutions to the Euler-Korteweg isothermal system. Our present work will be the starting point to improve the work by
Lallement and Villedieu (see \cite{La} and \cite{LaTrLaVi}) related
to disjunction term for triple point simulations: see \cite{BrCeCoGiLaNoRiRuViVi}.

\bigskip

The paper is divided in three parts: The first part introduces the
augmented version with full surface tension and discuss its connection
with the system derived in \cite{BCNV}. In the second part, we propose
a numerical scheme satisfying energy stability. Finally, we present
numerical illustrations based on our numerical scheme showing the
importance of considering our augmented system with the full surface
tension.

\bigskip

\section{Augmented version}

Extending ideas from \cite{NoVi} in the one dimensional case, an
augmented formulation of the shallow water equations \eqref{eq:Sys_SW_gen}
with $\mathcal{E_{\mathrm{cap}}}\left(\|\nabla h\|\right)=\frac{\sigma}{2}\left\Vert \nabla h\right\Vert ^{2}$
was proposed in \cite{BCNV} in the two dimensional setting: it is a second order system of PDEs which may be decomposed in two parts: A conservative hyperbolic part and a second order derivatives part which is skew symmetric with respect to the $L^{2}$ scalar product.
The additional quantity in \cite{BCNV} was given by ${\bf w}=\nabla\phi(h)$ with $\phi'(h)=\sqrt{{\sigma(h)}/{h}}$: it is thus colinear to $\nabla h$.

The main objective here is to consider a more general internal capillarity energy namely \eqref{E}. We now introduce our new formulation of \eqref{eq:Sys_SW_gen} which is valid in the fully decoupled case and provides a dual formulation of capillary terms which ensures a straightforward consistent energy
balance. To this end we introduce an additional unknown, denoted $\boldsymbol{v}$,
which is colinear to $\nabla h$ and satisfies
\[
    \frac{1}{2}h\left\Vert \boldsymbol{v}\right\Vert ^{2}=\sigma\left(h\right)\mathcal{E_{\mathrm{cap}}}\left(\|\nabla h\|\right)
\]
where $q=\|\boldsymbol{p}\|=\|\nabla h\|$. To do so, we define \textbf{v}
as
\begin{equation}
    \boldsymbol{v}=\alpha(q^{2}) \,\sqrt{\dfrac{\sigma(h)}{h}}\,  \boldsymbol{p}  \label{v}
\end{equation}
where the function $\alpha: s\mapsto \alpha(s)$ is given by
\[
    \alpha(s)=\sqrt{\frac{2\mathcal{E_{\mathrm{cap}}}
            \left(\sqrt{s}\right)}{s}}.
\]

\noindent {\bf Remark.}
Note that using the definition $\boldsymbol{v}$, we have the following
relations
\[
    \left\Vert \boldsymbol{v}\right\Vert ^{2}=\alpha^{2}(\left\Vert \boldsymbol{p}\right\Vert ^{2})\left\Vert \boldsymbol{p}\right\Vert ^{2}\dfrac{\sigma}{h},\]
\[   \frac{1}{2}\alpha^{2}(q^{2}) \, q^{2} \, \sigma(h)=\sigma(h)\mathcal{E_{\mathrm{cap}}}\left(q\right).
\]

\noindent {\bf Remark.}
Note that $\boldsymbol{v}$ has the dimension of a
velocity and transforms the capillary energy into some kinetic energy.
This interpretation of the capillary energy in terms of kinetic energy
in our augmented system defined below motivates surely the robustness
of our results.

\bigskip

Let us now write a system related to the unknowns $(h,\boldsymbol{u},\boldsymbol{v})$
where $\boldsymbol{v}$ is given by \eqref{v} with $\boldsymbol{p}=\nabla h$.
This will provide a system which combines a first order conservative and hyperbolic part on $(h,\boldsymbol{u},\boldsymbol{v})$ together with a second order part which has a skew-symmetric structure (for the $L^{2}$ scalar product). % for the full
%surface tension term extending the calculations made in \cite{BCNV}
%by D. Bresch, F. Couderc, P. Noble, J.--P. Vila. Such augmented version
%will be the starting point for numerical schemes satisfying energy
%stability and numerical computations showing the importance to consider
%the full surface tension compared to usual simplified versions.
More precisely, we have the following result.

\begin{lem} \label{NLaugmented}
    i) Let
    \begin{equation}
        \boldsymbol{U}=\left(\begin{array}{c}
                h               \\
                h\boldsymbol{u} \\
                h\boldsymbol{v}
            \end{array}\right),\qquad
        \mathcal{F} \left( \boldsymbol{ U}\right)=\left(\begin{array}{c}
                h\boldsymbol{u}                \\
                h\boldsymbol{u}\otimes\boldsymbol{u}
                +\pi\left(h\right) {\rm I_{d}} \\
                h\boldsymbol{v}\otimes\boldsymbol{u}
            \end{array}\right)\label{def_U_F}
    \end{equation}
    where ${\rm I_{d}}$ is the $d\times d$ identity matrix and
    $$\mathbf{M} = \mathcal{M}(h,\mathbf{v})(\mathbf{U})$$
    with, for all $U_1=(h_1,h_1\mathbf{u_1}, h_1 \mathbf{v_1})^t$,
    \begin{equation}
        \mathcal{M}(h,\mathbf{v})(\mathbf{U_1})=\left(\begin{array}{l}
                0                                                                                                                                                \\
                \operatorname{div}\left(h\nabla( f (h,\boldsymbol{v})\boldsymbol{v}_1)^{t}\right)-\nabla( \boldsymbol {g}(h,\boldsymbol{v})^{t}\boldsymbol{v_1}) \\
                - f (h,\boldsymbol{v})\operatorname{div}\left(h\nabla\boldsymbol{u}_1^{t}\right)\quad- \boldsymbol {g} (h,\boldsymbol{v})\operatorname{div}\boldsymbol{u}_1
            \end{array}\right)\label{eq:def_M}
    \end{equation}
    where $f(h,\boldsymbol{v})$ is a symmetric tensor and $\boldsymbol {g} \left(h,\boldsymbol{v}\right)$
    a vector field given by
    \[
    {f}(h,\boldsymbol{v})=\sqrt{\sigma(h)}\sqrt{h}\left(2\frac{\alpha'(q^{2})h}{\alpha(q^{2})^{2}\sigma(h)}\boldsymbol{v}\otimes\boldsymbol{v}+\alpha(q^{2}){\rm I_{d}}\right)
    \]
    \[
    \boldsymbol {g}\left(h,\boldsymbol{v}\right)=\left(\left(\frac{\sigma'(h)h}{2\sigma(h)}+\frac{1}{2}\right)\,+2\frac{\alpha'(q^{2})}{\alpha(q^{2})}q^{2}\right)h\boldsymbol{v}
    \]
    and
    \[
    \alpha(q^{2})=\frac{\sqrt{2\mathcal{E_{\mathrm{cap}}}\left(q\right)}}{q}\qquad\hbox{with}\qquad q={\mathcal{E}}_{{\rm cap}}^{-1}\left(\frac{h\|\boldsymbol{v}\|^{2}}{2\sigma(h)}\right).
    \]
    The augmented system reads
    \begin{equation}
        \partial_{t}  \boldsymbol{ U} +{\rm div}\left(\mathcal{F}
        \left(\boldsymbol{ U}\right)\right)=\boldsymbol{M}.\label{eq:Augmented_FHS-1}
    \end{equation}

    \smallskip{}

    \noindent ii) If $(h,\boldsymbol{u},\boldsymbol{v})$ is regular enough
    then it also satisfies the following energy balance
    \begin{eqnarray}
        \partial_{t}\left(\frac{1}{2}h\left\Vert \boldsymbol{u}\right\Vert ^{2}+E\right) &  & +\hskip1cm{\rm div}\left(\boldsymbol{u}\,(\frac{1}{2}h\left\Vert \boldsymbol{u}\right\Vert ^{2}+E+\pi)\right)\\
        &  & =\left(\operatorname{div}\left(h\boldsymbol{u}^{t}\nabla^{t}(  f(h,\boldsymbol{v})^{t}\boldsymbol{v})\right)-
        \operatorname{div}(h \nabla\boldsymbol{u}   f(h,\boldsymbol{v})\boldsymbol{v})\right)-\operatorname{div}\left(\boldsymbol{u}( f(h,\boldsymbol{v})^{t}\boldsymbol{v})\right)\nonumber
    \end{eqnarray}
    where $E= \Phi(h) + h\|\boldsymbol{v}\|^2/2.$

    \smallskip

    \noindent    iii) If $(h,\boldsymbol{u})$ is regular enough and the initial velocity
    $\boldsymbol{v}_{0}$ satisfies
    \[
        \boldsymbol{v}_{0}=\alpha(\left\Vert \nabla h_{0}\right\Vert ^{2})\sqrt{\frac{\sigma(h_{0})}{h_{0}}}\nabla h_{0}
    \]
    then $\boldsymbol{v}$ satisfies also
    \[
        \boldsymbol{v}=\alpha(\left\Vert \nabla h\right\Vert ^{2})\sqrt{\frac{\sigma(h)}{h}}\nabla h
    \]
    and $(h,\boldsymbol{u})$ solves the original equations with the full
    surface tension term given by \eqref{eq:Sys_SW_gen}--\eqref{P}. \end{lem}
\section{Energetically stable numerical scheme}
The augmented formulation \eqref{eq:Augmented_FHS-1} in Lemma \ref{NLaugmented} reads
\[
    \partial_{t}\boldsymbol{U}+{\rm div}\left(\mathcal{F} \left(\boldsymbol{U}\right)\right)=\boldsymbol{M}
\]
with definitions \eqref{def_U_F} and \eqref{eq:def_M} of $\boldsymbol{U}, \mathcal{F}$ and $\boldsymbol{M}$.
The first order part of the augmented formulation in the left-hand side
is the classical Euler barotropic model with an additional transport
equation. It admits an additional conservation law related to the total energy:
\[
    E_{tot}=\frac{\left\Vert h\boldsymbol{u}\right\Vert ^{2}}{2h}+\Phi\left(h\right)+\frac{\left\Vert h\boldsymbol{v}\right\Vert ^{2}}{2h}.
\]
whereas the entropy variable is
\[
    \left(\nabla_{\boldsymbol{U}}E_{tot}\right)^{t}=\boldsymbol{V}^t=\left(-\frac{1}{2}\left(\left\Vert \boldsymbol{u}\right\Vert ^{2}+\left\Vert \boldsymbol{v}\right\Vert ^{2}\right)+\Phi'\left(h\right),\boldsymbol{u}^{t},\boldsymbol{v}^{t}\right).
\]
This total energy is the total energy of the shallow-water equations
with surface tension whereas the potential energy associated to surface
tension is transformed into kinetic energy associated to the artificial
velocity ${\bf v}$. The full system admits also an energy equation:
{\setlength{\arraycolsep}{1pt}
\begin{eqnarray*}
    {\displaystyle \partial_{t}E_{tot}+\mbox{div}\left(\boldsymbol{u}\left(E_{tot}+\pi\left(h\right)\right)\right)} & = & \boldsymbol{V}^{t}\boldsymbol{M}\\
    {\displaystyle } & = & \mbox{div} \left(h\boldsymbol{u}^{t}\nabla(\nabla_{\mathbf{p}}E)\right)
    - \mbox{div}\left(h(\nabla_{\mathbf{p}}E)^{t}\nabla{\bf u}\right)\\
    {\displaystyle } &  & - \mbox{div} \left(\boldsymbol{u}\left(\boldsymbol{p}^{t}\nabla_{\mathbf{p}}E-\left(\sigma-h\sigma'\right)\mathcal{E_{\mathrm{cap}}}\left(q\right)\right)\right)
\end{eqnarray*}
} with the right-hand side in conservation form. One of the aim of
this paper is to design a numerical scheme that is free from a CFL
condition associated to surface tension. For that purpose, we follow
the strategy in \cite{BCNV} and introduce an IMplicit-EXplicit strategy
where the hyperbolic step is explicit in time whereas the step associated
to surface tension is implicit in time. The spatial discretization
is based on an entropy dissipative scheme for the first order part
whereas we mimic the skew symmetric structure found at the continuous
level to discretize the right hand side $\mathcal{M}$. We prove that
this strategy is energetically stable in the case of periodic boundary
conditions.

\subsection{IMplicit - EXplicit formulation}

\begin{figure}[!ht]
    \begin{centering}
        \includegraphics[width=7cm,height=4cm]{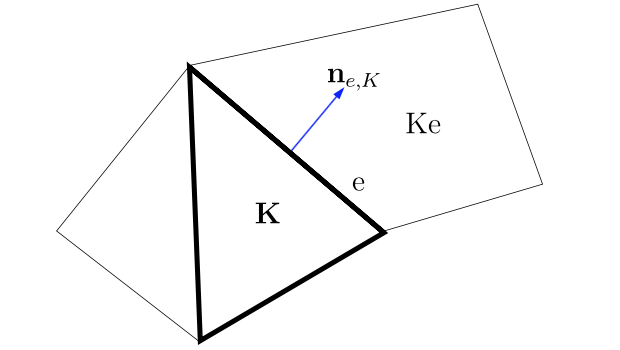}
        \par\end{centering}
    \caption{Notations for cell $\mathbf{K}$}
    \label{fig-mesh}
\end{figure}

Following \cite{BCNV}, we consider the following IMplicit-EXplicit
time discretization: the hyperbolic step is explicit
\begin{equation}
    \frac{\boldsymbol{{U}}^{n+1/2}-\boldsymbol{U}^{n}}{\Delta t}+{\rm div}\left(\mathcal{F}\left(\boldsymbol{U}^{n}\right)\right)=0
    \label{eq:Hyp_step}
\end{equation}
and the capillary skew symmetric second order step
\begin{equation}
    \frac{\boldsymbol{U}^{n+1}-\boldsymbol{U}^{n+1/2}}{\Delta t}=\boldsymbol{M}^{n+1}\label{eq:Cap_step}
\end{equation}
with
\[
    {\displaystyle \boldsymbol{M}^{n+1}=\left(\begin{array}{l}
                    0                                                                               \\
                    \operatorname{div}\left(h^{n+1}\nabla(f(h^{n+1},\boldsymbol{v}^{n+1/2})\boldsymbol{v}^{n+1})^{t}\right)
                    -\nabla(\boldsymbol{g}(h^{n+1},\boldsymbol{v}^{n+1/2})^{t}\boldsymbol{v}^{n+1}) \\
                    -f(h^{n+1},\boldsymbol{v}^{n+1/2})\operatorname{div}\left(h^{n+1}\nabla(\boldsymbol{u}^{n+1})^{t}\right)\quad
                    -\boldsymbol{g}(h^{n+1},\boldsymbol{v}^{n+1/2})\operatorname{div}\boldsymbol{u}^{n+1}
                \end{array}\right).}
\]
The second step is not fully implicit: instead it is semi-implicit
so that the problem to solve for $({\bf v}^{n+1},{\bf u}^{n+1})$
is \textit{linear}. This does not affect the order of the time discretization
since the time splitting is already first order in time. Let us now
consider the spatial discretization. We will use a generic Finite
Volume context. %%%%%
We introduce a spatial discretization of $\nabla$ and ${\rm div}$
operators with finite volume methods. For that purpose, we denote
$\mathbf{K}$ a cell of the mesh $T_{d}$, $\mathbf{e}\in\partial\mathbf{K}$
an edge of $\mathbf{K}$ and $\mathbf{K_{e}}$ a neighboring cell
of $\mathbf{K}$: see figure \ref{fig-mesh} for an illustration.
\noindent We use a classical entropy satisfying scheme of numerical
flux
\[
    \mathbf{G}_{e,K}^{n}=
    \mathbf{G}\left(\mathbf{U}_{K}^{n},\mathbf{U}_{K_{e}}^{n},n_{e,K}\right)
\]
where $n_{e,K}$ is the outward normal to the cell $K$ (of measure
$m_{K}$) at the edge $e$ (of measure $m_{e}$). We denote $\mathbf{U}_{K}$
the average of the vector $\mathbf{U}$ on the cell \textbf{K}. The hyperbolic
step then reads
\begin{equation}
    \mathbf{U}_{K}^{n+1/2}= \mathbf{U}_{K}^{n}
    -\frac{\Delta t}{m_{K}}\sum_{e\in\partial K}m_{e} \mathbf{G}_{e,K}^{n}\label{eq:Hyp_part}
\end{equation}
and we assume that it is entropy dissipative in the sense that it
satisfies the following discrete Entropy inequality
\begin{equation}
    E_{tot}\left(\mathbf{U}_{K}^{n+1/2}\right)\leq E_{tot}\left(\mathbf{U}_{K}^{n}\right)-\frac{\Delta t}{m_{K}}\sum_{e\in\partial K}m_{e}H_{e,K}^{n}\label{eq:Entropy_Hyp}
\end{equation}
where $H_{e,K}^{n}$ is the entropy numerical flux associated with
$G_{e,K}^{n}.$ In the particular case of Euler Barotropic equations
such numerical schemes exist and satisfy this inequality provided
a hyperbolic CFL condition of the type
\begin{equation}
    \max_{K}\frac{\Delta t}{m_{K}}m_{e}\left\Vert \nabla_{\mathbf{U}} F\left(\mathbf{U}_{K}^{n}\right)\right\Vert <a<1\label{eq:CFL_Hyp}
\end{equation}
is satisfied for some $a>0$. Moreover, under a similar CFL condition,
the positivity of the fluid $h$ is preserved and the total energy
$E_{tot}\left(\mathbf{U}\right)$ is strictly convex: this will be a useful
property to prove entropy stability for numerical schemes. The second
step is
\begin{equation}
    \mathbf{U}_{K}^{n+1}= \mathbf{U}_{K}^{n+1/2}+\Delta t \, \mathbf{M}_{K}^{n+1}\label{eq:Sys_Cap_dis}
\end{equation}
with
\begin{equation}
    \mathbf{M}_{K}^{n+1}=\left(\begin{array}{c}
            0                                                                                                                                                                                                                                                                   \\
            -\nabla_{3,\Delta}\left(\mathbf{g}(h_{K}^{{n+1}},\boldsymbol{v}_{K}^{n+1/2})^{t}\boldsymbol{v}_{K}^{n+1}\right)+{\rm div}_{1,\Delta}\left(h_{K}^{{n+1}}\nabla_{1,\Delta}\left(f(h_{K}^{{n+1}},\boldsymbol{v}_{K}^{n+1/2})\boldsymbol{v}_{K}^{n+1}\right)^{T}\right) \\[3mm]
            -\mathbf{g}(h_{K}^{{n+1}},\boldsymbol{v}_{K}^{n+1/2}){\rm div}_{3,\Delta}\left(\boldsymbol{u}_{K}^{n+1}\right)-f(h_{K}^{{n+1}},\boldsymbol{v}_{K}^{n+1/2}){\rm div}_{1,\Delta}\left(h_{K}^{{n+1}}\nabla_{1,\Delta}\left(\boldsymbol{u}_{K}^{n+1}\right)^{T}\right)
        \end{array}\right)\label{eq:def_M_disc}
\end{equation}
where $\nabla_{3,\Delta}$, ${\rm div}_{1,\Delta}$, $\nabla_{1,\Delta}^{T}$,
${\rm div}_{3,\Delta}$ are linear discrete operators approximating
the corresponding ones in the definition of the operator $\mathcal{M}$
and that will be defined hereafter. In particular ${\rm div}_{3,\Delta}$
shall be chosen as the dual discrete operator of $\nabla_{3,\Delta}$
in the following sense :
\begin{equation}
    (a,\nabla_{3,\Delta}\left(\varphi\right))_{T_{d}}=-\left({\rm div}_{3,\Delta}\left(a\right),\varphi\right)_{T_{d}}\label{eq:defGD3}
\end{equation}
for any smooth function $\varphi$ and $a$ defined on the mesh $T_{d}$
where we have used the discrete scalar product below
\[
    (a,b)_{T_{d}}=\sum_{K\in T_{d}}m_{k}\langle a_{K},b_{K}\rangle_{\mathbb{R}^{d}}.
\]
One possible choice is taking the classical approximation of flux
in the finite volume context which leads to
\[
    {\rm div}_{3,\Delta}\left(\boldsymbol{a}\right)=\frac{1}{m_{K}}\sum_{e\in\partial K}m_{e}\frac{1}{2}\left(\boldsymbol{a}_{K_{e}}+\boldsymbol{a}_{K}\right).n_{e,K}=\frac{1}{2m_{K}}\sum_{e\in\partial K}m_{e}\boldsymbol{a}_{K_{e}}.n_{e,K}
\]
and the corresponding (weak) approximation of $\nabla_{3,\Delta}\left(\varphi\right)$
\[
    \nabla_{3,\Delta}\left(\varphi\right)=\frac{1}{2m_{K}}\sum_{e\in\partial K}m_{e}\frac{1}{2}\left(\varphi_{K}\frac{m_{K_{e}}}{m_{K}}-\varphi_{K_{e}}\right)n_{e,K}.
\]
In the context of finite difference approximations, we may consider
the discrete analogous of the ${\rm div}$ operator
\begin{equation}
    div_{3,\Delta}\left(\boldsymbol{a}\right)_{ij}=\frac{a_{i+1,j}^{x}-a_{i-1,j}^{x}}{2\Delta x}+\frac{a_{i,j+1}^{y}-a_{i,j-1}^{y}}{2\Delta y}\label{eq:deg_div3}
\end{equation}
which leads to
\begin{equation}
    \nabla_{3,\Delta}\left(\varphi\right)_{ij}=\frac{\varphi_{i+1,j}-\varphi_{i-1,j}}{2\Delta x}\boldsymbol{n}_{x}+\frac{\varphi_{i,j+1}-\varphi_{i,j-1}}{2\Delta y}\boldsymbol{n}_{y}.\label{eq:def_G3}
\end{equation}

\begin{rem} In the case of general finite volume discretization on
    any mesh, the question of finding consistent second order operators
    is not so simple and requires some refined tools such as renormalisation
    or adhoc discrete gradient (see eg {\rm \cite{ADEM,LVII,JB}}). \end{rem}

In the next section, we focus on the definition of the discrete divergence
and gradients operators ${\rm div}_{1,\Delta}$ and $\nabla_{1,\Delta}$
so as to ensure the energy stability.

\subsection{Energy Stability of first order schemes}
Let us now analyse the stability properties of the above scheme. The
hyperbolic step is entropy stable in the sense that
\[
    \sum_{K}E_{tot}\left(\mathbf{U}_{K}^{n+1/2}\right)m_{K}\leq\sum_{K}E_{tot}\left(\mathbf{U}_{K}^{n}\right)m_{K}.
\]
since it is a direct consequence of entropy inequality (\ref{eq:Entropy_Hyp}).
Let us now focus on the ``capillary time step'' and the definition
of ${\rm div}_{1,\Delta}$ and $\nabla_{1,\Delta}$. In order to get
more compact form of discrete operators, let us define
\begin{equation}
    \begin{array}{c}
        \left(\partial_{x,1\Delta}\left(m\right)\right)_{i+1/2,j}=\frac{m_{i+1,j}-m_{i,j}}{\Delta x},\;\left(\partial_{x,1\Delta}^{0}\left(p\right)\right)_{i,j}=\frac{p_{i+1/2,j}-p_{i-1/2,j}}{\Delta x},\;\left(\partial_{x,1\Delta}^{00}\left(m\right)\right)_{i,j}=\frac{m_{i+1,j}-m_{i-1,j}}{2\Delta x}, \\
        \left(\partial_{y,1\Delta}\left(m\right)\right)_{i,j+1/2}=\frac{m_{i,j+1}-m_{i,j}}{\Delta y},\;\left(\partial_{y,1\Delta}^{0}\left(p\right)\right)_{i,j}=\frac{p_{i,j+1/2}-p_{i,j-1/2}}{\Delta y},\;\left(\partial_{y,1\Delta}^{00}\left(m\right)\right)_{i,j}=\frac{m_{i,j+1}-m_{i,j-1}}{2\Delta y}
    \end{array}\label{eq:def_der_disc_DF}
\end{equation}

\begin{equation}
    {\rm div}_{1,\Delta}\left(h\nabla_{1,\Delta}\boldsymbol{m}^{T}\right)=\left(\begin{array}{c}
            \left(\partial_{x,1\Delta}^{0}\left(h\partial_{x,1\Delta}m^{x}\right)\right)+\left(\partial_{y,1\Delta}^{00}\left(h\partial_{x,1\Delta}^{00}m^{y}\right)\right) \\
            \left(\partial_{x,1\Delta}^{00}\left(h\partial_{y,1\Delta}^{00}m^{x}\right)\right)+\left(\partial_{y,1\Delta}^{0}\left(h\partial_{y,1\Delta}m^{y}\right)\right)
        \end{array}\right)\label{eq:def_d1D_compact}
\end{equation}
We thus have the following property : \begin{lem} \label{Dual_discret_compact}
    Let us suppose that ${\rm div}_{1,\Delta}\left(h\nabla_{1,\Delta}\boldsymbol{m}^{T}\right)$
    is defined as (\ref{eq:def_d1D_compact}): Then we have
    \begin{equation}
        \left(\boldsymbol{u},{\rm div}_{1,\Delta}\left(h\nabla_{1,\Delta}\boldsymbol{m}^{T}\right)\right)_{T_{d}}=\left(\boldsymbol{m},div_{1,\Delta}\left(h\nabla_{1,\Delta}\boldsymbol{u}^{T}\right)\right)_{T_{d}}\label{eq:ident_dual_cap}
    \end{equation}
    where
    \[
        \left(\boldsymbol{a},\boldsymbol{b}\right)_{T_{d}}=\sum_{i,j}\Delta y\Delta x\langle\boldsymbol{a}_{ij},\boldsymbol{b}_{ij}\rangle_{\mathbb{R}^{d}}.
    \] \end{lem}

\noindent \textbf{Proof of Lemma \ref{Dual_discret_compact}.} Thanks
to definitions (\ref{eq:def_der_disc_DF})--(\ref{eq:def_d1D_compact})
we have
\begin{eqnarray*}
    \left(\partial_{x,1\Delta}^{0}\left(h\partial_{x,1\Delta}m^{x}\right)\right)_{i,j} &  & =\frac{1}{\left(\Delta x\right)^{2}}\left(h_{i+1/2,j}\left(m_{i+1,j}^{x}-m_{i,j}^{x}\right)-h_{i-1/2,j}\left(m_{i,j}^{x}-m_{i-1,j}^{x}\right)\right)\\
    \left(\partial_{y,1\Delta}^{00}\left(h\partial_{x,1\Delta}^{00}m^{y}\right)\right)_{i,j} &  & =\frac{1}{4\Delta y\Delta x}\left(h_{i,j+1}\left(m_{i+1,j+1}^{y}-m_{i-1,j+1}^{y}\right)-h_{i,j-1}\left(m_{i+1,j-1}^{y}-m_{i-1,j-1}^{y}\right)\right)\\
    \left(\partial_{x,1\Delta}^{00}\left(h\partial_{y,1\Delta}^{00}m^{x}\right)\right)_{i,j} &  & =\frac{1}{4\Delta y\Delta x}\left(h_{i+1,j}\left(m_{i+1,j+1}^{x}-m_{i+1,j-1}^{x}\right)-h_{i-1,j}\left(m_{i-1,j+1}^{x}-m_{i-1,j-1}^{x}\right)\right)\\
    \left(\partial_{y,1\Delta}^{0}\left(h\partial_{y,1\Delta}m^{y}\right)\right)_{i,j} &  & =\frac{1}{\left(\Delta y\right)^{2}}\left(h_{i,j+1/2}\left(m_{i,j+1}^{y}-m_{i,j}^{y}\right)-h_{i,j-1/2}\left(m_{i,j}^{y}-m_{i,j-1}^{y}\right)\right)
\end{eqnarray*}
where we take
\[
    h_{i+1/2,j}=\frac{1}{2}\left(h_{i+1,j}+h_{i,j}\right), \qquad h_{i,j+1/2}=\frac{1}{2}\left(h_{i,j+1}+h_{i,j}\right).
\]
It follows
\begin{eqnarray*}
    &  & \left(\boldsymbol{u},{\rm div}\left(h\nabla\boldsymbol{m}^{T}\right)\right)_{T_{d}}=\sum\Delta y\Delta xu_{i,j}^{x}\frac{1}{\left(\Delta x\right)^{2}}\left(h_{i+1/2,j}\left(m_{i+1,j}^{x}-m_{i,j}^{x}\right)-h_{i-1/2,j}\left(m_{i,j}^{x}-m_{i-1,j}^{x}\right)\right)\\
    &  & +\sum\Delta y\Delta x\left(u_{i,j}^{x}\frac{1}{4\Delta y\Delta x}\left(h_{i,j+1}\left(m_{i+1,j+1}^{y}-m_{i-1,j+1}^{y}\right)-h_{i,j-1}\left(m_{i+1,j-1}^{y}-m_{i-1,j-1}^{y}\right)\right)\right)\\
    &  & +\sum\Delta y\Delta x\frac{1}{4\Delta y\Delta x}u_{i,j}^{y}\left(h_{i+1,j}\left(m_{i+1,j+1}^{x}-m_{i+1,j-1}^{x}\right)-h_{i-1,j}\left(m_{i-1,j+1}^{x}-m_{i-1,j-1}^{x}\right)\right)\\
    &  & +\sum\Delta y\Delta x\frac{1}{\left(\Delta y\right)^{2}}u_{i,j}^{y}\left(h_{i,j+1/2}\left(m_{i,j+1}^{y}-m_{i,j}^{y}\right)-h_{i,j-1/2}\left(m_{i,j}^{y}-m_{i,j-1}^{y}\right)\right).
\end{eqnarray*}
We compute successively
\begin{eqnarray*}
    &  & \sum\frac{\Delta y\Delta x}{\left(\Delta x\right)^{2}}u_{i,j}^{x}\left(h_{i+1/2,j}\left(m_{i+1,j}^{x}-m_{i,j}^{x}\right)-h_{i-1/2,j}\left(m_{i,j}^{x}-m_{i-1,j}^{x}\right)\right)\\
    &  & =\sum\frac{\Delta y\Delta x}{\left(\Delta x\right)^{2}}\left(\left(u_{i-1,j}^{x}h_{i-1/2,j}m_{i,j}^{x}-u_{i,j}^{x}h_{i+1/2,j}m_{i,j}^{x}\right)-\left(u_{i,j}^{x}h_{i-1/2,j}m_{i,j}^{x}-u_{i+1,j}^{x}h_{i+1/2,j}m_{i,j}^{x}\right)\right)\\
    &  & =\sum\frac{\Delta y\Delta x}{\left(\Delta x\right)^{2}} \, m_{i,j}^{x}\left(\left(u_{i+1,j}^{x}-u_{i,j}^{x}\right)h_{i+1/2,j}-\left(u_{i,j}^{x}-u_{i-1,j}^{x}\right)h_{i-1/2,j}\right)\\
    &  & =\sum\frac{\Delta y\Delta x}{\left(\Delta x\right)^{2}} \, m_{i,j}^{x}\left(h_{i+1/2,j}\left(\partial_{x,1\Delta}u^{x}\right)_{i+1/2,j}-h_{i-1/2,j}\left(\partial_{x,1\Delta}u^{x}\right)_{i-1/2,j}\right)\\
    &  & =\sum\Delta y\Delta x \, m_{i,j}^{x}\partial_{x,1\Delta}^{0}\left(h\left(\partial_{x,1\Delta}u^{x}\right)\right)_{i,j}
\end{eqnarray*}
and
\begin{eqnarray*}
    &  & \sum\frac{\Delta y\Delta x}{4\Delta y\Delta x}u_{i,j}^{x}\left(h_{i,j+1}\left(m_{i+1,j+1}^{y}-m_{i-1,j+1}^{y}\right)-h_{i,j-1}\left(m_{i+1,j-1}^{y}-m_{i-1,j-1}^{y}\right)\right)\\
    &  & =\sum\frac{\Delta y\Delta x}{4\Delta y\Delta x}\left(u_{i,j}^{x}h_{i,j+1}m_{i+1,j+1}^{y}-u_{i,j}^{x}h_{i,j+1}m_{i-1,j+1}^{y}-u_{i,j}^{x}h_{i,j-1}m_{i+1,j-1}^{y}+u_{i,j}^{x}h_{i,j-1}m_{i-1,j-1}^{y}\right)\\
    &  & =\sum\frac{\Delta y\Delta x}{4\Delta y\Delta x} \, m_{i,j}^{y}\left(u_{i-1,j-1}^{x}h_{i-1,j}-u_{i+1,j-1}^{x}h_{i+1,j}-u_{i-1,j+1}^{x}h_{i-1,j}+u_{i+1,j+1}^{x}h_{i+1,j}\right)\\
    &  & =\sum\frac{\Delta y\Delta x}{4\Delta y\Delta x} \, m_{i,j}^{y}\left(\left(u_{i+1,j+1}^{x}-u_{i+1,j-1}^{x}\right)h_{i+1,j}-\left(u_{i-1,j+1}^{x}-u_{i-1,j-1}^{x}\right)h_{i-1,j}\right)\\
    &  & =\sum\Delta y\Delta x \, m_{i,j}^{y}\left(\partial_{y,1\Delta}^{00}\left(h\partial_{x,1\Delta}^{00}u^{x}\right)\right)_{i,j}.
\end{eqnarray*}
So that with
\[
    \begin{array}{c}
        \sum\frac{\Delta y\Delta x}{4\Delta y\Delta x}u_{i,j}^{y}\left(h_{i+1,j}\left(m_{i+1,j+1}^{x}-m_{i+1,j-1}^{x}\right)-h_{i-1,j}\left(m_{i-1,j+1}^{x}-m_{i-1,j-1}^{x}\right)\right) \\
        =\sum\Delta y\Delta xm_{i,j}^{x}\left(\partial_{y,1\Delta}^{00}\left(h\partial_{y,1\Delta}^{00}u^{y}\right)\right)_{i,j}
    \end{array}
\]
and
\[
    \begin{array}{c}
        \sum\Delta y\Delta x\frac{1}{\left(\Delta y\right)^{2}}u_{i,j}^{y}\left(h_{i,j+1/2}\left(m_{i,j+1}^{y}-m_{i,j}^{y}\right)-h_{i,j-1/2}\left(m_{i,j}^{y}-m_{i,j-1}^{y}\right)\right) \\
        =\sum\Delta y\Delta xm_{i,j}^{y}\partial_{y,1\Delta}^{0}\left(h\left(\partial_{y,1\Delta}u^{y}\right)\right)_{i,j}.
    \end{array}
\]
We get finally
\[
    \begin{array}{c}
        \left(\boldsymbol{u},{\rm div}_{1,\Delta}\left(h\nabla_{1,\Delta}\boldsymbol{m}^{T}\right)\right)_{T_{d}}=\left(\boldsymbol{m},\left(\begin{array}{c}
                    \left(\partial_{x,1\Delta}^{0}\left(h\partial_{x,1\Delta}u^{x}\right)\right)+\left(\partial_{y,1\Delta}^{00}\left(h\partial_{x,1\Delta}^{00}u^{y}\right)\right) \\
                    \left(\partial_{x,1\Delta}^{00}\left(h\partial_{y,1\Delta}u^{x}\right)\right)+\left(\partial_{y,1\Delta}^{0}\left(h\partial_{y,1\Delta}u^{y}\right)\right)
                \end{array}\right)\right) \\
        =\left(\boldsymbol{m},div_{1,\Delta}\left(h\nabla_{1,\Delta}\boldsymbol{u}^{T}\right)\right)_{T_{d}}.
    \end{array}
\]

\begin{prop} \label{Stab_Cap} Let us suppose that ${\rm div}_{1,\Delta}\left(h\nabla_{1,\Delta}\boldsymbol{m}^{T}\right)$
    satisfies identity \eqref{eq:ident_dual_cap} of Lemma \ref{Dual_discret_compact},
    then the capillary step
    \[
        \mathbf{U}_{K}^{n+1}= \mathbf{U}_{K}^{n+1/2}
        +\Delta t\, \mathbf{M}_{K}^{n+1}
    \]
    admits a unique solution which satisfies an energy inequality:
    \begin{equation}
        \sum_{K}E_{tot}\left(\mathbf{U}_{K}^{n+1}\right)m_{K}
        \leq\sum_{K}E_{tot}\left(\mathbf{U}_{K}^{n+1/2}\right)m_{K}.\label{eq:cons_E_Cap}
    \end{equation}
\end{prop}

\noindent \textbf{Proof of Proposition \ref{Stab_Cap}.} Let us first
prove that the system \eqref{eq:Sys_Cap_dis} admits a unique solution.
Indeed, one can write $\mathbf{M}_{K}^{n+1}=  \mathcal{M}(h^{n+1},\boldsymbol{v}^{n+1/2})(\mathbf{U}_{K}^{n+1})$
and $\mathcal{M}(h^{n+1},{\bf v}^{n+1/2})$ satisfies
$$(\mathbf{U}, \mathcal{M}(h^{n+1},\boldsymbol{v}^{n+1/2})(\mathbf{U}))_{T_{d}}=0 \hbox{ for all } \mathbf{U}$$
from which we deduce that $ \mathcal{M}(h^{n+1},\boldsymbol{v}^{n+1/2})$ is
a skew-symmetric matrix for the scalar product $(\:.\:,\:.\:)_{T_{d}}$.
Thus its eigenvalues are purely imaginary and ${\rm Id}-\Delta t\,  \mathcal{M}(h^{n+1},{\bf v}^{n+1/2})$
is invertible. Now, thanks to identity (\ref{eq:Sys_Cap_dis}) and
the convexity of $E_{tot}$ (the fluid height $h$ is assumed $h>0$):
\[
    E_{tot}\left(\mathbf{U}_{K}^{n+1}\right)\leq E_{tot}\left(\mathbf{U}_{K}^{n+1/2}\right)-\Delta t\,\nabla_{\mathbf{U}}E_{tot}\left(\mathbf{U}_{K}^{n+1}\right)^{T}
    \mathbf{M}_{K}^{n+1}.
\]
Denote $f^{n+1/2}=f(h_{K}^{^{n+1}},\boldsymbol{v}_{K}^{n+1/2}),\;
    \mathbf{g}^{n+1/2}= \mathbf{g}(h_{K}^{{n+1}},\boldsymbol{v}_{K}^{n+1/2})$
and $DE:=\sum_{K}\,\nabla_{\mathbf{U}}E_{tot}
    \left(\mathbf{U}_{K}^{n+1}\right)^{T}\mathbf{M}_{K}^{n+1}m_{K}$.

\[
    R=-\left(\boldsymbol{u}_{K}^{n+1},\nabla_{3,\Delta}\left(g^{n+1/2}\boldsymbol{v}_{K}^{n+1}\right)\right)_{T_{d}}-\left(g^{n+1/2}\boldsymbol{v}_{K}^{n+1},{\rm div}_{3,\Delta}\left(\boldsymbol{u}_{K}^{n+1}\right)\right)_{T_{d}}
\]
and
\[
    D=\left(\boldsymbol{u}_{K}^{n+1},{\rm div}_{1,\Delta}\left(h_{K}^{{n+1}}\nabla_{1,\Delta}\left(f^{n+1/2}\boldsymbol{v}_{K}^{n+1}\right)^{T}\right)\right)_{T_{d}}-\left({\rm div}_{1,\Delta}\left(h_{K}^{{n+1}}\nabla_{1\Delta}(\boldsymbol{u}_{K}^{n+1})^{T}\right),f^{n+1/2}\boldsymbol{v}_{K}^{n+1}\right)_{T_{d}}
\]
We easily get that $DE=R+D$. As a consequence of definition \ref{eq:defGD3},
we get directly $R=0$, and, as a consequence of lemma \ref{Dual_discret_compact},
we get $D=0$. It follows that
\[
    \sum_{K}E_{tot}\left(\mathbf{U}_{K}^{n+1}\right)m_{K}
    \leq\sum_{K}E_{tot}\left(\mathbf{U}_{K}^{n+1/2}\right)m_{K}.
\]

\noindent We thus have proved the following stability result.

\begin{prop} Consider the scheme \eqref{eq:Hyp_part}--\eqref{eq:Sys_Cap_dis}--\eqref{eq:def_M_disc}
    with discretization \eqref{eq:def_d1D_compact} of capillary terms,
    then provided a CFL condition of the type \eqref{eq:CFL_Hyp} is satisfied,
    the fluid height $h$ is positive and the scheme satisfies energy
    stability
    \[
        \sum_{K}E_{tot}\left(\mathbf{U}_{K}^{n+1}\right)m_{K}
        \leq\sum_{K}E_{tot}\left(\mathbf{U}_{K}^{n}\right)m_{K}.
    \]
\end{prop}

%%Règles plus générale
This stability result can be extended to a more general numerical
framework and other time discretizations. By taking discrete dual
operators with similar rules as \eqref{eq:defGD3} namely
\[
    \left(\boldsymbol{w},{\rm div}_{1,\Delta}\left(\mathcal{T}\right)\right)_{T}=-\left(\nabla_{1,\Delta}\boldsymbol{w},\mathcal{T}\right)_{T_{d}}
\]
We thus get {\setlength{\arraycolsep}{1pt}
        \begin{eqnarray*}
            D &  & =\left(\boldsymbol{u}_{K}^{n+1},{\rm div}_{1,\Delta}\left(h_{K}^{{n+1}}\nabla_{1,\Delta}\left(f^{n+1/2}\boldsymbol{v}_{K}^{n+1}\right)^{T}\right)\right)_{T_{d}}-\left({\rm div}_{1,\Delta}\left(h_{K}^{{n+1}}\nabla_{1\Delta}(\boldsymbol{u}_{K}^{n+1})^{T}\right),f^{n+1/2}\boldsymbol{v}_{K}^{n+1}\right)_{T_{d}}\\
            &  & =-\left(\nabla_{1,\Delta}\boldsymbol{u}_{K}^{n+1},h_{K}^{{n+1}}\nabla_{1,\Delta}\left(f^{n+1/2}\boldsymbol{v}_{K}^{n+1}\right)^{T}\right)_{T_{d}}-\left({\rm div}_{1,\Delta}\left(h_{K}^{{n+1}}\nabla_{1\Delta}(\boldsymbol{u}_{K}^{n+1})^{T}\right),f^{n+1/2}\boldsymbol{v}_{K}^{n+1}\right)_{T_{d}}\\
            &  & =-\left(\nabla_{1,\Delta}\left(f^{n+1/2}\boldsymbol{v}_{K}^{n+1}\right),h_{K}^{{n+1}}\nabla_{1,\Delta}(\boldsymbol{u}_{K}^{n+1})^{T}\right)_{T_{d}}-\left({\rm div}_{1,\Delta}\left(h^{^{n+1}}\nabla_{1\Delta}^{t}\boldsymbol{u}_{K}^{n+1}\right),f^{n+1/2}\boldsymbol{v}_{K}^{n+1}\right)_{T_{d}}\\
            &  & =\left(f^{n+1/2}\boldsymbol{v}_{K}^{n+1},{\rm div}_{1,\Delta}\left(h_{K}^{{n+1}}\nabla_{1\Delta}(\boldsymbol{u}_{K}^{n+1})^{T}\right)\right)_{T_{d}}-\left({\rm div}_{1,\Delta}\left(h_{K}^{{n+1}}\nabla_{1\Delta}(\boldsymbol{u}_{K}^{n+1})^{T}\right),f^{n+1/2}\boldsymbol{v}_{K}^{n+1}\right)_{T_{d}}\\
            &  & =0
        \end{eqnarray*}
    } Condition (\ref{eq:ident_dual_cap}) of Lemma \ref{Dual_discret_compact} is valid and also insures energy stability result of Proposition \ref{Stab_Cap}.\\

\noindent One could also consider alternative time discretization
like the fully implicit scheme for the capillary step:
\begin{equation}
    \mathbf{U}^{n+1})= \mathbf{U}^{n+1/2}
    +\Delta t \, \mathcal{M} \left(h^{n+1}, \mathbf{U}^{n+1}\right) (\mathbf{U}^{n+1})\label{eq:Full_NL_Cap}
\end{equation}
This system could be solved through an iterative scheme:
\begin{equation}
    \mathbf{U}^{n+1,p+1}=\mathbf{U}^{n+1/2}
    +\Delta t  \, \mathcal{M} \left(h^{n+1},\mathbf{U}^{n+1,p}\right)(\mathbf{U}^{n+1,p+1}),
    \quad{\mathbf{U}}^{n+1,0}={\mathbf{U}}^{n+1/2}.\label{eq:Linearsysfull}
\end{equation}
The linear system \eqref{eq:Linearsysfull} admits a unique solution
which, moreover, satisfies the energy estimate
\[
    \sum_{K}E_{tot}\left(\mathbf{U}_{K}^{n+1,p}\right)m_{K}\leq\sum_{K}E_{tot}\left(\mathbf{U}_{K}^{n+1/2}\right)m_{K},\quad\forall p\geq0.
\]
If $\Delta t$ is small enough so that $\|\delta t\,  \mathcal{M}(h_{K}^{n+1},{\bf v}_{K}^{n+1/2})\|<1$,
the sequence $({\bf U}_{K}^{n+1,p})_{p\in\mathbb{N}}$ converges to
${\bf U}_{K}^{n+1}$ which, in turn, satisfies
\[
    \sum_{K}E_{tot}\left(\mathbf{U}_{K}^{n+1}\right)m_{K}
    \leq\sum_{K}E_{tot}\left(\mathbf{U}_{K}^{n+1/2}\right)m_{K}\leq\sum_{K}E_{tot}\left(\mathbf{U}_{K}^{n}\right)m_{K}.
\]

As a result, the IMplicit-EXplicit scheme build on a time discretization
with explicit steps for the hyperbolic part and implicit steps for
the capillary part are entropy stable. This provides a method to design
higher order in time IMplicit-EXplicit schemes which are build on
fully implicit time discretizations.

\section{Numerical Simulations}
We present in this section various numerical simulations to illustrate
the benefits of the proposed extended model. The new extended system  composed by a conservative hyperbolic part and a second order derivative part which is skew symmetric for the $L^2$ scalar product is crucial
to develop an appropriate splitting method.

We are able to carry out extremely fast simulations of capillary wave propagation in comparison to direct numerical simulations of the original Navier-Stokes equations (DNS). On the one hand, this is due to the vertical integration along the fluid height which reduces the dimension of the problem and withdraw the initial free surface problem. On the other hand, the implicit treatment of surface tension removes the classical restrictive capillary time step, empirically based to the fastest ``eligible'' wave speed
whose wavelength is the grid size. We will illustrate both the overall
stability of the numerical method and the interest of considering
the full surface tension source term.

Global energy dissipation will be shown on time discretizations that
are first order accurate. The time discretization is of IMplicit-EXplicit
type: for the hyperbolic part, an explicit Euler time-stepping scheme
has been used, associated with a Rusanov flux,
%\[
%G_{e,K}^{n}=G\left(\boldsymbol{U}_{e,K}^{n},\boldsymbol{U}_{e,K_{e}}^{n},\boldsymbol{n}_{e,K}\right)={\displaystyle \frac{F\left(\boldsymbol{U}_{e,K}^{n}\right)+F\left(\boldsymbol{U}_{e,K_{e}}^{n}\right)}{2}}+\max_{K,K_{e}}\left(\left|\boldsymbol{u}.\boldsymbol{n}_{e,K}\right|+{\displaystyle \sqrt{g_{r}h}}\right)\frac{\boldsymbol{U}_{e,K_{e}}^{n}-\boldsymbol{U}_{e,K}^{n}}{2},
%\]
\[
    G_{e,K}^{n}=G\left(\mathbf{U}_{e,K}^{n},\mathbf{U}_{e,K_{e}}^{n},\boldsymbol{n}_{e,K}\right)={\displaystyle \frac{F\left(\mathbf{U}_{e,K}^{n}\right)+F\left(\mathbf{U}_{e,K_{e}}^{n}\right)}{2}}-\max_{K,K_{e}}\left(\left|\boldsymbol{u}.\boldsymbol{n}_{e,K}\right|+{\displaystyle \sqrt{g_{r}h}}\right)\frac{\mathbf{U}_{e,K_{e}}^{n}-\mathbf{U}_{e,K}^{n}}{2},
\]

using the rotational invariance and considering second-order in space
MUSCL reconstructions denoted by $\mathbf{U}_{e,K}^{n}$ and $\mathbf{U}_{e,K_{e}}^{n}$
of the primitive variables (without limitation as very smooth solution
will be considered here) whereas an implicit Euler time-stepping
scheme is used for the capillary step, by considering a simpler linearized
resolution of the initial fully nonlinear problem of coupled equations.
While other reconstruction method are possible and can lead to
higher order accuracy, the MUSCL method is a second order method that is easy
to implement and has the advantage to be suitable with structured and
unstructured mesh. The model has been coded on an unstructured environment,
and it is planned to run simulations with it once we have solved the
discrete duality problem in such context.

It should be noted that a global second-order solver can be derived
by considering an appropriate IMEX time-stepping scheme to combine
the explicit and implicit steps but this strategy is costly as it
requires to solve the full nonlinear problem, that can be achieved
using Newton-like method or simply iterating on the linearized version
of the initial full nonlinear problem of coupled equations until convergence.

\subsection{Numerical Set Up}

We consider a rectangular domain $\left[0,l_{x}\right]\times\left[0,l_{y}\right]$
divided into $n_{x}\times n_{y}$ cells considering uniform discretization
steps $\Delta x$ and $\Delta y$ respectively in each direction.
In a Finite Volume framework, the $h_{i,j}$ and $\boldsymbol{{u}_{i,j}}$
discrete unknowns are associated classically to the mean of respectively
a scalar field $h$ and a vector field $\boldsymbol{u}$ over the appropriate
cell. In order to avoid any specific treatment of boundary conditions,
we have only considered \textit{periodic boundary conditions}.

We have carried out numerical simulations of the augmented version
of the shallow water equations in two situations: the quadratic capillary
case and the fully nonlinear capillary case. In the quadratic case,
the system (\ref{def_U_F},\ref{eq:def_M}) is written with,
\begin{equation}
    E\left(h,\nabla h\right)={\displaystyle g_{r}\frac{h}{2}+\frac{1}{2}\frac{\sigma}{\rho}\|\nabla h\|^{2}},\label{eq:quad_case}
\end{equation}
meaning,
\[
    \mathcal{E_{\mathrm{cap}}}\left(q\right)={\displaystyle \frac{1}{2}q^{2}},\quad\kappa\left(h\right)={\displaystyle \frac{\sigma}{\rho}},\quad\alpha\left(q^{2}\right)=1,
\]
and
\[
    {\displaystyle f\left(h,\boldsymbol{v}\right)={\displaystyle \sqrt{h}\sqrt{\frac{\sigma}{\rho}}\,\mathrm{I_d}},\quad \boldsymbol{g}\left(h,\boldsymbol{v}\right)={\displaystyle {\displaystyle \frac{h\boldsymbol{v}}{2}}}},
\]
where $g_{r}$, $\sigma$ and $\rho$ are respectively the constant
gravity acceleration, the surface tension coefficient and the constant
density of the flow. In the fully nonlinear capillary case, the system
(\ref{def_U_F},\ref{eq:def_M}) is defined with,
\begin{equation}
    E\left(h,\nabla h\right)={\displaystyle g_{r}\frac{h}{2}+\frac{\sigma}{\rho}(\sqrt{1+\|\nabla h\|^{2}}-1)},\label{eq:full_nl_case}
\end{equation}
meaning,
\[
    \mathcal{E_{\mathrm{cap}}}\left(q\right)=\sqrt{1+q^{2}}-1,\quad\kappa\left(h\right)=\frac{\sigma}{\rho},\quad\alpha\left(q^{2}\right)=\sqrt{2}\left(1+\sqrt{1+q^{2}}\right)^{-1/2},
\]
and
\[
    f\left(h,\boldsymbol{v}\right)=\sqrt{h}{\displaystyle \left(1+{\displaystyle \frac{\rho h}{4\sigma}\left\Vert \boldsymbol{v}\right\Vert ^{2}}\right)^{-1/2}}\left(\mathrm{I_d}-\left(1+{\displaystyle \frac{\rho h}{2\sigma}\left\Vert \boldsymbol{v}\right\Vert ^{2}}\right)^{-1}{\displaystyle \frac{\rho h}{4\sigma}}\boldsymbol{v}\otimes\boldsymbol{v}\right),\quad \boldsymbol{g}\left(h,\boldsymbol{v}\right)={\displaystyle {\displaystyle \frac{h\boldsymbol{v}}{2}}\left(1+{\displaystyle \frac{\rho h}{2\sigma}\left\Vert \boldsymbol{v}\right\Vert ^{2}}\right)^{-1}}.
\]
We recall that the expression of $\boldsymbol{v}$ as a function of
$\alpha$ and $\kappa$ is given in Equation \eqref{v}.

\subsection{One-dimensional simulation with Gaussian initial data\label{subsec:One-dimensional-Gaussian}}
\begin{figure}[!ht]
    \centering{}\includegraphics[width=0.4\textwidth]{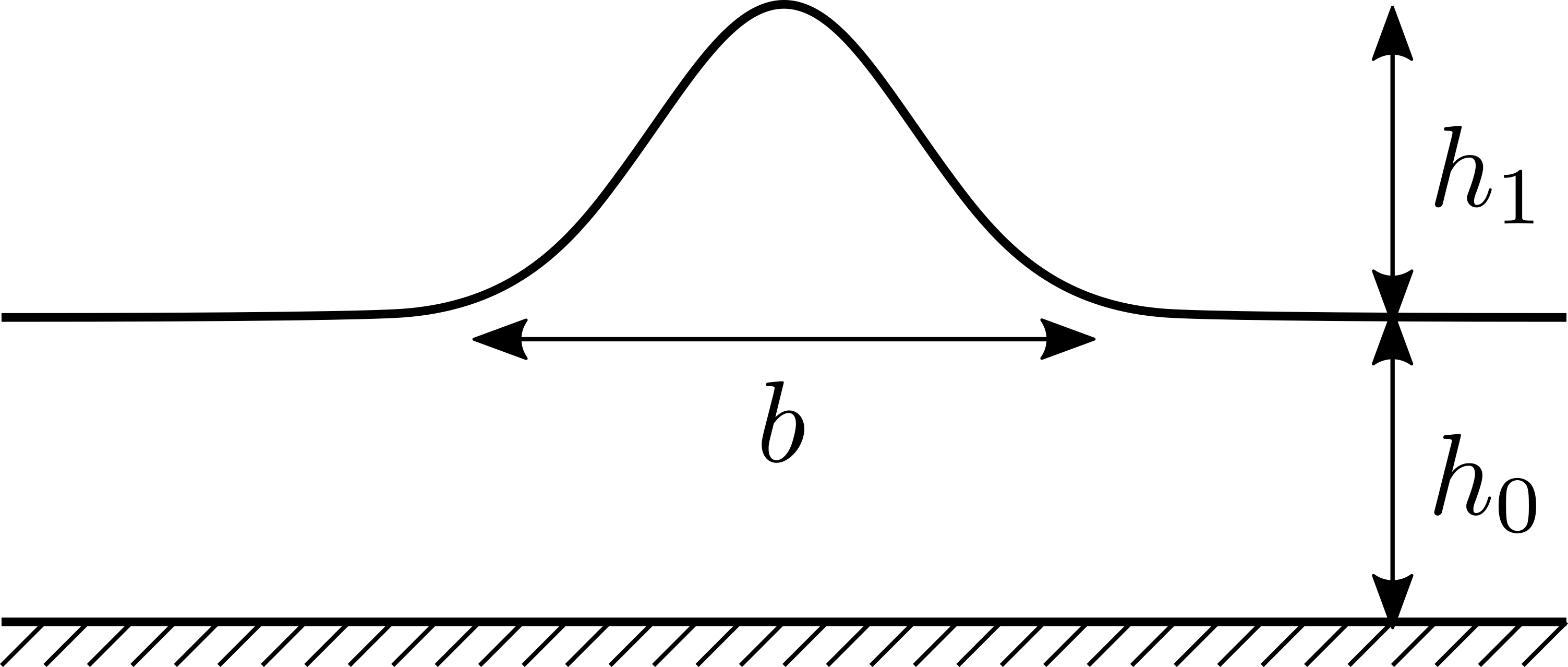}\caption{One-dimensional sketch of the Gaussian deformation of a layer of water
        where $b$ is the full width at tenth of maximum (FWTM).\label{fig:1d-gaussian}}
\end{figure}

We consider a one-dimensional Gaussian-shaped deformation of
the free surface of a water layer, as illustrated in Figure \eqref{fig:1d-gaussian}.
This deformation produces both gravity and capillary waves whose relative
influence is measured by the E\"{o}tv\"{o}s number, also called
Bond number,
\begin{equation}
    \mathrm{Eo}=\mathrm{Bo}={\displaystyle \frac{\rho g_{r}h^{2}}{\sigma}},\label{eq:Eotvos}
\end{equation}
as long as the shape of the Gaussian is close to the shape of a drop,
i.e. its curve peak height is comparable to its width. We set physical
parameters to the conventional values for water at room temperature
and are summarized in Table \eqref{tab:Physical-parameters}.

\begin{table}[!ht]
    \begin{centering}
        \begin{tabular}{l}
            $\sigma=0.0728$ \SI{}{\newton\per\meter}\tabularnewline
            $\rho=1000$ \SI{}{\kilogram\per\cubic\meter}\tabularnewline
            $\nu=10^{-6}$ \SI{}{\square\meter\per\second}\tabularnewline
            $g_{r}=9.81$ \SI{}{\meter\per\square\second}\tabularnewline
        \end{tabular}
        \par\end{centering}
    \caption{Physical parameters for the simulations.\label{tab:Physical-parameters}}
\end{table}
\noindent The initial Gaussian-shaped deformation of the water layer
parametrizes the initial surface elevation as,
\begin{equation}
    {\displaystyle h(x,t=0)=h_{0}+h_{1}e^{{\displaystyle -\frac{x^{2}}{2\left(b/b_{0}\right)^{2}}}},}
\end{equation}
with $b_{0}=4.29193$ allowing to consider approximately the full
width at tenth of maximum as the length $b$ represented in the Figure \eqref{fig:1d-gaussian}.
As the E\"{o}tv\"{o}s number Eq.\eqref{eq:Eotvos}is set to 1,
such that gravity and capillary waves are generated in the same time
order, this gives a water deformation peak elevation $h_{1}=2.725$
\SI{}{\milli\metre}. The layer of water elevation at rest is
set to $h_{0}=h_{1}$ whereas the full width at tenth of maximum is
set to $b=1.5\,h_{1}$. The computational domain is set to $\left[-50\SI{}{\milli\metre},50\SI{}{\milli\metre}\right]$
and the simulation time to $5$ \SI{}{\milli\second} in order
to produce significant waves in order to compare the results with
the two models with respectively a linearized capillary contribution
and a full nonlinear capillary contribution. Finally, the initial
velocity is set to zero and the auxiliary variable $\boldsymbol{v}$
is initialized through the formulas according to the two models considered.
In practice, it is not needed to compute exactly $\nabla h$, a
simple discretization using a classic centered scheme for example
is sufficient and used here in practice.

\begin{figure}[!ht]
    \begin{centering}
        \includegraphics[width=0.5\textwidth]{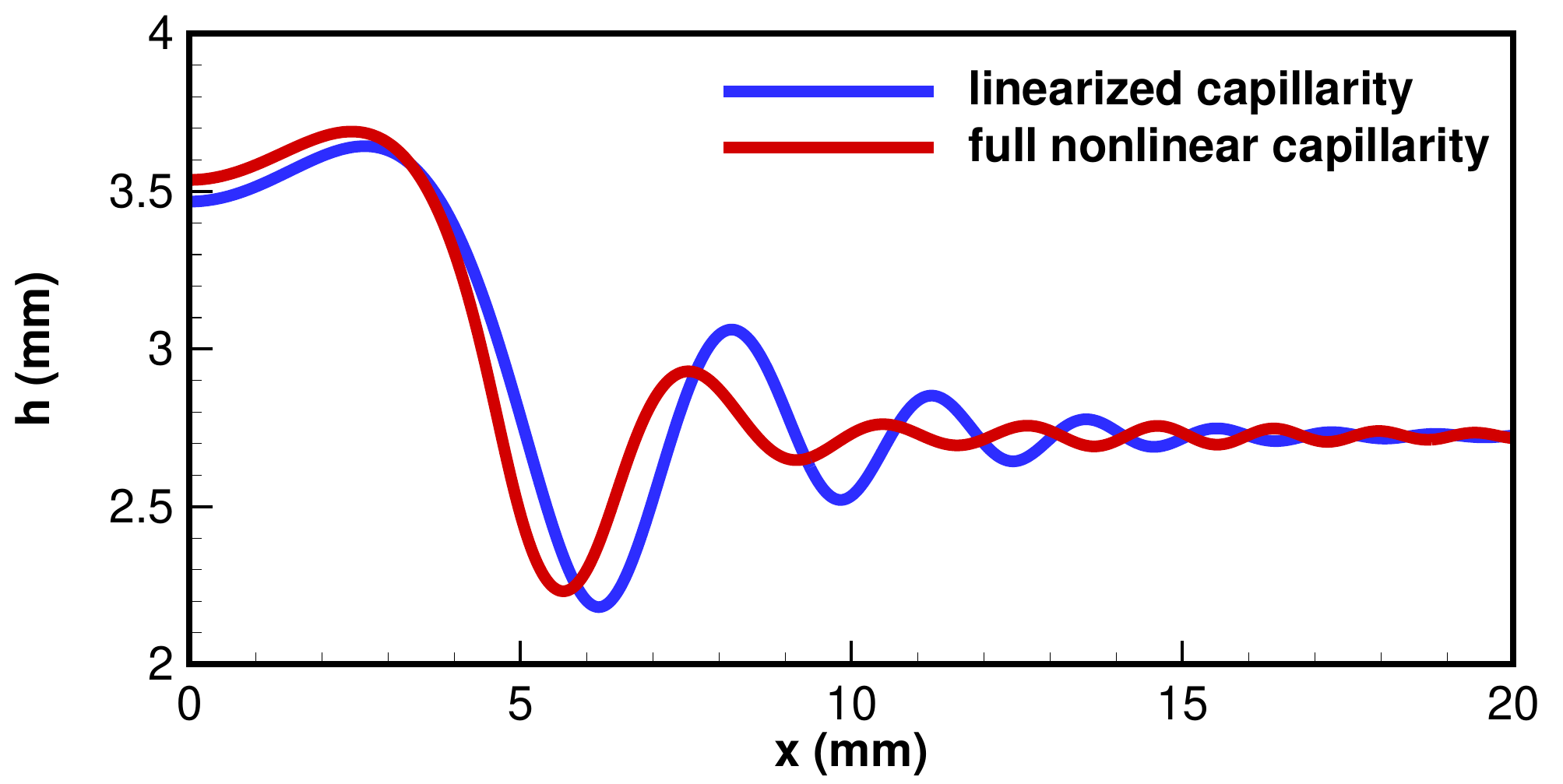}\includegraphics[width=0.5\textwidth]{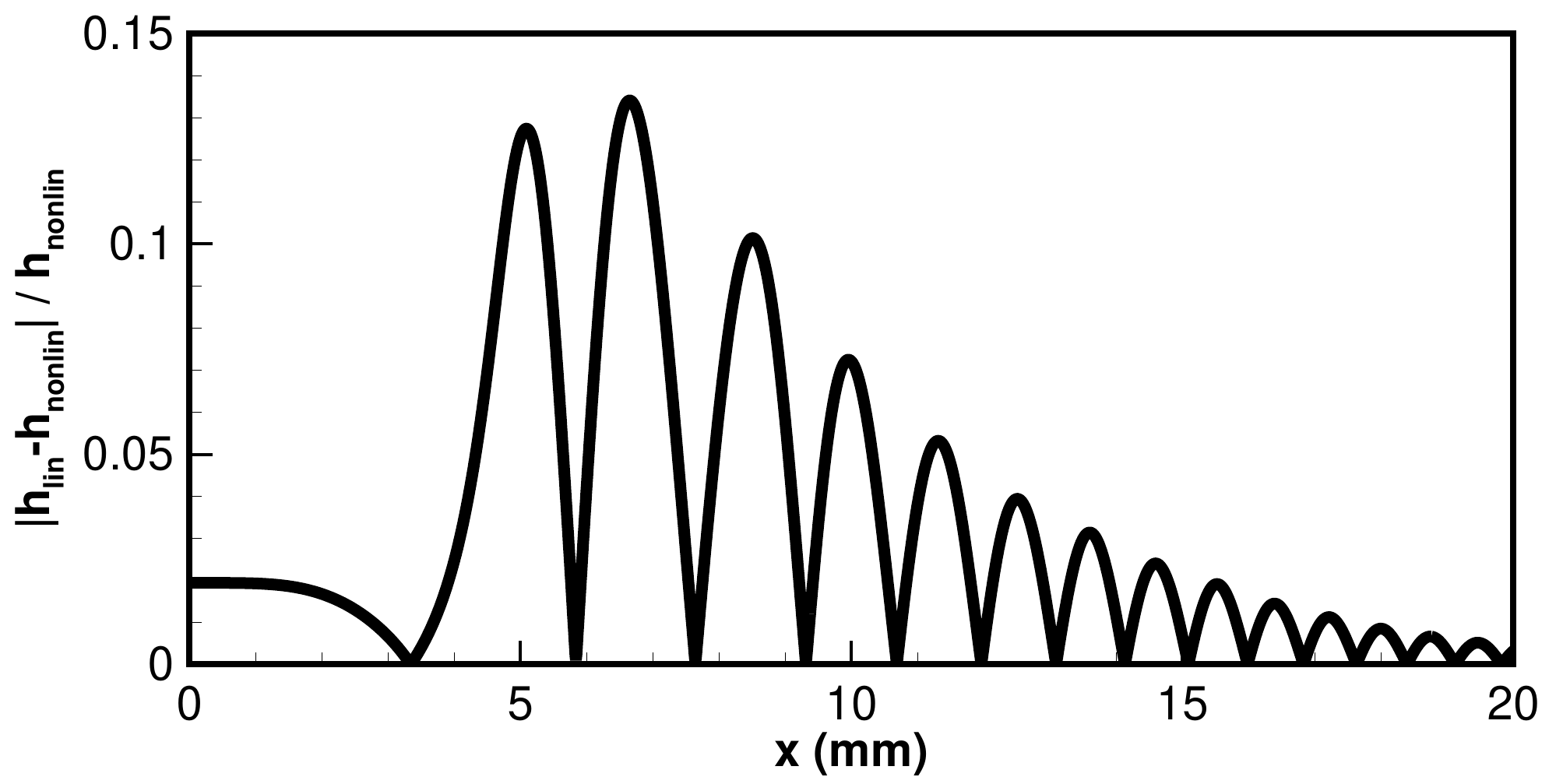}
        \par\end{centering}
    \begin{centering}
        \includegraphics[width=0.5\textwidth]{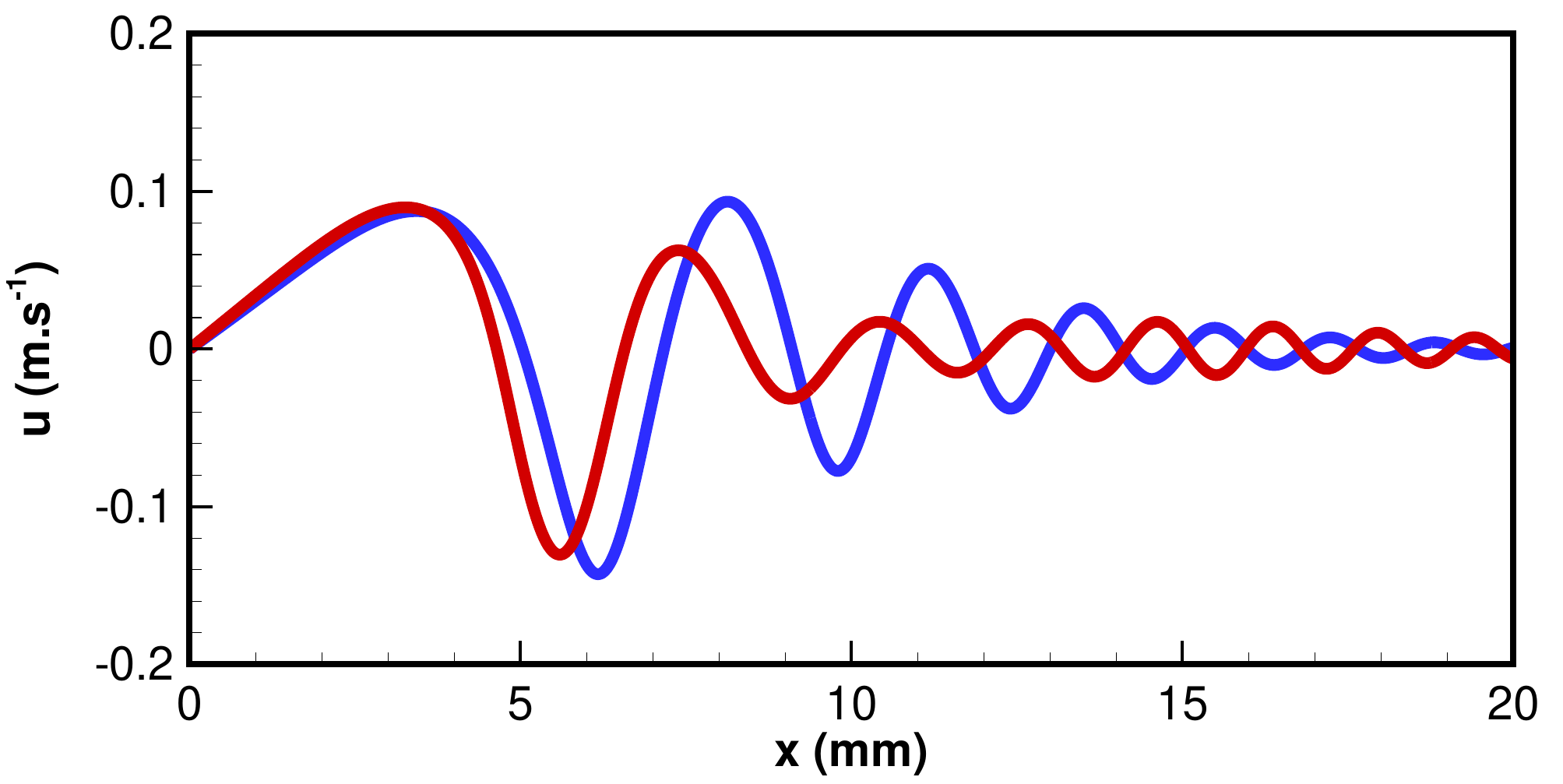}\includegraphics[width=0.5\textwidth]{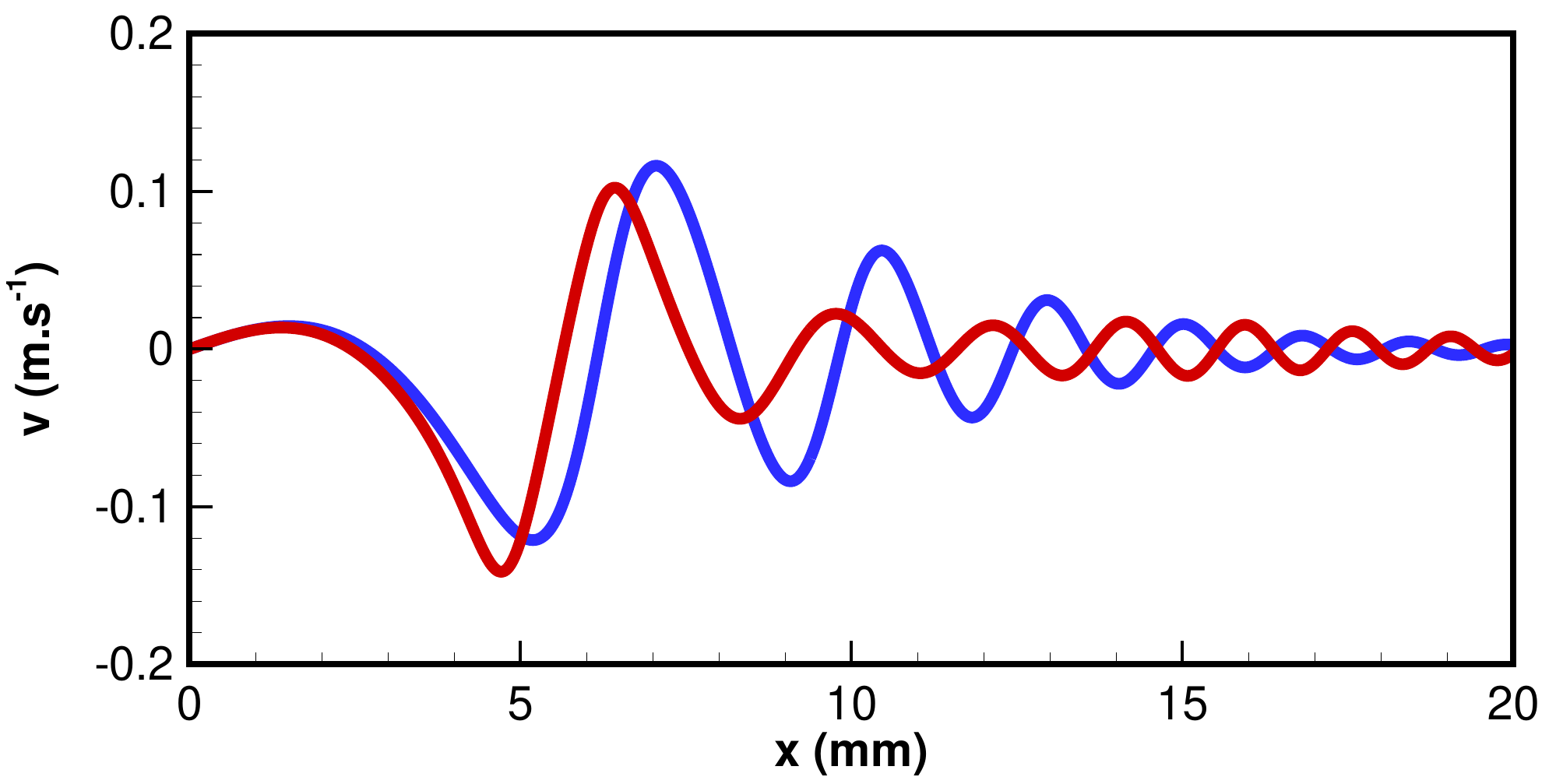}
        \par\end{centering}
    \caption{Very fine (51200 cells) resolved numerical simulations of capillary-gravity waves
        considering a one-dimensional Gaussian-shaped deformation of a layer
        of water using the proposed augmented shallow-water model Eq.(\ref{def_U_F},\ref{eq:def_M})
        and formulas Eq.(\ref{eq:quad_case},\ref{eq:full_nl_case}). Only
        a window of the real computational domain is plotted since the simulations
        are symmetric around zero and the waves not significant far away from
        zero; (top-left) Water height $h$; (top-right) Relative difference
        between the two water heights; (bottom-left) Velocity $\boldsymbol{u}$;
        (bottom-right) Auxiliary velocity $\boldsymbol{v}$.\label{fig:1d-diff-lin-nonlin}}
\end{figure}
It is presented in Figure (\ref{fig:1d-diff-lin-nonlin}) the very
fine resolved results for the water height $h$, the velocity $\boldsymbol{u}$
and the auxiliary velocity $\boldsymbol{v}$ considering the two proposed
models. For the physical parameters ans space scaling chosen, there
is a significant difference between the two models since the gradient
of the water height $\nabla h$ is sufficiently large to observe such
a behaviour. The computation of the relative difference between the
water height of each model shows an approximate maximal difference
of $14\:\%$. This is not only due to the difference in the capillary
wave amplitude, but also to an important phase shift.

\begin{figure}[!ht]
    \begin{centering}
        \includegraphics[width=0.5\textwidth]{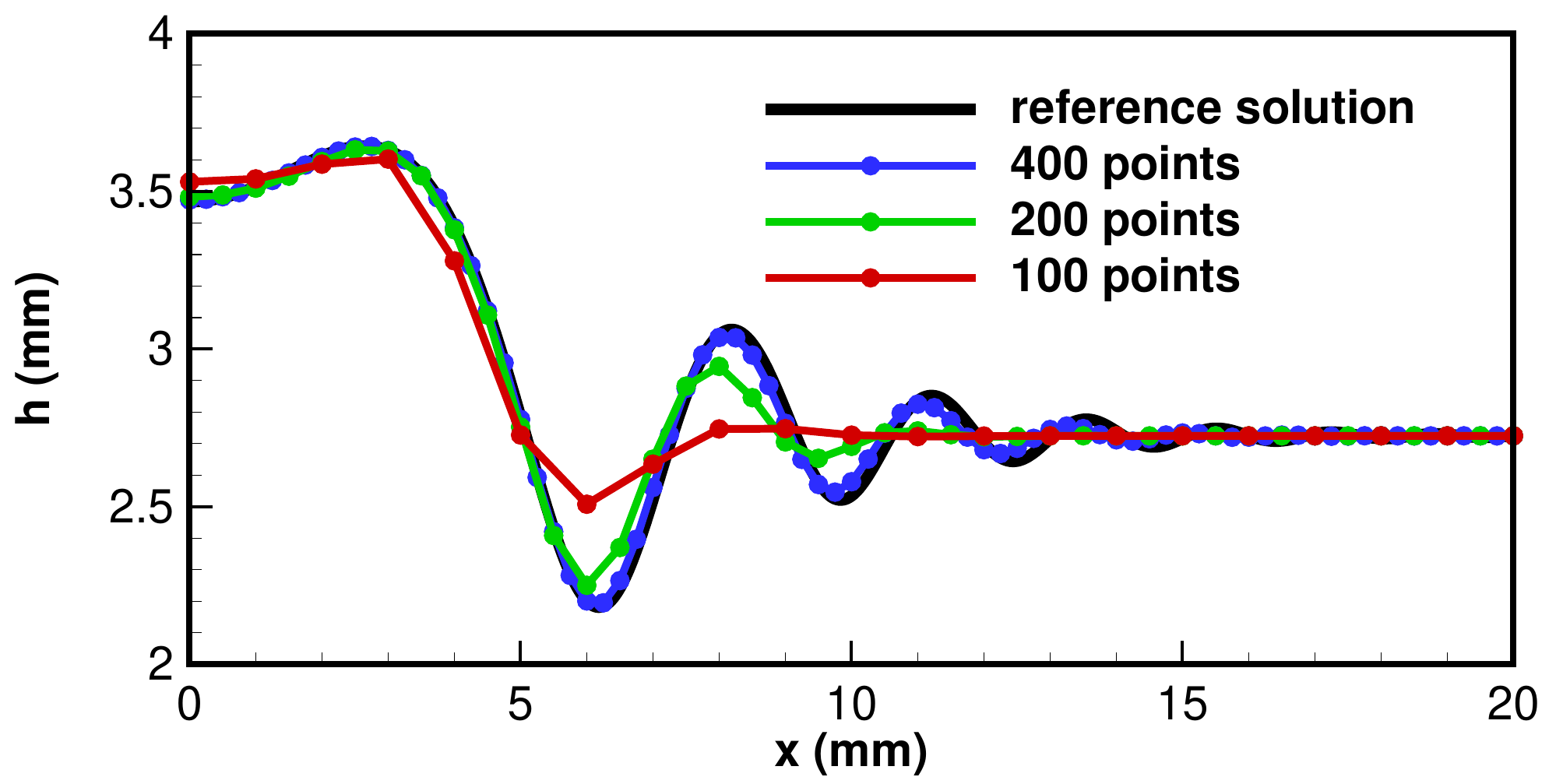}\includegraphics[width=0.5\textwidth]{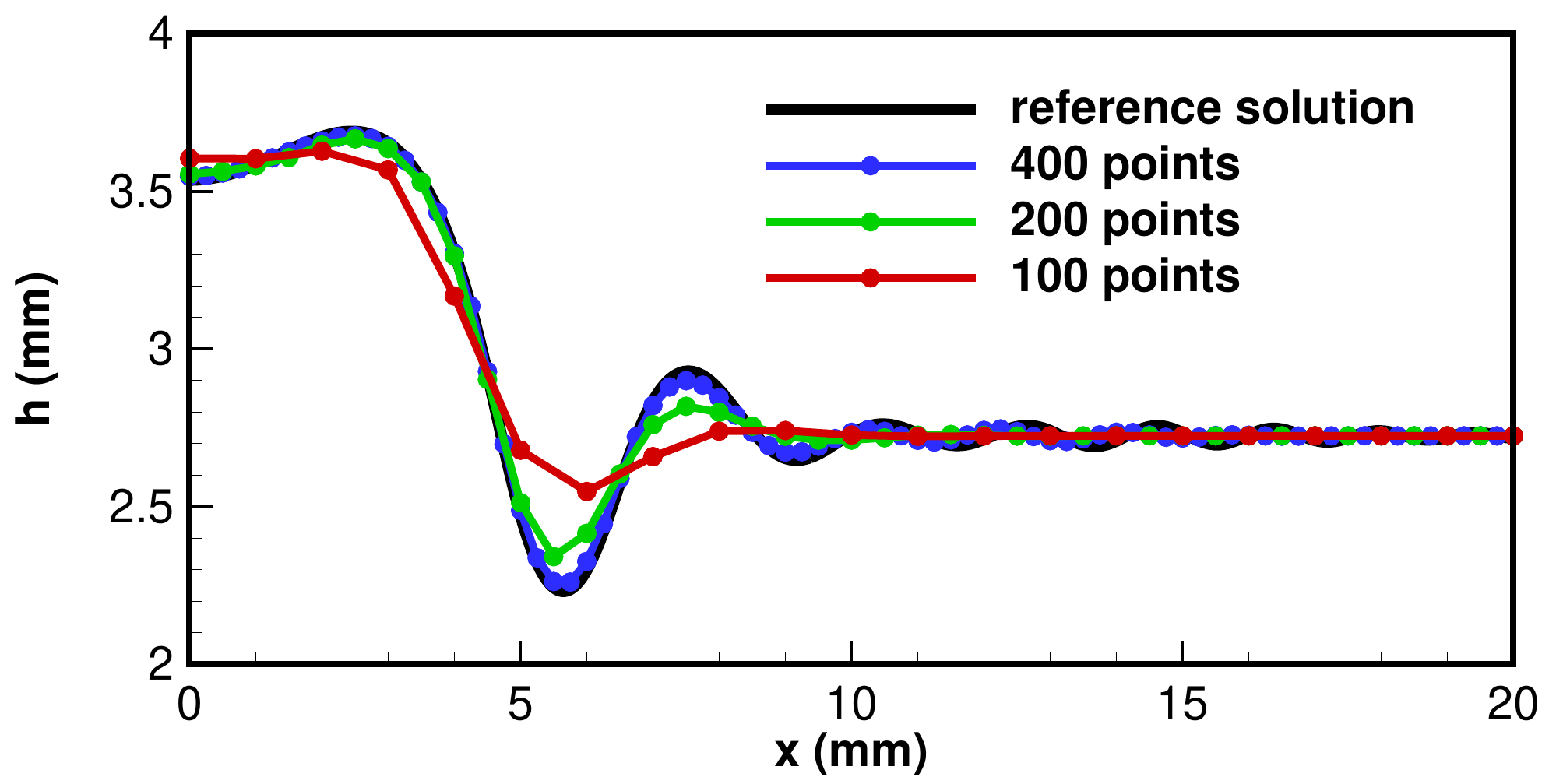}
        \par\end{centering}
    \begin{centering}
        \includegraphics[width=0.5\textwidth]{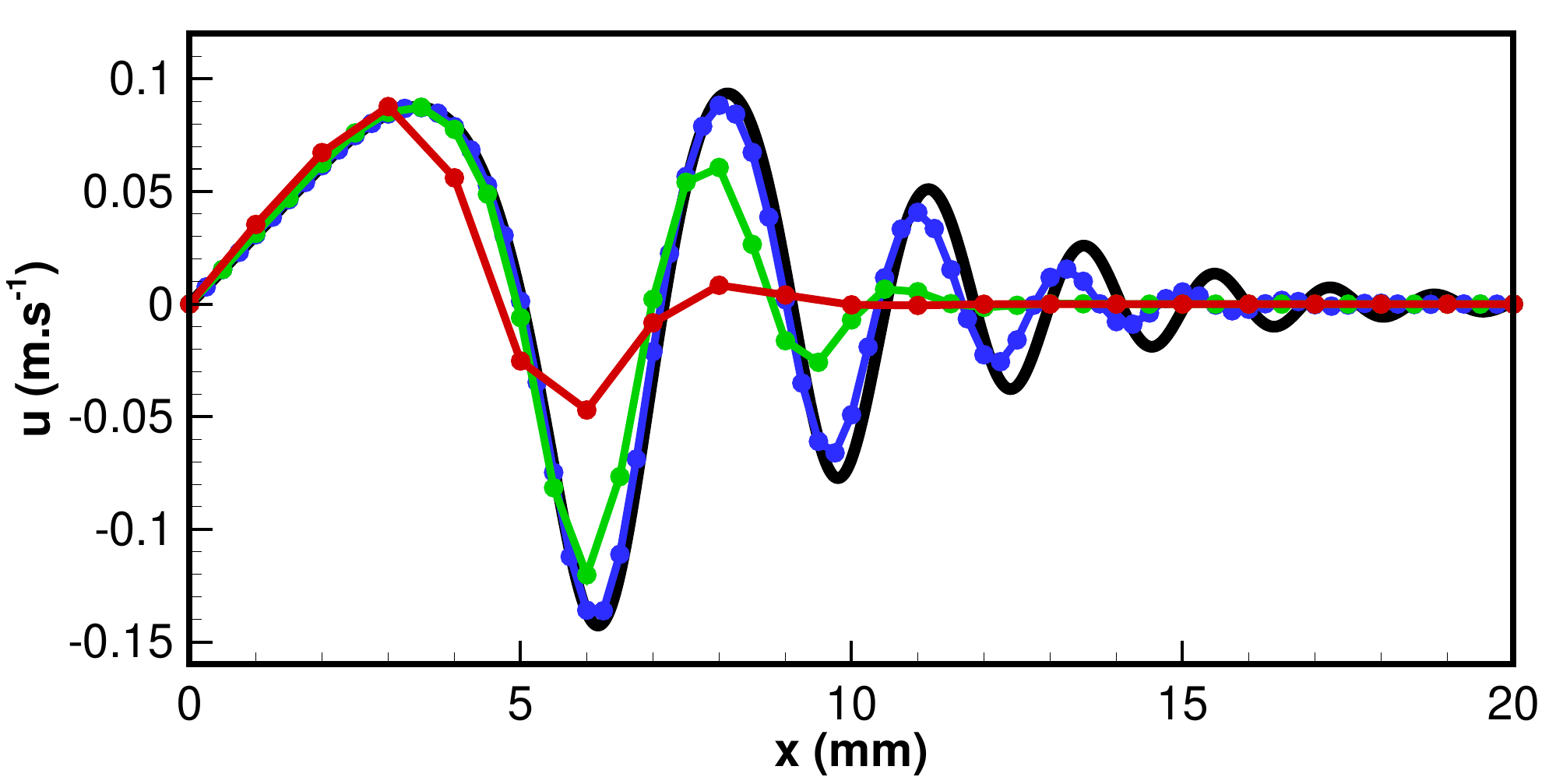}\includegraphics[width=0.5\textwidth]{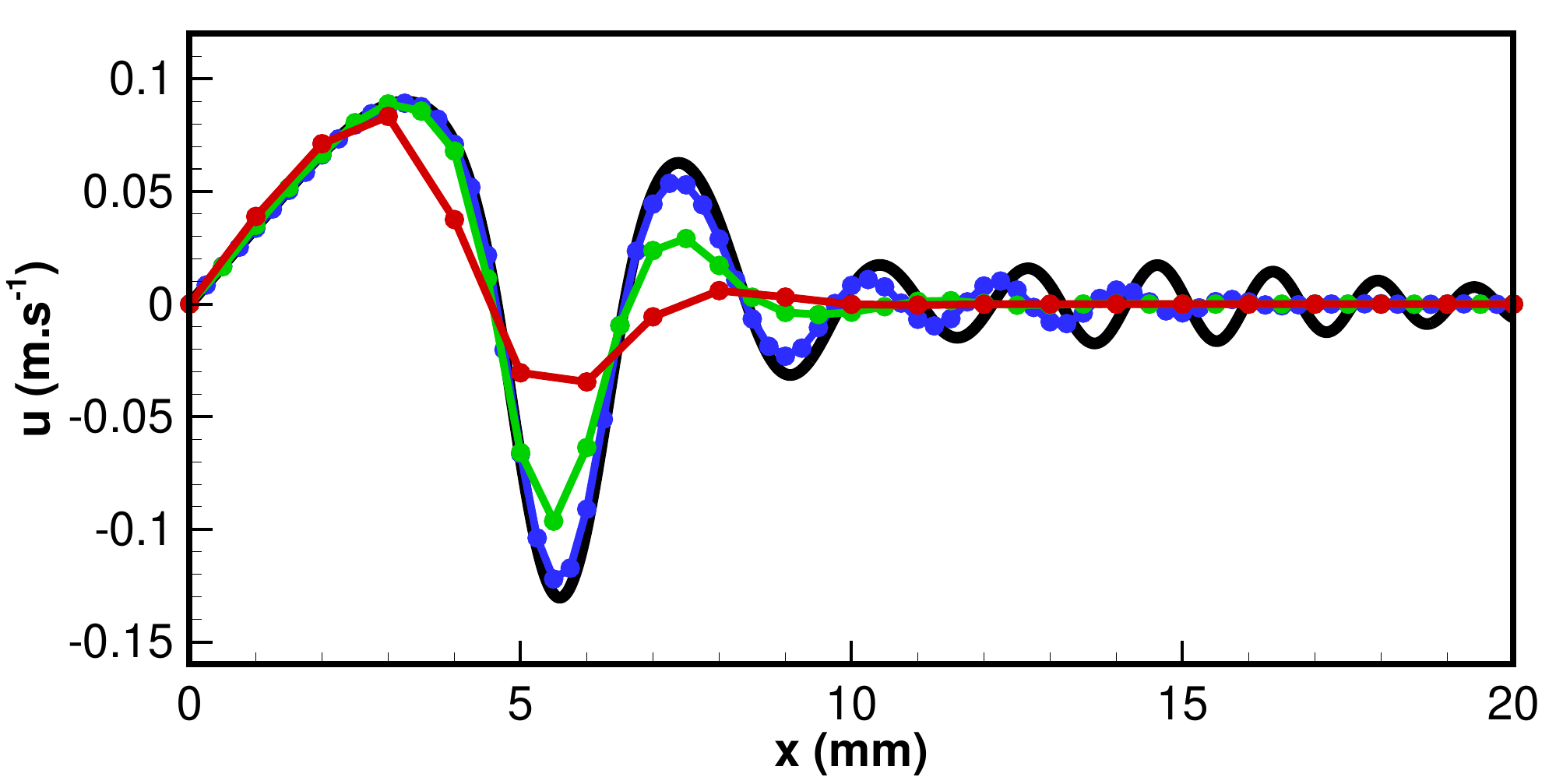}
        \par\end{centering}
    \begin{centering}
        \includegraphics[width=0.5\textwidth]{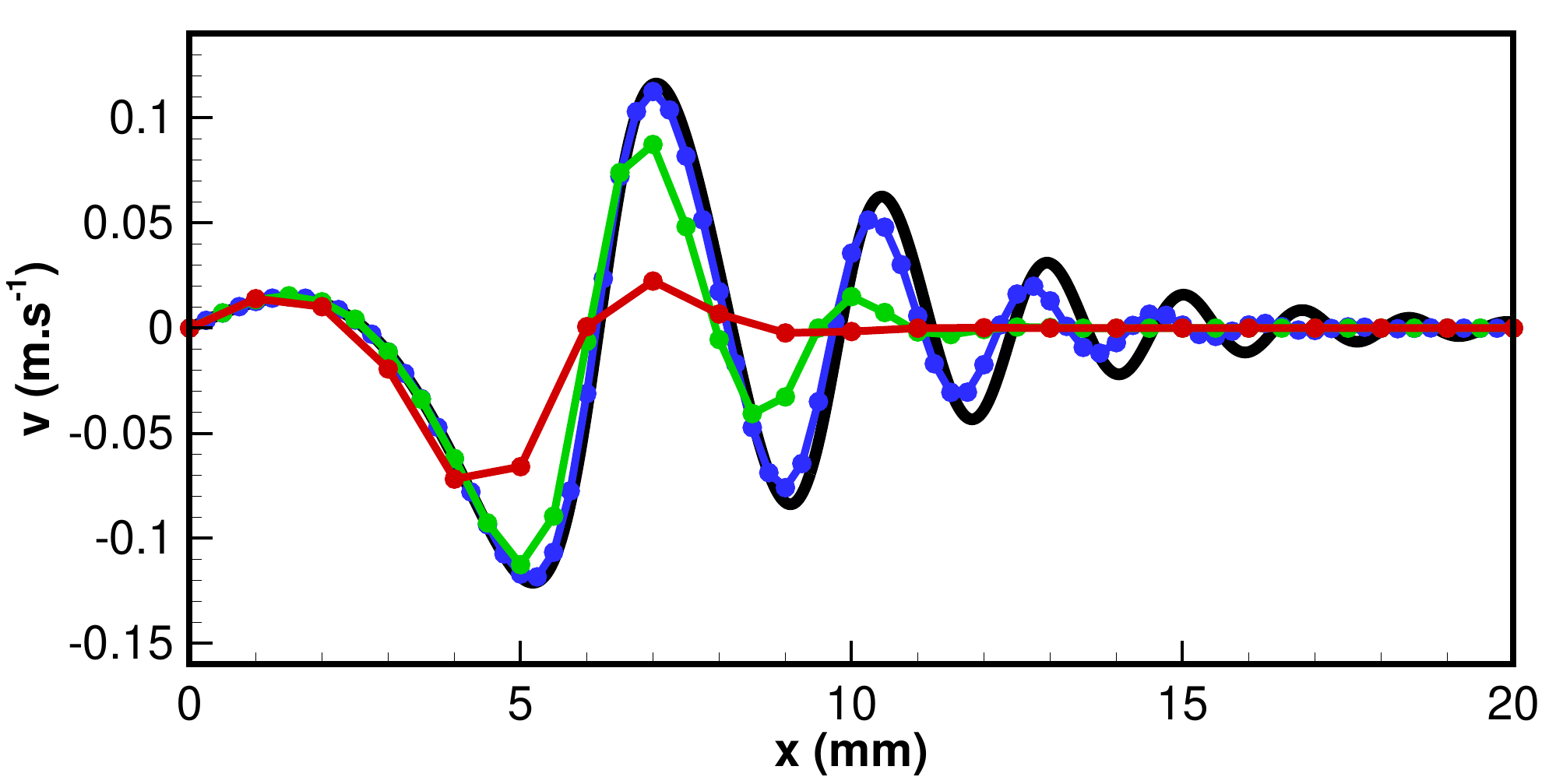}\includegraphics[width=0.5\textwidth]{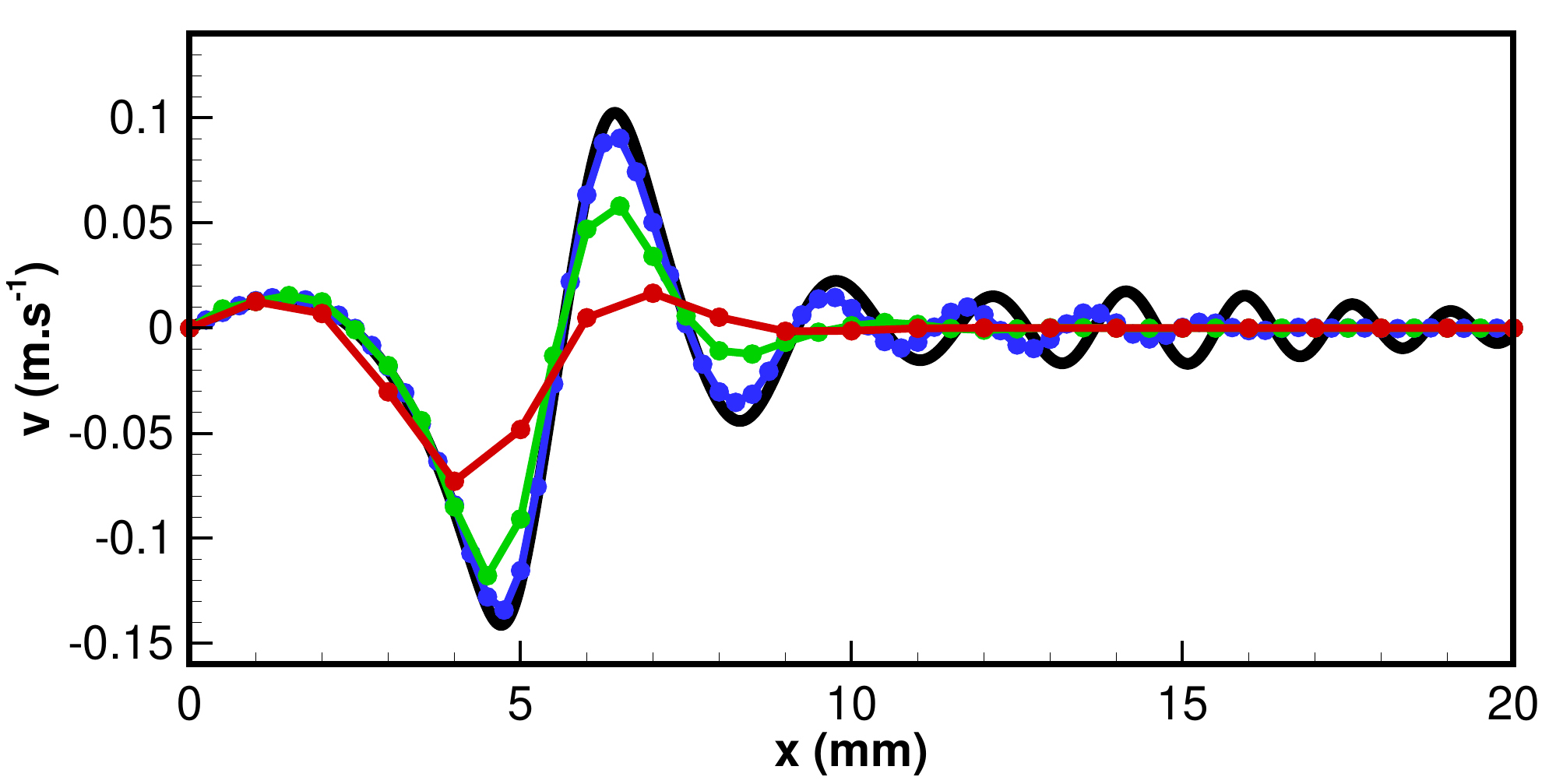}
        \par\end{centering}
    \caption{Same simulations than for Figure \eqref{fig:1d-diff-lin-nonlin},
        but for the first grid sizes considered in order to materialize the
        numerical solution quality with a growing number of discretization
        points in the characteristic wavelegth; (left) with the linearized
        capillary contribution; (right) with the full nonlinear capillary
        contribution.\label{fig:1d-simulations}}
\end{figure}
The computational simulation time being $5$ \SI{}{\milli\second},
it can also be observed that the capillary waves phase velocity are
much larger than the fluid velocity. This can be easily explained
by studying the dispersion relation, developed around a layer with
a height $h_{0}$ and a zero velocity, giving a wave speed,

\begin{equation}
    c\approx u\pm{\displaystyle \sqrt{g_{r}h_{0}+{\displaystyle \frac{h_{0}\sigma}{\rho}k^{2}}}}.
\end{equation}

where $k$ denotes the wave number of a plane wave. The ratio between
the capillary wave speed and gravity wave speed is then approximately
equal to ${\displaystyle \sqrt{\sigma/g_{r}\rho}\:2\pi/\lambda}\approx0.017/\lambda$,
where $\lambda$ is the characteristic wavelength of the surface elevation.
As the Fourier transform of an initial Gaussian-shape deformation
is again a Gaussian, there are wavelengths as small as the machine
accuracy allows to capture. Thus, for plane waves with a wavelength
of $0.17$ \SI{}{\milli\metre}, the capillary wave speed is
100 times faster than the gravity wave speed. This is the reason why
we have chosen a CFL number based on the maximal absolute eigenvalue
of the hyperbolic Jacobian matrix at an arbitrary value of $0.01$
in order to capture the propagation of the capillary waves.
Whereas proposed numerical discretization allows to work with higher CFL numbers
close to 1 due to the implicit resolution of the source terms modelling
the full contribution of the surface tension, the induced larger time
steps imply a numerical time capturing low pass filter regarding the
capillary waves. Another numerical viewpoint of using CFL numbers
close to 1 is that the induced linear system resolution becomes more
difficult due to a growing condition number of the resulted matrix
with larger time steps. In other words, the numerical resolution is
computationally more expensive whereas less physical phenomenon of
the capillary action is captured.

\begin{figure}[!hb]
    \begin{centering}
        \includegraphics[width=0.5\textwidth]{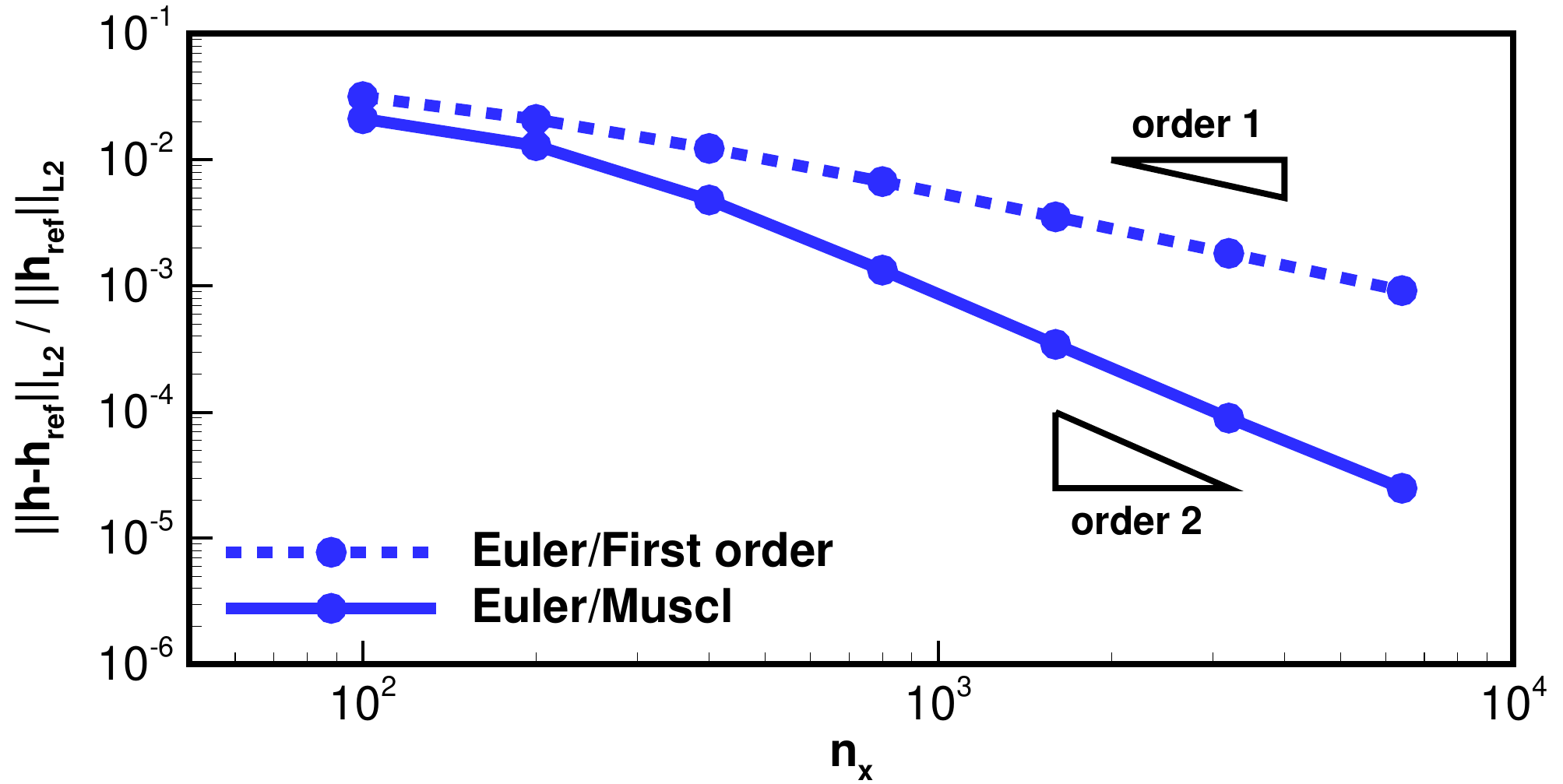}\includegraphics[width=0.5\textwidth]{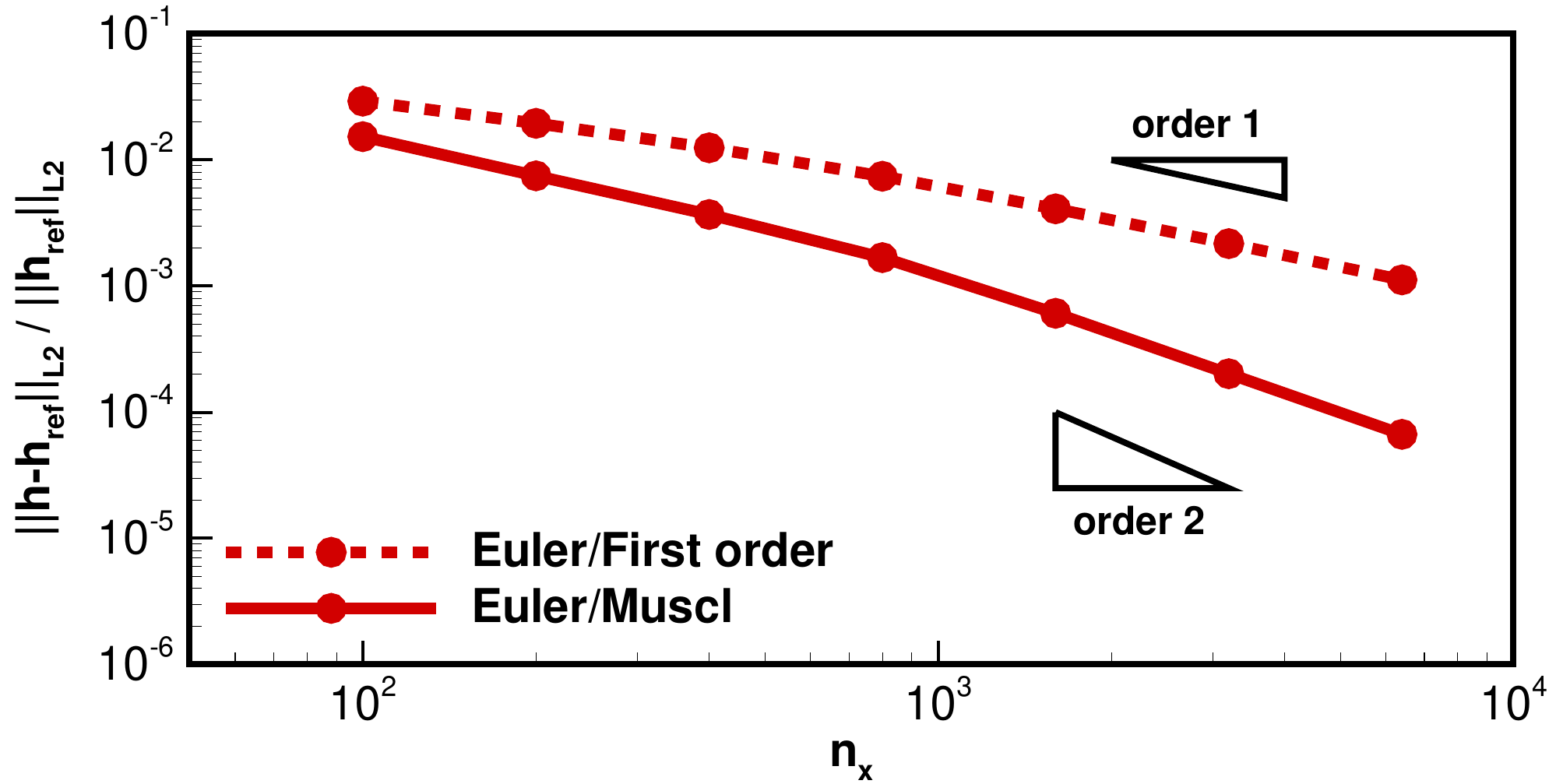}
        \par\end{centering}
    \caption{Relative error for the water height $h$ in the $L_{2}$ norm as a
        function of the grid size, computed at the end of the simulation given
        a reference solution $h_{\text{ref}}$ computed with $51200$ points;
        (left) with the linearized capillary contribution; (right) with the
        full nonlinear capillary contribution.\label{fig:conv_h}}
\end{figure}
\begin{figure}[!hb]
    \begin{centering}
        \includegraphics[width=0.5\textwidth]{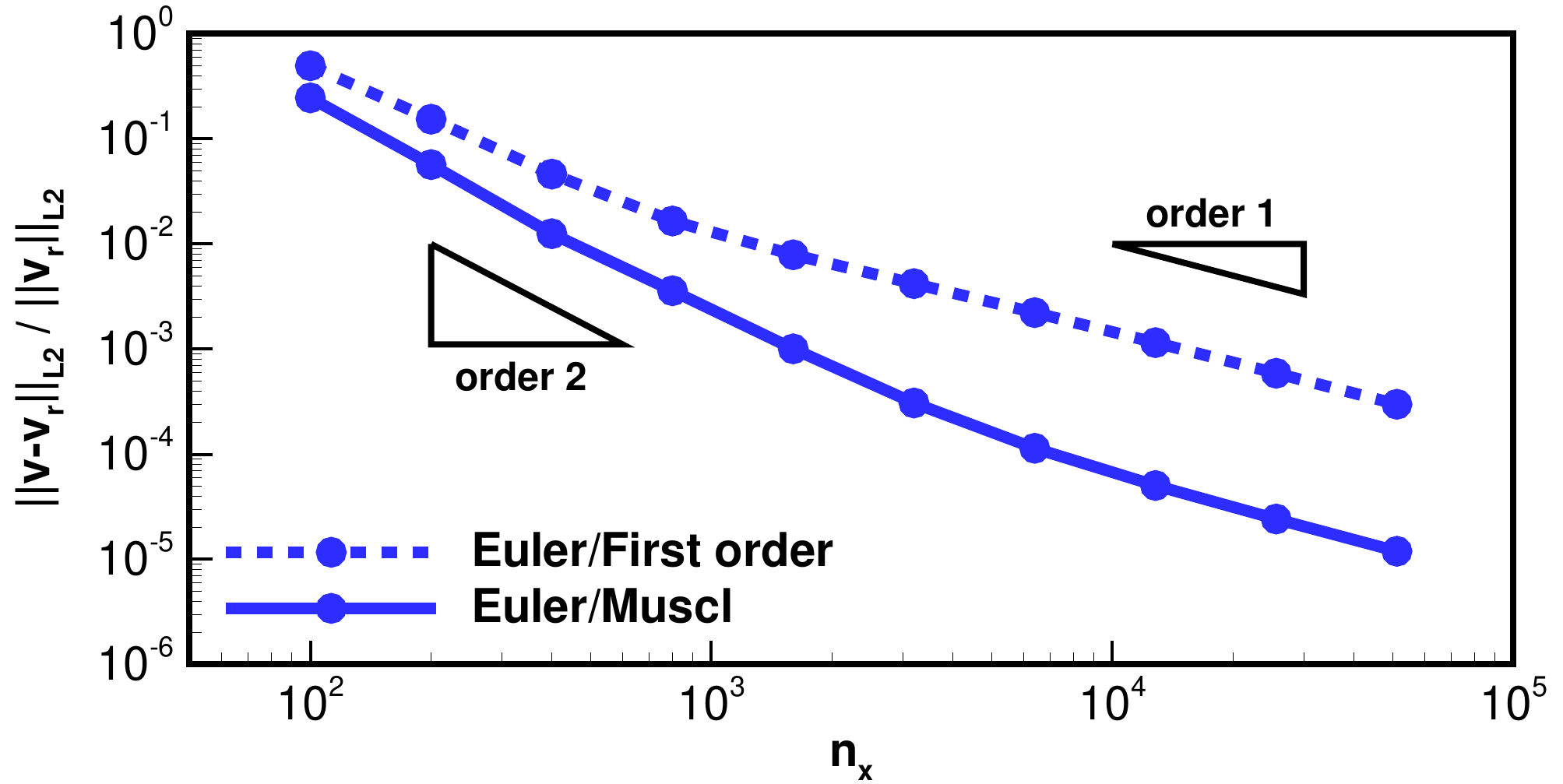}\includegraphics[width=0.5\textwidth]{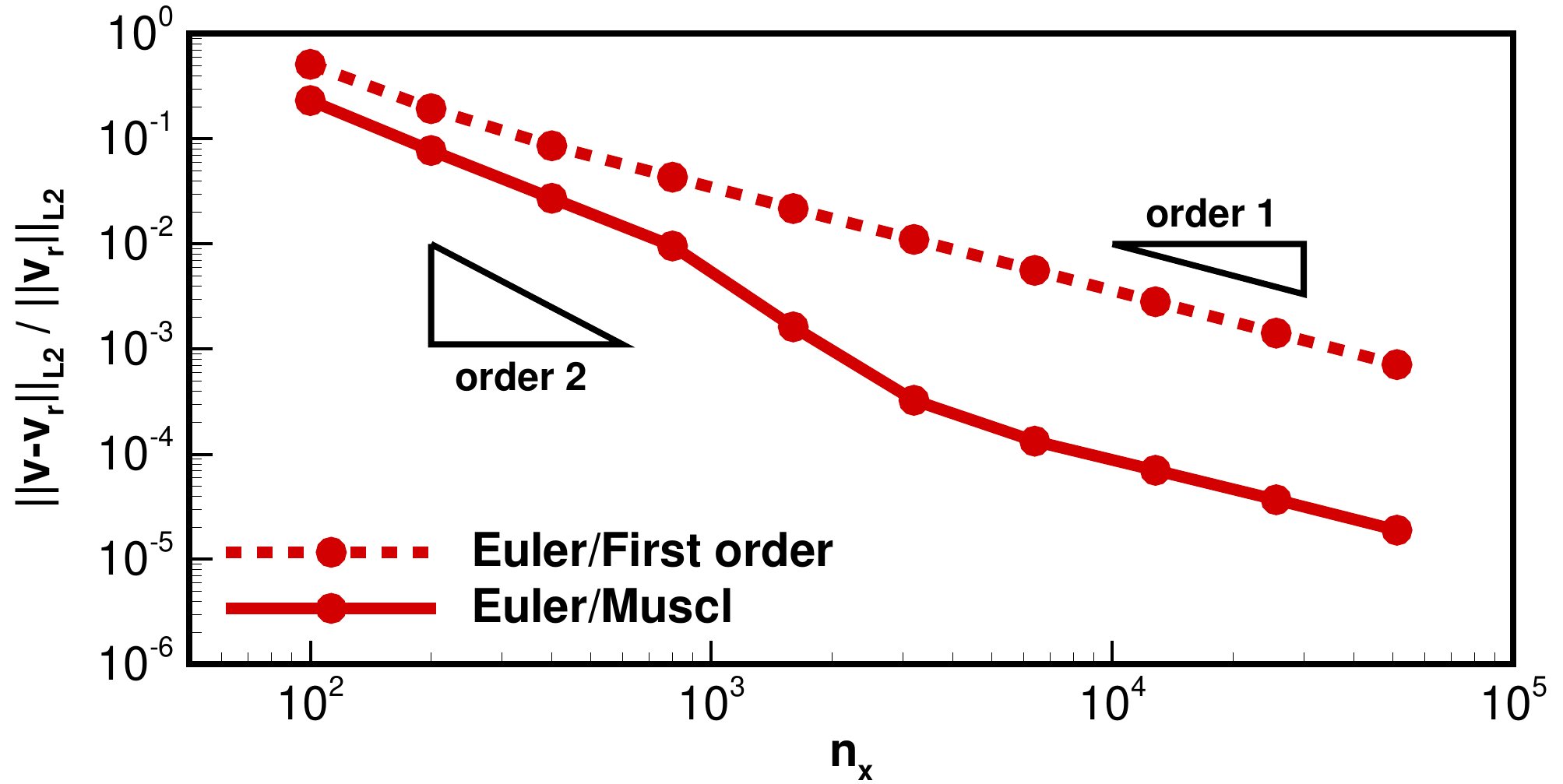}
        \par\end{centering}
    \caption{Relative error for the auxiliary velocity $\boldsymbol{v}$ in the
        $L_{2}$ norm as a function of the grid size, comparatively to the
        velocity $\boldsymbol{v_{r}}$ recomputed from $h$ and its gradient
        $\nabla h$ for the same grid size; (left) with the linearized capillary
        contribution; (right) with the full nonlinear capillary contribution.\label{fig:conv_v}}
\end{figure}

A convergence study has been made for these same parameters, considering
different grid resolution in space, with a CFL number fixed to $0.01$.
The complete results for the water elevation $h$, the fluid velocity
$\boldsymbol{u}$ and the auxiliary velocity $\boldsymbol{v}$ and
for the first grid sizes of $100$, $200$ et $400$ points for are
given in Figure (\ref{fig:1d-simulations}), for the linearized capillary
contribution version of the model as well as the full capillary contribution
version, in order to materialize the numerical solution quality. The
relative error for the water height $h$ in the $L_{2}$ norm has
been plotted in Figure (\ref{fig:conv_h}), computed at the end of
the simulation given a previous reference solution computed with $51200$
points. This has been made with both first and second order schemes
in space (without and with MUSCL reconstructions, with no limitation
as the solution is very smooth). The benefit of the MUSCL reconstruction
can be clearly noticed, especially as soon as the meshes are of medium
size, when the characteristic wavelength is meshed by more than approximately
$10$ points. However, an asymptotic convergence of 1 should be found
increasing mesh grid sizes due to the use of a first order time-stepping
scheme. But the finest mesh used of $6400$ points is not yet fine
enough to find it. It is validating partially the choice to use a
simple split explicit/implicit Euler time-stepping scheme rather than
a more sophisticated IMEX time-stepping method. Indeed, an IMEX time-stepping
scheme at second order requires more than ten times of computational
time in the present case, due to the mandatory resolution of the full
nonlinear problem, knowing that the consistency error in space is
predominant over the error in time.

The relative error for the auxiliary velocity $\boldsymbol{v}$ in
the $L_{2}$ norm as a function of the grid size has been plotted
in Figure \eqref{fig:conv_v}, comparatively to the velocity $\boldsymbol{v_{r}}$
recomputed from $h$ and its gradient $\nabla h$ for the same grid
size. The purpose is to check if the velocity field $\boldsymbol{v}$
once advected in time is still the one that carries the capillary
energy as defined by the Eq.\eqref{v}. And we can verify that this
is the case as it naturally converges with the grid size as the numerical
consistency errors and the residual error in the linear system resolution
deviate $\boldsymbol{v}$ from the ``right'' solution. But even
for very coarse meshes, the relative error is relatively low and of
course even more with MUSCL reconstructions. Also note that the relative
error is slightly more important when the full nonlinear capillary
contribution version is used.

\begin{figure}[!hb]
    \begin{centering}
        \includegraphics[width=0.5\textwidth]{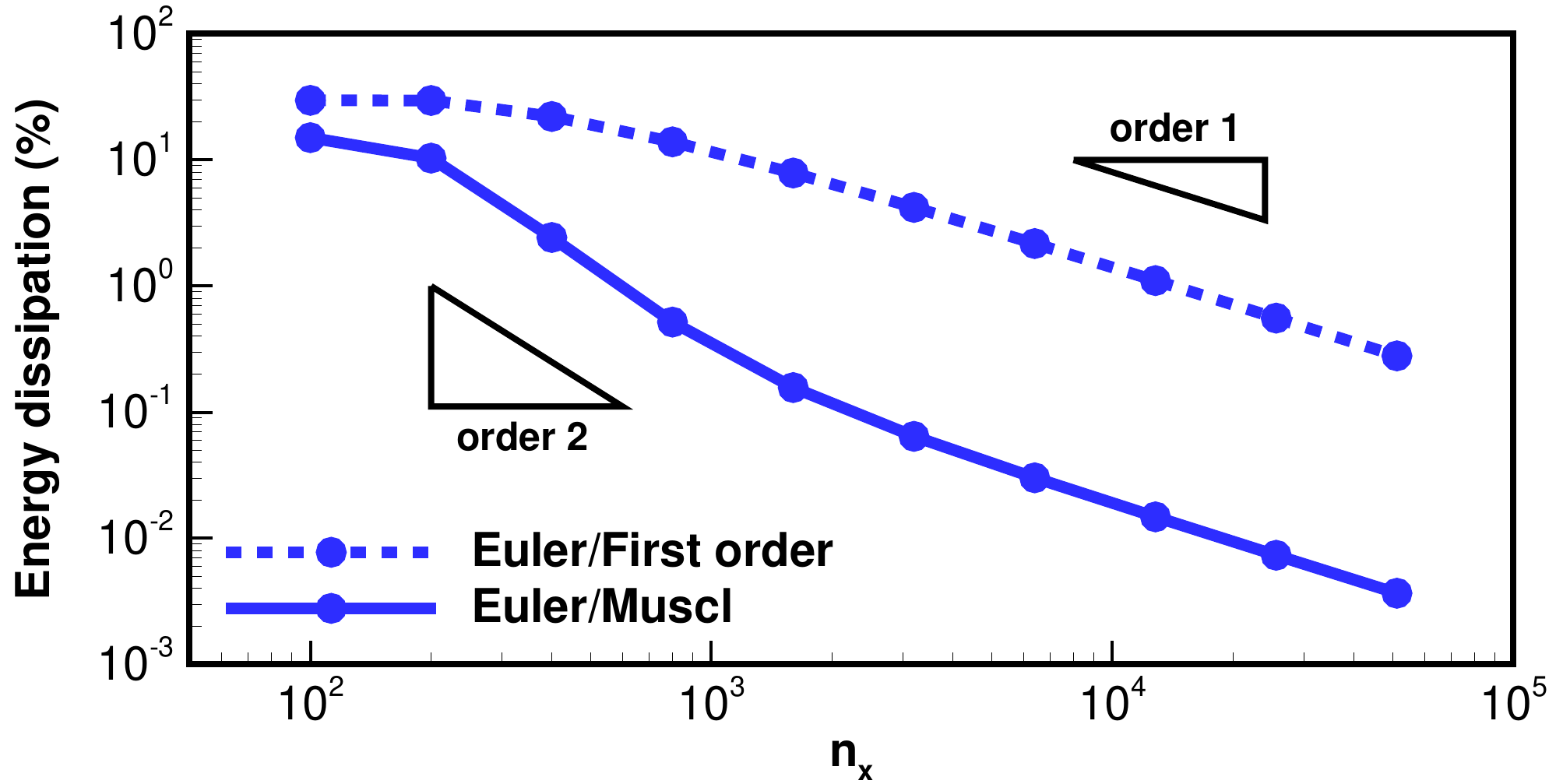}\includegraphics[width=0.5\textwidth]{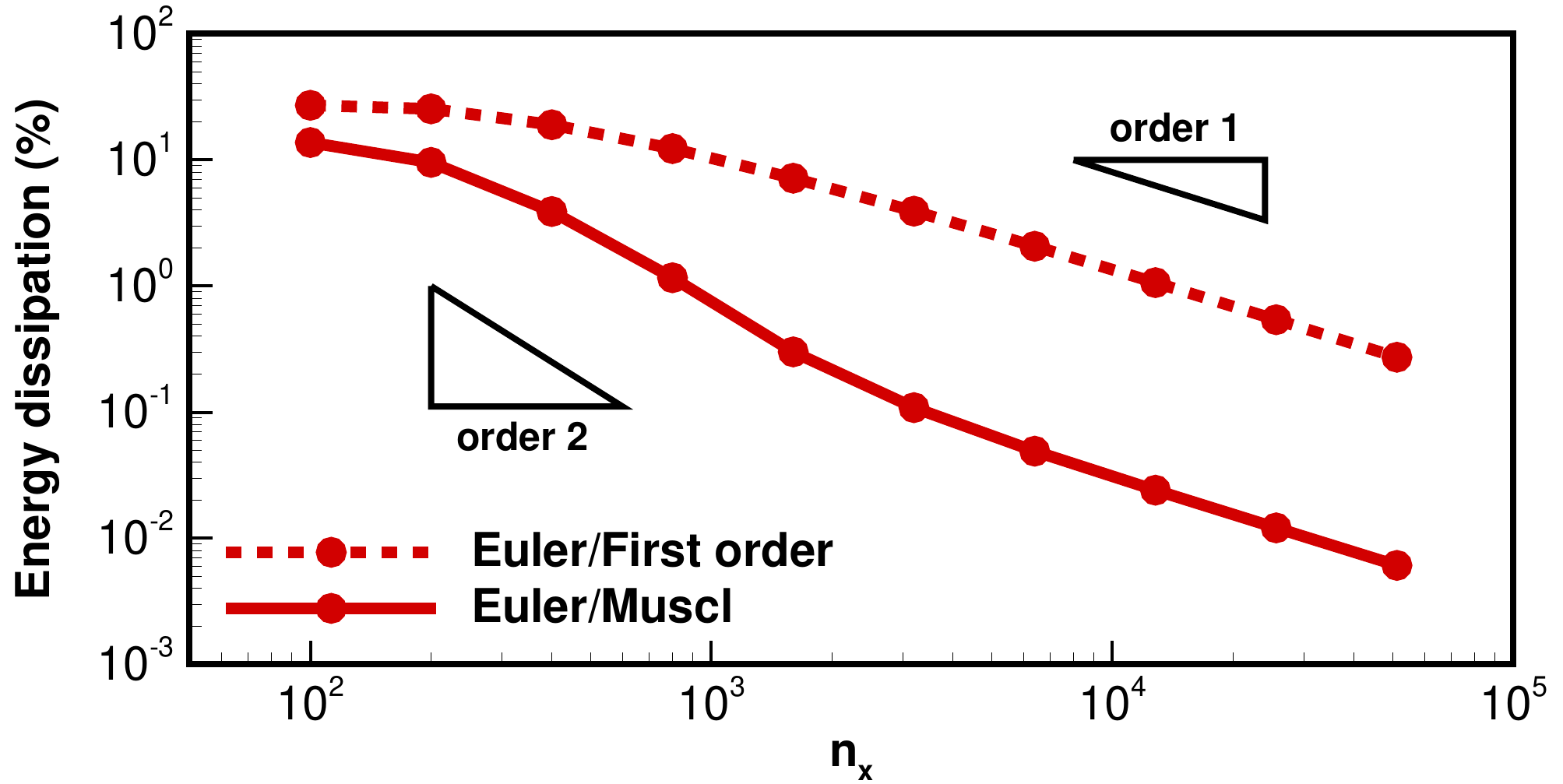}
        \par\end{centering}
    \caption{Amount of energy dissipation in percents as a function of the grid
    size, computed at the end of the simulation following the formula
    $E=gh^{2}/2+h\left\Vert \boldsymbol{u}\right\Vert ^{2}/2+h\left\Vert \boldsymbol{v}\right\Vert ^{2}/2$;
    (left) with the linearized capillary contribution; (right) with the
    full nonlinear capillary contribution.\label{fig:conv_E}}
\end{figure}
\begin{figure}[!hb]
    \centering{}\includegraphics[width=0.5\textwidth]{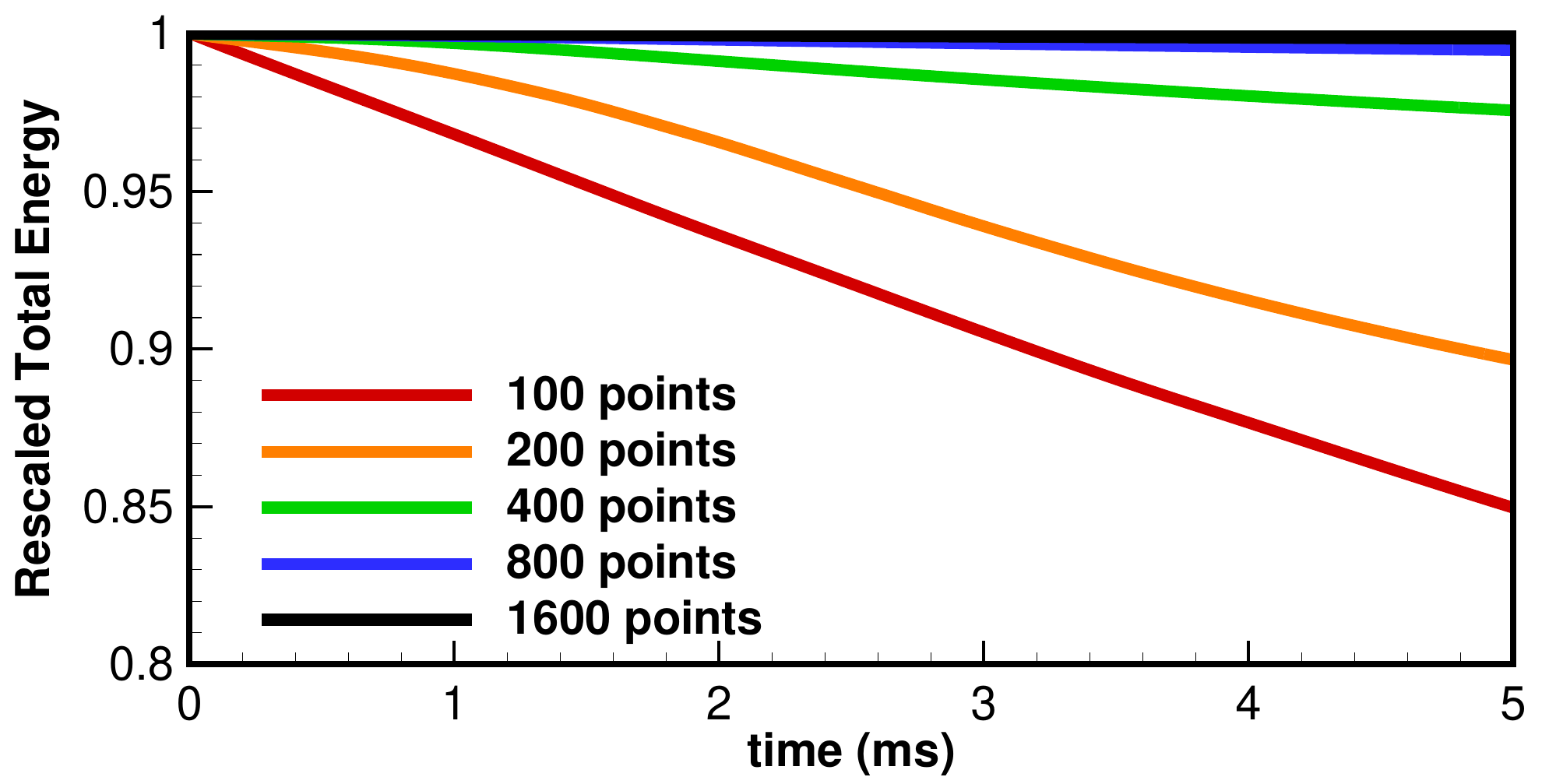}\includegraphics[width=0.5\textwidth]{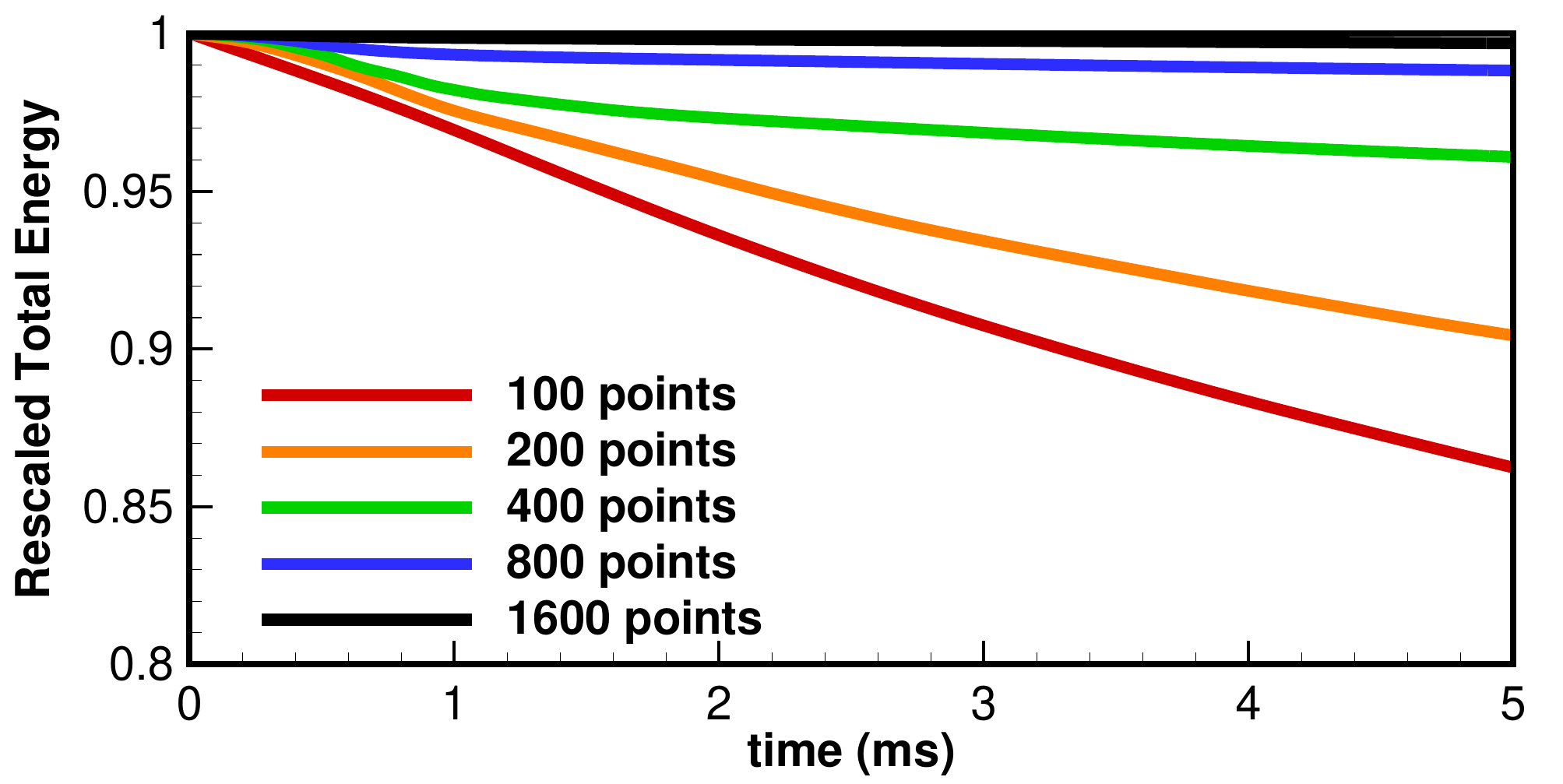}
    \caption{Evolution in time of the energy for the first grid sizes considered;
        (left) with the linearized capillary contribution; (right) with the
        full nonlinear capillary contribution.\label{fig:E-in-time}}
\end{figure}
The amount of energy dissipation in percents as a function of the
grid size is shown in Figure \eqref{fig:conv_E}. We recall that
the energy is strictly dissipated at each time step as it has been
demonstrated previously. It can be verified that this is indeed the
case in practice by looking at Figure \eqref{fig:E-in-time}. The
energy can be interpreted as an $L_{2}$ norm with the advantage to
check in one measure all the contributions in the numerical system,
rather than to check separately the convergence in a chosen norm for
the water height $h$ and the two velocities $\boldsymbol{u}$ and
$\boldsymbol{v}$. For very coarse meshes, representing few points
in the characteristic wavelength, see Figure (\ref{fig:1d-simulations}),
approximately 10\% of energy dissipation is found, which is relatively
acceptable with regard to the grid resolution used. Using MUSCL reconstruction
for finer meshes, an extra rate of convergence greater than 2 is reached
before falling to the theoretical asymptotic rate of 1 for very fine
meshes. Whereas, without MUSCL reconstructions, the convergence rate
begin at a value lower than 1, giving quickly significant differences,
before reaching asymptotically the same theoretical convergence rate
of 1 for very fine meshes. It gives finally an important order of
magnitude difference of approximately 2 when the characteristic wavelength
is sufficiently meshed with more than $10$ points.

\subsection{Two-dimensional simulations with Gaussian initial data\label{subsec:Two-dimensional-Gaussian}}

A two-dimensional version of the same previous problem (\ref{subsec:One-dimensional-Gaussian})
is now considered. The initial Gaussian-shaped deformation of a layer
of water is materialised initialising the water elevation by,

\begin{equation}
    h(x,y,t=0)=h_{0}+h_{1}e^{{\displaystyle -\frac{x^{2}+y^{2}}{2\left(b/b_{0}\right)^{2}}}}
\end{equation}

The physical parameters are the same than ones summarised in the Tab.(\ref{tab:Physical-parameters}),
as well as the space scaling with an E\"tv\"os number again chosen to
1, giving a layer of water deformation elevation $h_{1}=h_{0}=2.725$
\SI{}{\milli\metre}, and a full width at tenth of maximum again
fixed to $b=1.5\,h_{1}$. The computational domain is set to $\left[\SI{-50}{\milli\metre},\SI{50}{\milli\metre}\right]^{2}$
and the simulation time to \SI{5}{\milli\second}.

\begin{figure}[!ht]
    \begin{centering}
        \includegraphics[width=0.65\textwidth]{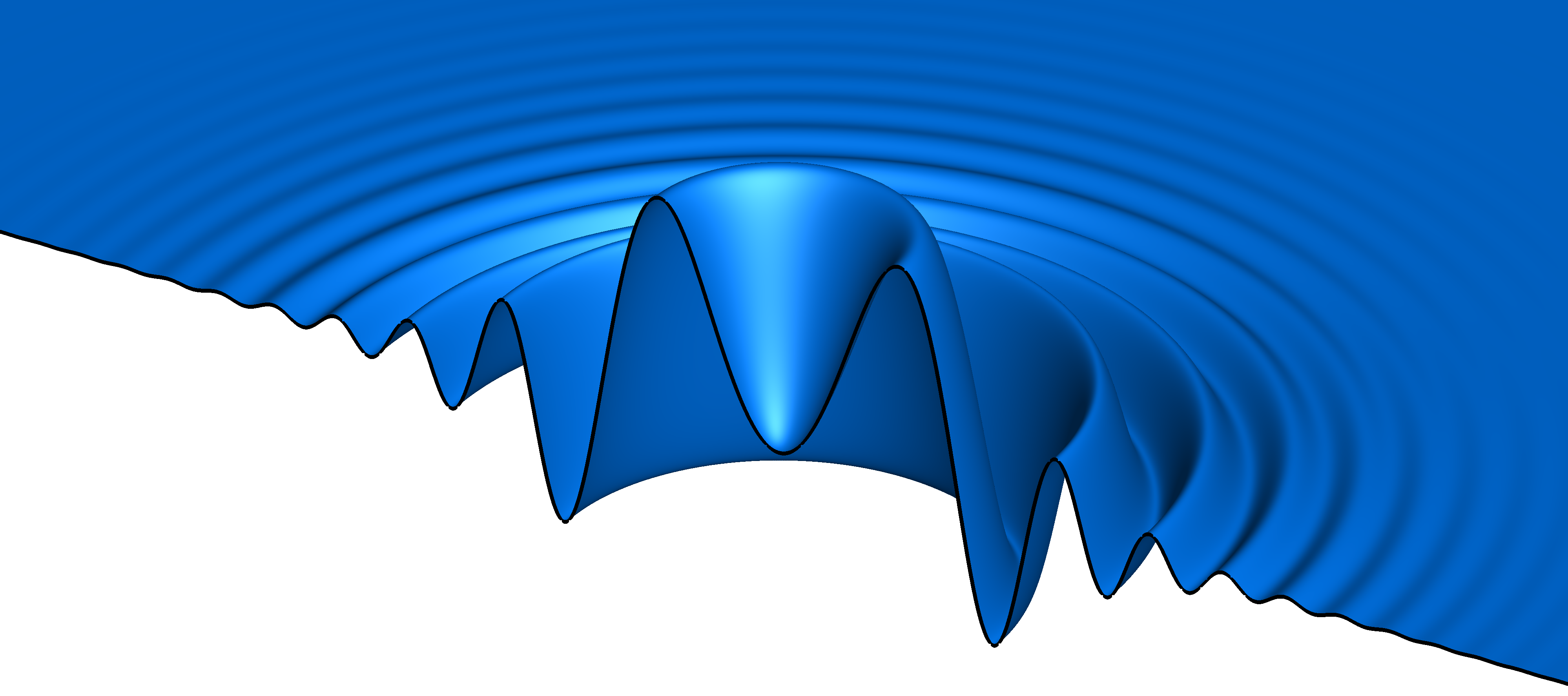}
        \par\end{centering}
    \begin{centering}
        \vspace{5mm}
        \par\end{centering}
    \begin{centering}
        \includegraphics[width=0.65\textwidth]{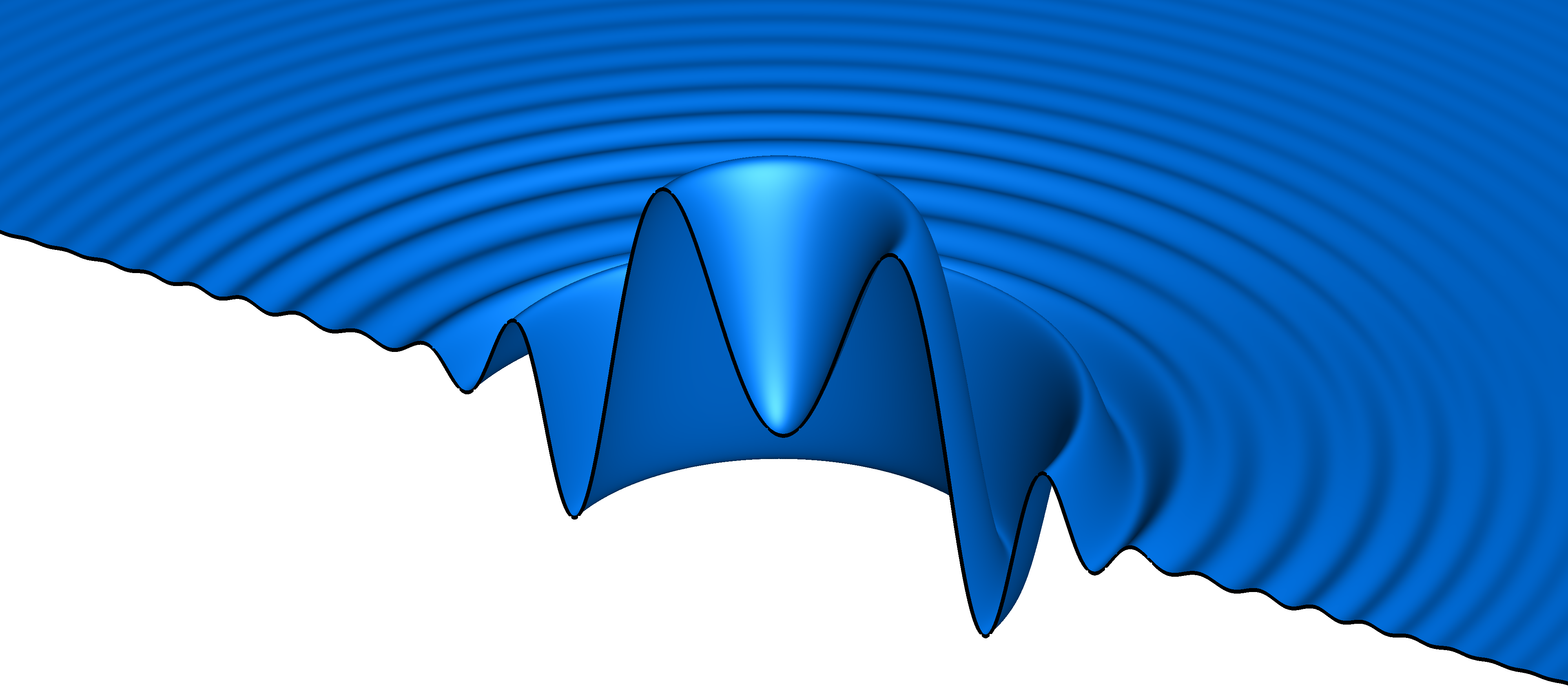}
        \par\end{centering}
    \caption{Numerical simulations of capillary-gravity waves considering a two-dimensional
        Gaussian-shaped deformation of a layer of water, using the proposed
        augmented model with parameters Eq.(\eqref{eq:quad_case},\eqref{eq:full_nl_case})
        and a $1600\times1600$ cells grid; (up) with the linearized capillary
        contribution; (down) with the full nonlinear capillary contribution.\label{fig:2d-view}}
\end{figure}
\begin{figure}[!ht]
    \begin{centering}
        \includegraphics[width=0.4\textwidth]{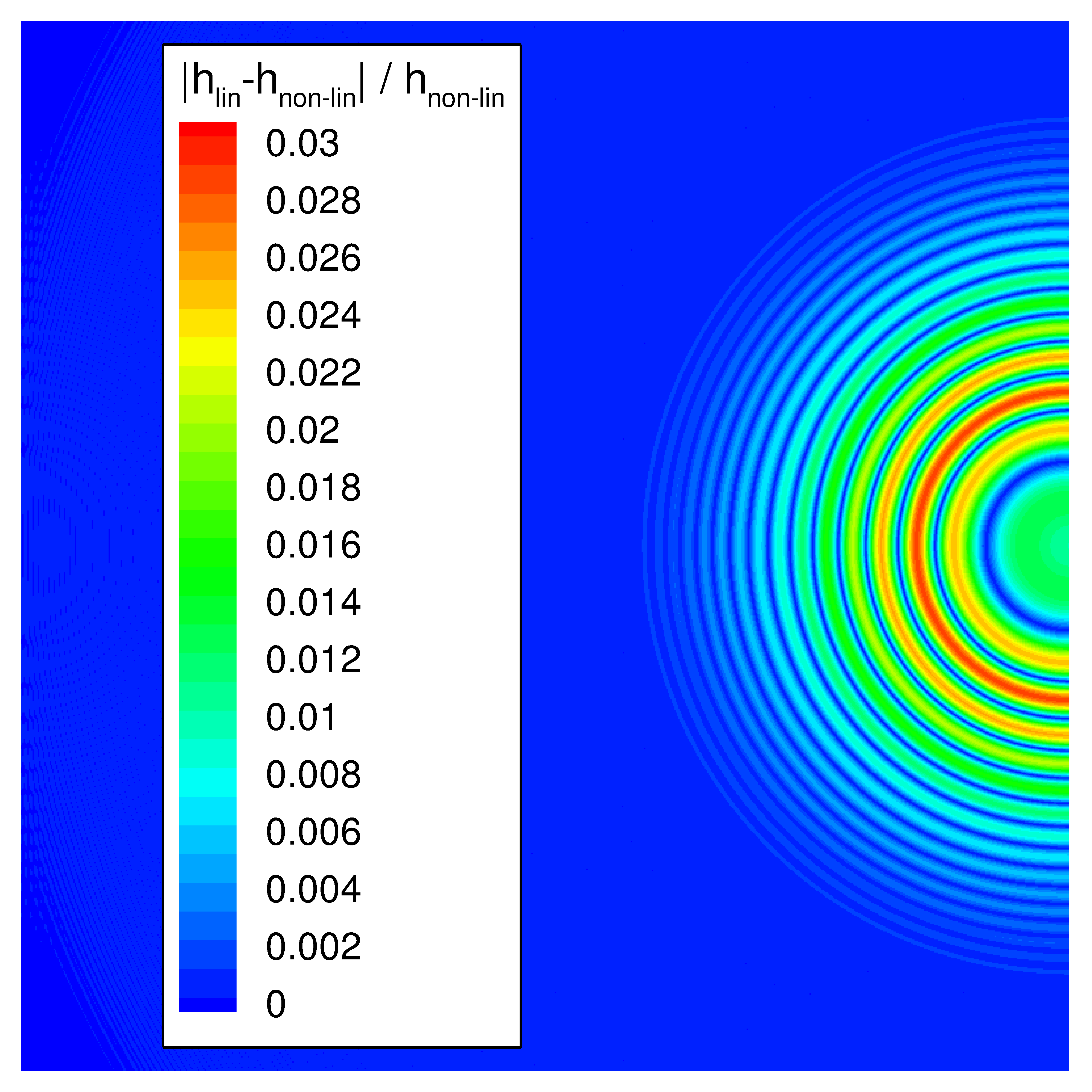}\includegraphics[width=0.4\textwidth]{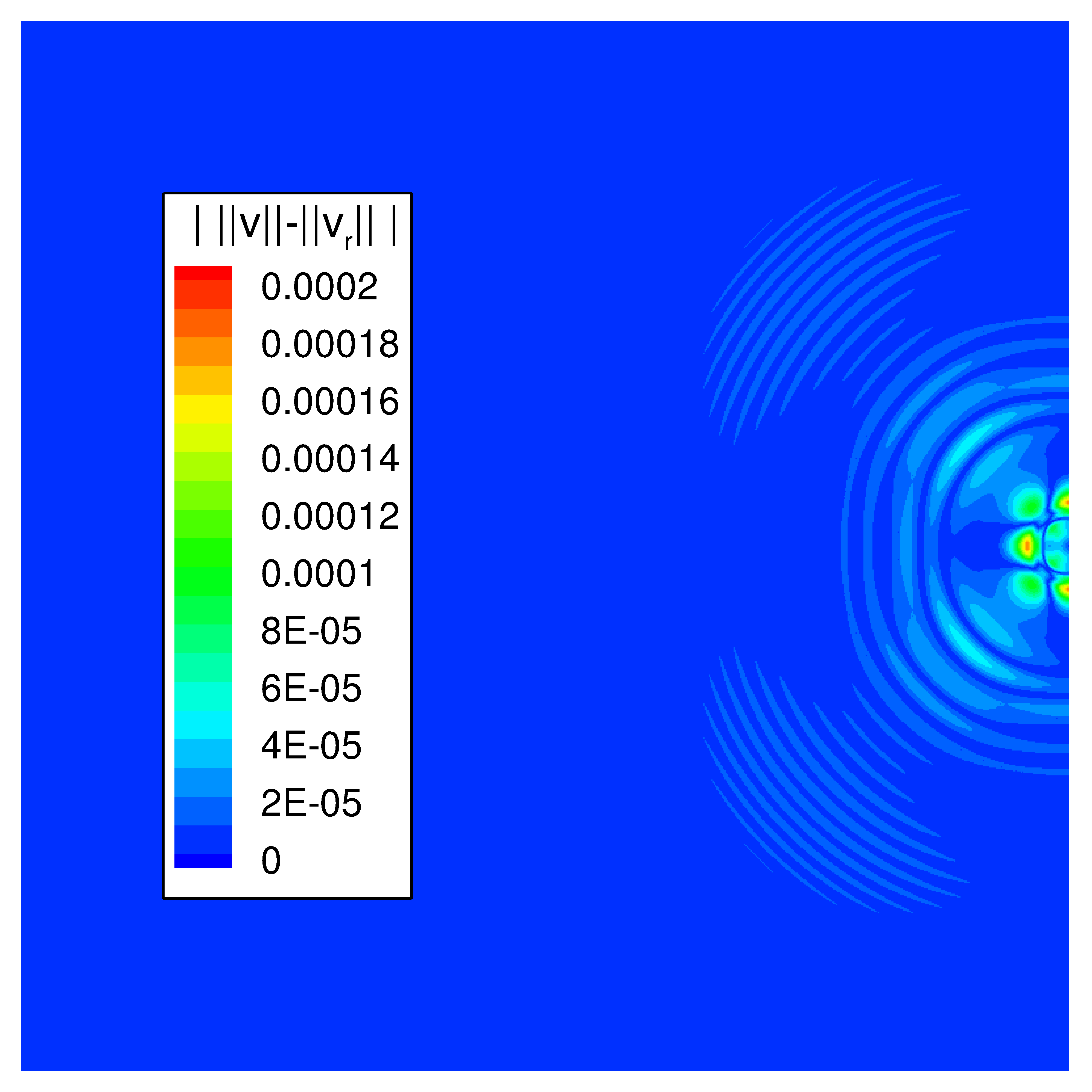}
        \par\end{centering}
    \caption{Same simulation than for the Figure \eqref{fig:2d-view} with half of
        the domain in each direction; (left) relative difference between the
        water height for the two models; (right) absolute difference between
        the auxiliary velocity magnitude $\left\Vert \boldsymbol{v}\right\Vert $
        and the velocity magnitude $\left\Vert \boldsymbol{v_{r}}\right\Vert $
        recomputed from $h$ and its gradient $\nabla h$ for the same grid
        size.\label{fig:2d_hdiff_1600}}
\end{figure}
\begin{figure}[!ht]
    \begin{centering}
        \includegraphics[width=0.5\textwidth]{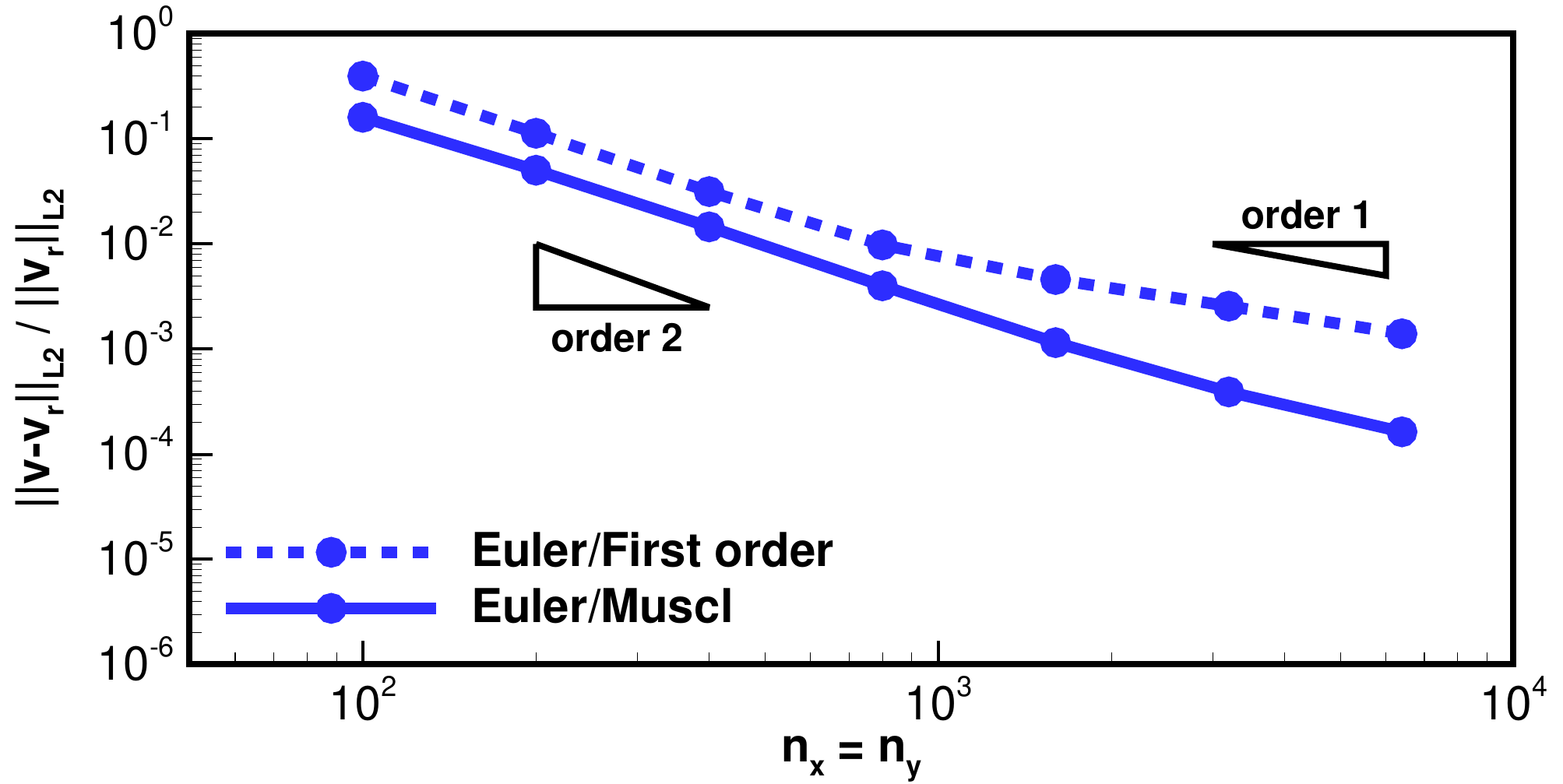}\includegraphics[width=0.5\textwidth]{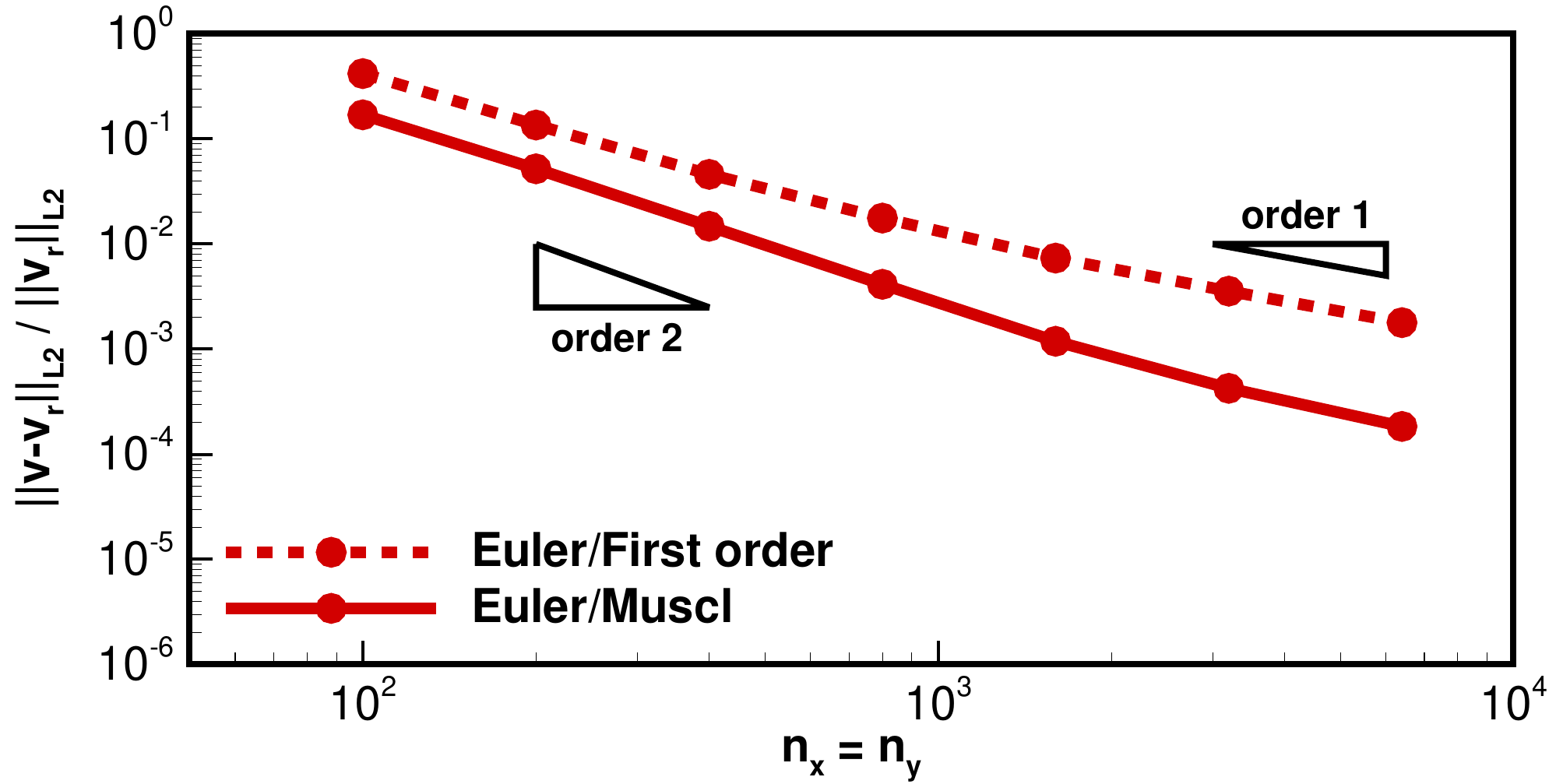}
        \par\end{centering}
    \caption{Relative error for the auxiliary velocity magnitude $\left\Vert \boldsymbol{v}\right\Vert$
        in the $L_{2}$ norm as a function of the grid size, comparatively
        to the same velocity magnitude $\left\Vert \boldsymbol{v_{r}}\right\Vert$
        recomputed from $h$ and its gradient $\nabla h$ for the same grid
        size; (left) with the linearized capillary contribution; (right) with
        the full nonlinear capillary contribution.\label{fig:conv_v_2d}}
\end{figure}
\begin{figure}[!ht]
    \begin{centering}
        \includegraphics[width=0.5\textwidth]{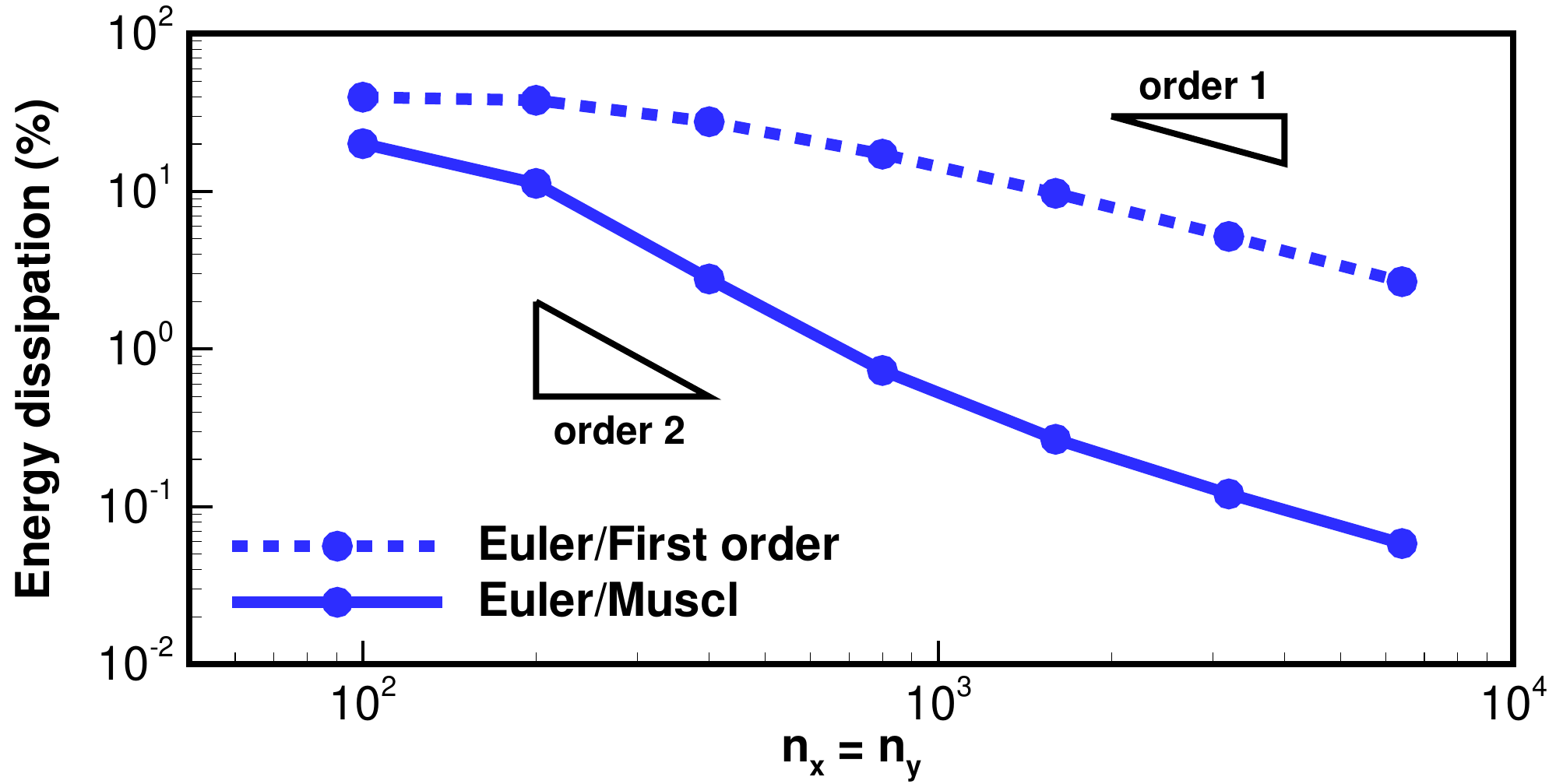}\includegraphics[width=0.5\textwidth]{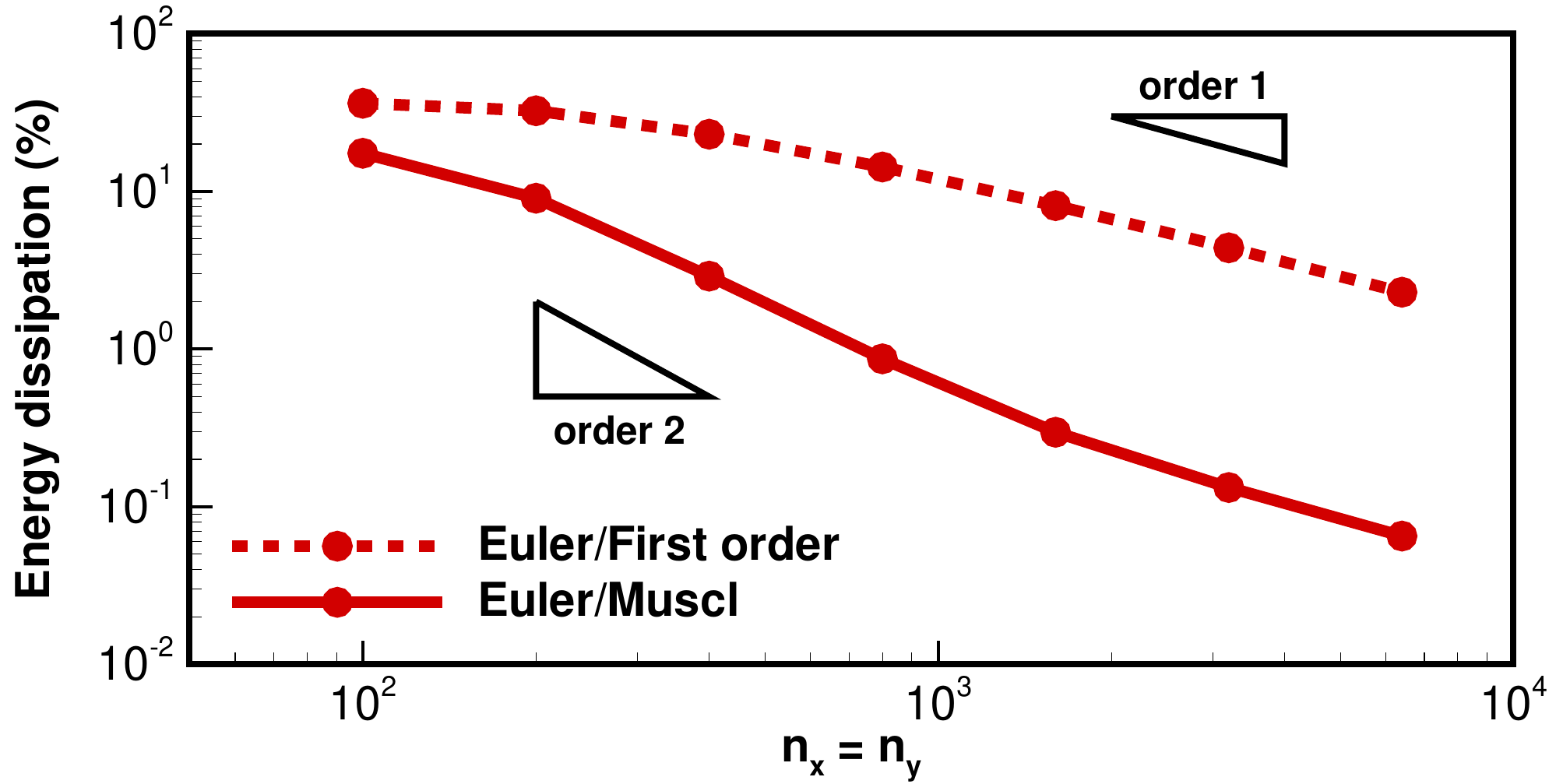}
        \par\end{centering}
    \caption{Amount of energy dissipation in percents as a function of the grid
    size, computed at the end of the simulation following the formula
    $E=gh^{2}/2+h\left\Vert \boldsymbol{u}\right\Vert ^{2}/2+h\left\Vert \boldsymbol{v}\right\Vert ^{2}/2$;
    (left) with the linearized capillary contribution; (right) with the
    full nonlinear capillary contribution.\label{fig:E_conv_2d}}
\end{figure}
It can be observed in Figure (\ref{fig:2d-view}) the axisymmetric
propagation of capillary-gravity waves using the proposed augmented
shallow-water model Eq.(\ref{def_U_F},\ref{eq:def_M}) with formulas
Eq.(\ref{eq:quad_case},\ref{eq:full_nl_case}). The initial peak
collapse on his own weight and generate a train of capillary waves.
The difference between the two models for the water height $h$ is
plotted in Figure (\ref{fig:2d_hdiff_1600}). It can be observed
a maximum difference of 3\%. If the same phenomena as in the one-dimensional
case can be observed, the initial chosen shape of the Gaussian deformation
of the water layer gives lower differences since axisymmetry reduces
the capillary waves speed difference and the phase shift generated.

The two-dimensional version behaves like the one-dimensional one,
as it can be seen in Figure (\ref{fig:conv_v_2d}) and in Figure (\ref{fig:E_conv_2d}).
In the same way, the relative error for the auxiliary velocity magnitude
$\left\Vert \boldsymbol{v}\right\Vert$ comparatively to the recomputed
one $\left\Vert \boldsymbol{v_{r}}\right\Vert$ converges with the grid
size, and is relatively low even for coarse meshes, with a lower error
using MUSCL reconstructions. In addition, it is also shown in Figure \eqref{fig:2d_hdiff_1600}
a map of the absolute difference showing no particular region with
a big peak. The amounts of energy dissipation are only slightly higher
than for the one-dimensional case. Again, the MUSCL reconstructions
provide an extra rate of convergence near 2 for medium grid sizes
before falling asymptotically to the theoretical rate of 1, which finally
gives approximately two orders of magnitude of difference if not using
it.

This validate the application of the proposed augmented shallow water
model Eq.(\ref{def_U_F},\ref{eq:def_M}) for both one- and two-dimensional
cases, with very similar numerical convergence behaviour.

\subsection{Liquid film falling on a vertical wall}

The purpose is now to apply our extended shallow water model with
surface tension to a more realistic case of simulation of a liquid
film falling on a vertical wall. We are able to compare the resulting
numerical solutions with the solutions computed from the \texttt{Auto07p}
software (\cite{auto07p}). The model used is the one developed by
G. Richard and al. \cite{RiGiRuVi} written here in non-dimensional form,

\begin{equation}
    \begin{array}{lll}
         &  & \partial_{t}h+\partial_{x}\left(hu\right)  =  0                                                                                                                                                                                                                                                                                             \\
        \\
         &  & \partial_{t}\left(hu\right)+\partial_{x}\left(hu^{2}+{\displaystyle \frac{2}{225}\lambda^{2}h^{5}+\frac{h^{2}\cos\theta}{2\mathrm{Fr}^{2}}}\right)  =  {\displaystyle \frac{1}{\mathrm{Re}}\left(\lambda h-\frac{3u}{h}\right)+\frac{9}{2\mathrm{Re}}\partial_{x}\left(h\partial_{x}u\right)+\frac{1}{\mathrm{We}}h\partial_{x}\mathcal{K}}
    \end{array}\label{eq:model_falling_film}
\end{equation}
where $\mathcal{K}$ is the curvature, and its expression is $\partial_{xx}h$
in the linearized case and $\partial_{x}\left(\left.\partial_{x}h\right/\sqrt{1+\partial_{x}h^{2}}\right)$
in the full nonlinear case. Scales are the film thickness $h_{N}$
and the Nusselt speed $u_{N}=\left.gh_{N}^{2}\sin\theta\right/\left(3\nu\right)$
. Dimensionless numbers are the Reynolds number $\mathrm{Re}=\left.h_{N}u_{N}\right/\nu=\left.gh_{N}^{3}\sin\theta\right/\left(3\nu^{2}\right)$,
the Froude number $\mathrm{Fr}=\left.u_{N}\right/\sqrt{gh_{N}}$,
the Weber number $\mathrm{We}=\left.\rho h_{N}u_{N}^{2}\right/\sigma$
and $\lambda=\left.\mathrm{Re}\sin\theta\right/\mathrm{Fr^{2}}=3$.
The Kapitza dimensionless number $\mathrm{Ka}={\displaystyle \left(\sigma/\rho\right)\left(g\sin\theta\right)^{-1/3}\nu^{-4/3}}$
can be added such that $\mathrm{Ka}={\displaystyle \mathrm{Re}^{4/3}\mathrm{Fr}^{2/3}\mathrm{We}^{-1}}$.

We have constructed solitary wave solutions to \eqref{eq:model_falling_film}
using the \texttt{Auto07p} software with the constraint of a constant
averaged thickness $\left\langle h\right\rangle =1$. The system of
partial differential equations simplifies into ordinary differential
equations in the moving frame of reference $\xi=x-ct$, where $c$
refers to the phase speed of the waves. Travelling wave solutions to \eqref{eq:model_falling_film}
correspond to limit cycles of the resulting autonomous dynamical system in a phase space of dimension three spanned by the thickness $h$ and its first and second derivatives. Limit cycles are found by continuation starting from a Hopf bifurcation of a fixed point corresponding to the uniform-film solution $h=1$. Solitary waves are next obtained by increasing the period of the limit cycles.
The parameters retained are
a vertical wall $\theta=\pi/2$, a Reynolds number $\mathrm{Re}=80$,
a Kapitza number $\mathrm{Ka}=1000$ and a length $L=400\:h_{N}$.
The liquid corresponds to $\nu=0.9310\,10^{-6}\,m^{2}.s^{-1}$, $\rho=994.3\,kg$
and $\sigma=0.019322\,N.m^{-1}$ and the gravitational acceleration
to $g=9.81\,m.s^{-2}$. The Nusselt film thickness is $h_{N}=0.27659\,mm$
so that the length of the numerical domain is $L=110.636\,mm$. The
grid comprises $8000$ mesh points with spatial adaptation. The spatial discretization uses
the method of orthogonal collocation using piecewise polynomials with four collocation points per
mesh intervals (2000 mesh intervals). The collocation points are placed to equidistribute the
local discretization error in the three-dimensional phase space, which ensures refinement of the mesh at locations of steep gradients. \texttt{Auto07p} uses a predictor-corrector algorithm based on a Keller's pseudo-arc length continuation method that enables to detect bifurcations and folds. Details of the algorithm can be found online at \texttt{http://indy.cs.concordia.ca/auto/}. Convergence was checked by varying the error tolerances and the number of mesh cells.

Numerical simulations have been carried out using the extended shallow
water model system (\ref{def_U_F},\ref{eq:def_M}) defined in exactly
the same way as in (\ref{eq:quad_case},\ref{eq:full_nl_case}). The only differences
are a different pressure term in the hyperbolic system changing the
Rusanov flux and the introduction of additional source terms treated
implicitly simply introducing them in the original linear system arising
from the surface tension treatment. A periodic boundary condition
and an initial sinusoidal deformation of the film liquid at rest,
taking care to a constant mass flow rate, are prescribed. As we are
only interested in whether the model is capable of reproducing the
solutions given by the \texttt{Auto07p} software, only simulations
with a very fine mesh size of $25600$ cells are presented here. We
can observe an almost perfect agreement between the solutions validating
the good behaviour of the proposed extended model regarding a produced
exact solution from the resolution of ordinary differential equations.
It is not presented here but the liquid film deformations are perfectly
stable in time.

\begin{figure}[!ht]
    \begin{centering}
        \includegraphics[width=0.5\textwidth]{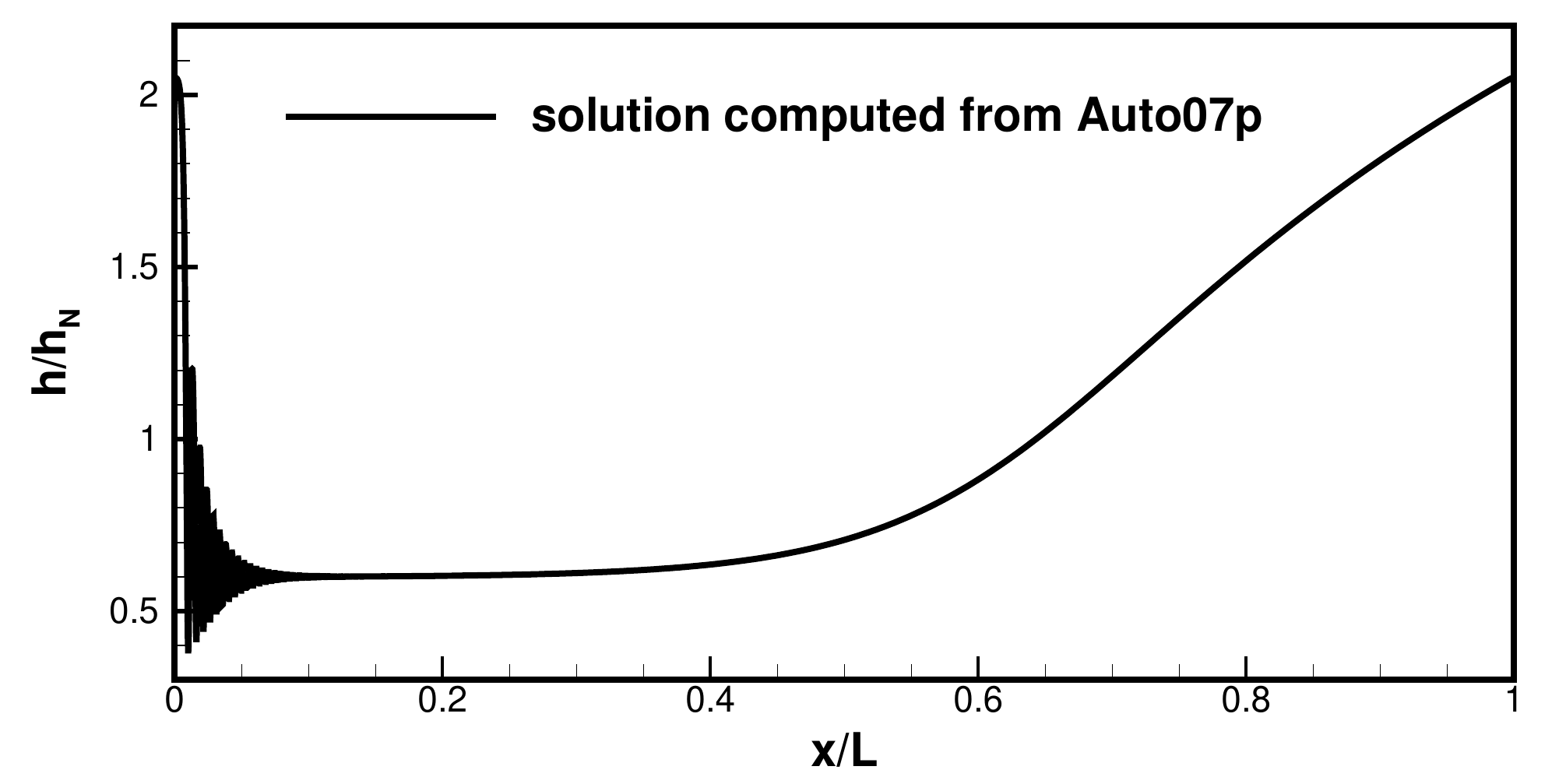}\includegraphics[width=0.5\textwidth]{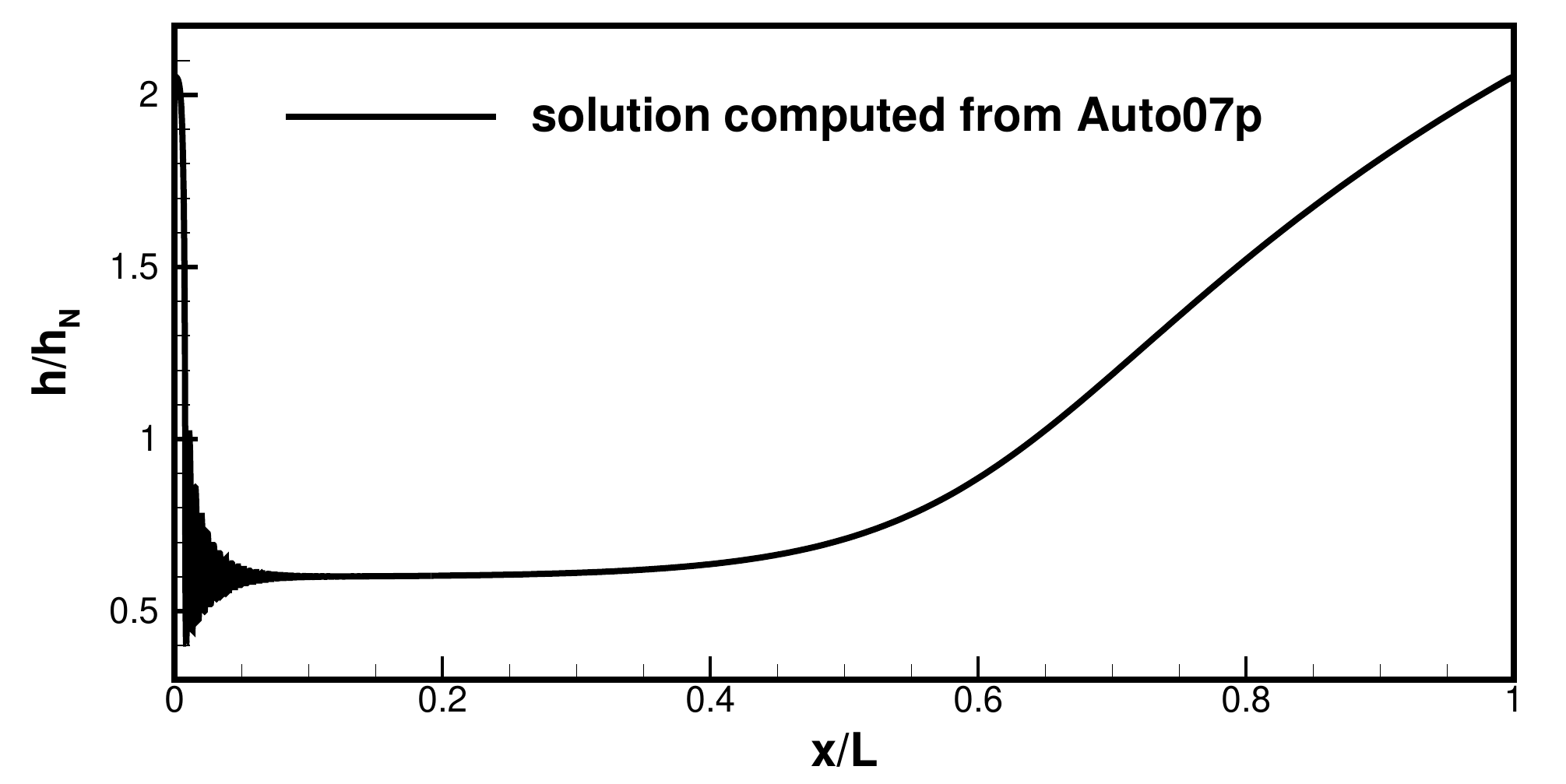}
        \par\end{centering}
    \begin{centering}
        \includegraphics[width=0.5\textwidth]{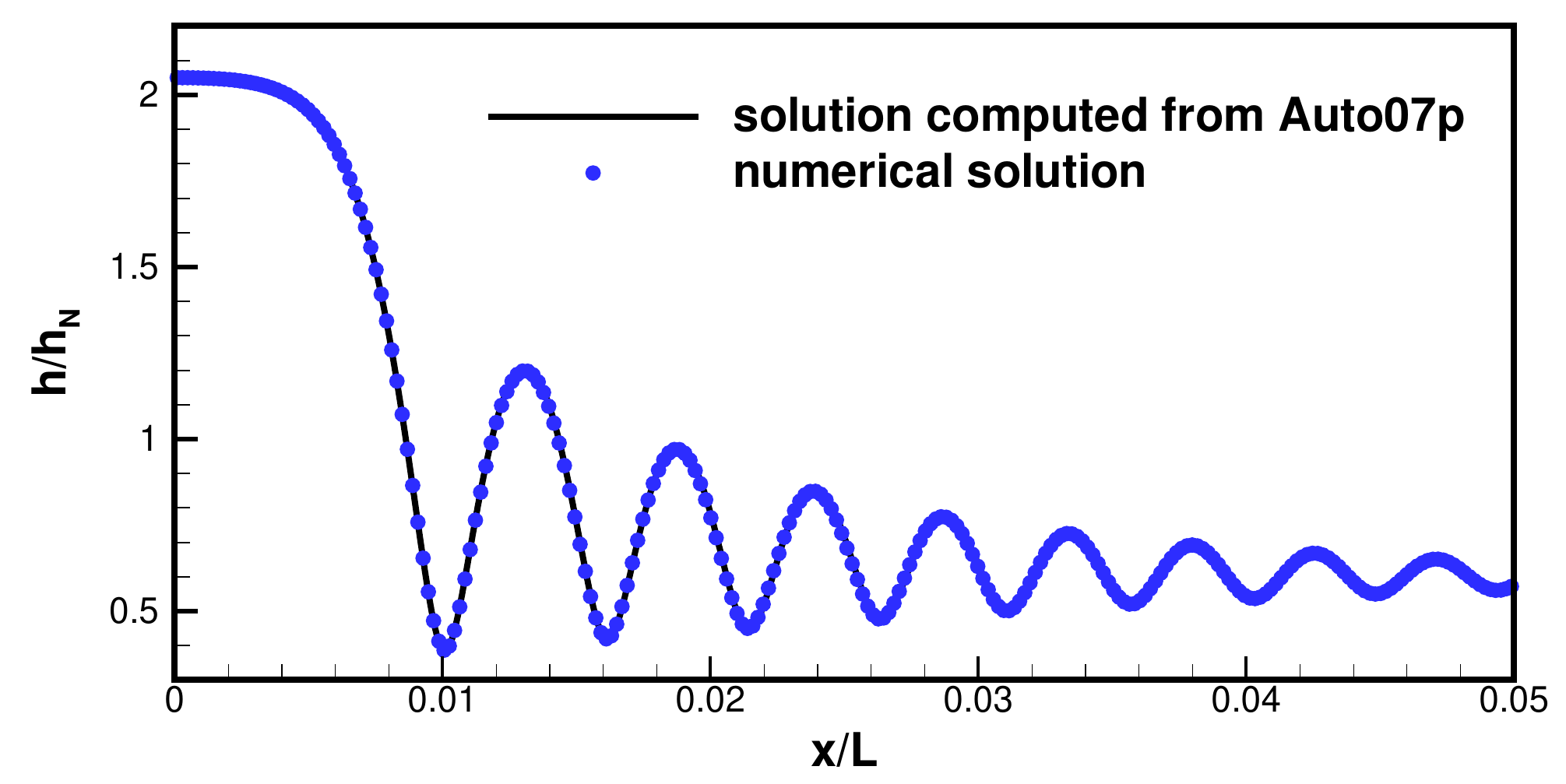}\includegraphics[width=0.5\textwidth]{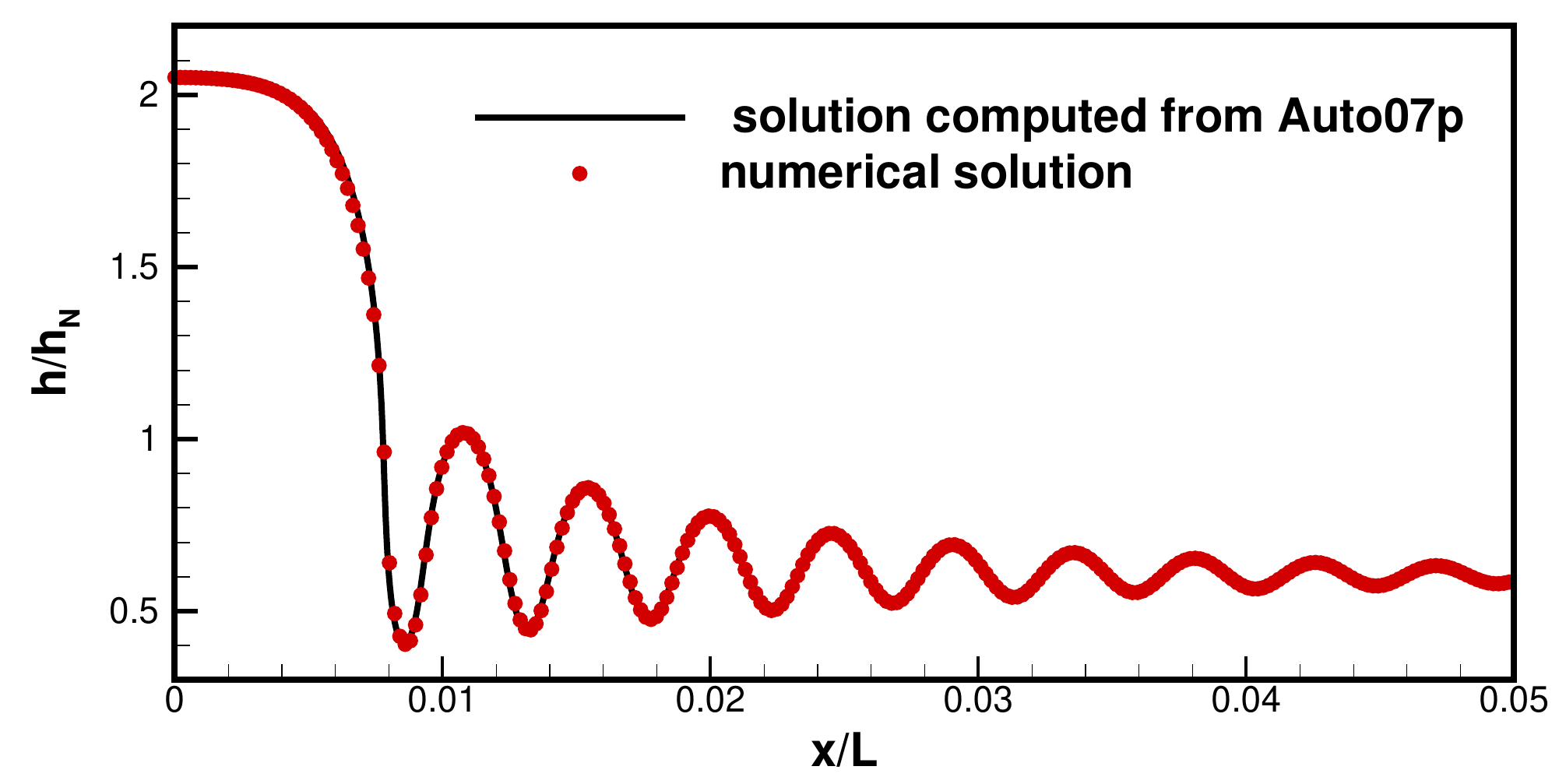}
        \par\end{centering}
    \caption{Liquid film falling on a vertical wall at $\mathrm{Re}=80$ and $\mathrm{Ka}=1000$;
        (top) Complete solutions computed with the \texttt{Auto07p} software;
        (bottom) Comparison between the numerical solutions (one point per
        ten considering $25600$ cells) and the solution computed with the
        \texttt{Auto07p} software with a focus on the capillary ripples; (left)
        with linearized curvature $\mathcal{K}$ in (\ref{eq:model_falling_film});
        (right) with full nonlinear curvature $\mathcal{K}$ in (\ref{eq:model_falling_film}).}
\end{figure}

\section{Conclusion and Perspectives}

In this paper, we have introduced a new extended version of the shallow
water equations with surface tension which may be decomposed in two parts:
a conservative hyperbolic part and a second order derivative part which
is skew-symmetric with respect to the $L^{2}$ scalar product.
This extended form is suitable for an appropriate splitting method which
allows for large gradients of fluid height. The formulation is valid for any non
linear form of the capillary energy as a functional of $\|\nabla h\|$.
In the case of deep gradient of the free surface the new formulation make
possible to deal with complete nonlinear capillary models.

The formulation allows to deal with the capillary terms as a semi-linear
skew symmetric problem, and thus associate a semi implicit resolution
of the capillary terms, relaxing the time step restriction $\Delta t<Ch^{2}$
in the original formulation \cite{NoVi} and \cite{BCNV} which
are restricted to quadratic forms of the capillary energy.
We expect that this property (semi linear implicit treatment of capillary
terms) will ensure some robustness and will be successful in the extension
of this methodology to wetting problems as in the work of
J. Lallement, P. Trontin, C. Laurent and P. Villedieu published in
\cite{LaTrLaVi} where ad hoc extension with some disjunction pressure
are proposed.
We recently developed various 3 or 4 equations models (see \cite{RiRuVi}), for thin film models which are extension of the Korteweg system introduced here.
We expect that our formalism also apply to such systems. Nevertheless
the question of non-linear stability remains to study. In such case
the natural extension of our proofs relies on entropy estimate (related
to the so called enstrophy introduced in \cite{RiRuVi}).

Finally it will be interesting to study possible applications to the
solution of Benney type equations obtained as relaxed system of the
Shallow water system with friction source terms such as those studied
in \cite{RiGiRuVi}. This reads, starting from a simplified 2D version,
\[
    \left\{ \begin{array}{cc}
        \partial_{t}h+{\rm div}\left(h\mathbf{u}\right)=0                                                                                                                                                                                                                                                          & (i)  \\
        \partial_{t}\left(h\mathbf{u}\right)+{\rm div}\left(h\mathbf{u}\otimes\mathbf{u}\right)+\nabla P=\frac{1}{\varepsilon Re}\left(\boldsymbol{\lambda}h-\frac{3\mathbf{u}}{h}\right)-{\rm div}\left(\nabla h\otimes\nabla_{\mathbf{p}}E\right)+\nabla\left(h{\rm div}\left(\nabla_{\mathbf{p}}E\right)\right) & (ii)
    \end{array}\right.
\]
gives after relaxation when $\varepsilon\rightarrow0$
\[
    \mathbf{u}=\boldsymbol{\lambda}\frac{h^{2}}{3}-\varepsilon\frac{h}{3}Re\left(-h^{2}{\rm div}\left(\boldsymbol{\lambda}\frac{h^{3}}{3}\right)+{\rm div}\left(h\boldsymbol{\lambda}\frac{h^{2}}{3}\otimes \boldsymbol{\lambda}\frac{h^{2}}{3}\right)+\nabla P+{\rm div}\left(\nabla h\otimes\nabla_{\mathbf{p}}E\right)-\nabla\left(h{\rm div}\left(\nabla_{\mathbf{p}}E\right)\right)\right)
\]
and
\[
    \partial_{t}h+{\rm div}\left(h\mathbf{u}\right)=0.
\]

\section{Appendix}
In this section, we present the Proof of Lemma  \ref{NLaugmented}.

\noindent \textbf{Part i) and iii)} \textit{Equation satisfied by
    $\boldsymbol{v}$.} Let us first recall that $\boldsymbol{v}=\alpha(q^{2})\sqrt{\frac{\sigma(h)}{h}}\nabla h$
and therefore
\[
    h\boldsymbol{v}=\alpha(q^{2}) \, \sqrt{\sigma(h)} \, h^{3/2} \, \frac{\nabla h}{h}:=\alpha(q^{2})G(h)\boldsymbol{a},
\]
with $G(h)=\sqrt{\sigma(h)}h^{3/2}$ and $\boldsymbol{a}=\nabla(\log(h))$.
In order to write an evolution equation on $h\boldsymbol{v}$, the
first step is to calculate evolution equations on ${\bf a}$, $G(h)$
and $\alpha(q^{2})$. For that purpose, we consider the mass conservation
law written as
\begin{equation}
    {\displaystyle \partial_{t}h+{\bf u}^{t}\nabla h+h\,{\rm div}({\bf u})=0.}\label{eq:mass}
\end{equation}
By dividing \eqref{eq:mass} by $h$ and differentiating with respect
to $x_{i},i=1,2$, one finds
\begin{equation}
    {\displaystyle \partial_{t}{\bf a}+\nabla({\bf u}^{t}{\bf a})+\nabla(\rm{div}({\bf u}))=0.\label{eqa}}
\end{equation}
By multiplying \eqref{eq:mass} by $G'(h)$, one finds
\begin{equation}
    {\displaystyle \partial_{t}G(h)+{\bf u}^{t}\nabla G(h)=-h G'(h)\, \mathrm{div}(\mathbf{\bf u}),\quad {G}'(h)=\frac{\sigma'(h)h^{3/2}}{2\sqrt{\sigma(h)}}+\frac{3}{2}\sqrt{\sigma(h)\,h}.\label{h32}}
\end{equation}
From \eqref{eq:mass}, we find that $\nabla h$ satisfies
\begin{equation}
    \partial_{t}\nabla h+(\boldsymbol{u}^{t}\nabla)\nabla h=-{\rm div}(h\nabla\boldsymbol{u}^{t})-\nabla h\,{\rm div}\boldsymbol{u}.\label{eq:nablah}
\end{equation}
The derivatives of $q=\|\nabla h\|$ are given by
\[
    q \, \partial_{t}q=(\nabla h)^{t}\partial_{t}\nabla h,\quad q \, \partial_{i}q=(\nabla h)^{t}\partial_{i}\nabla h,\quad i=1,2.
\]
This allows to calculate the equation on $\alpha(q^{2})$. Indeed,
we can write:
\[
    \begin{array}{ccl}
        {\displaystyle \partial_{t}\alpha(q^{2})+\boldsymbol{u}^{t}\nabla\alpha(q^{2})} & = & {\displaystyle \alpha'(q^{2})\left(\partial_{t}q^{2}+({\bf u}^{t}\nabla)q^{2}\right)}                      \\
                                                                                        & = & {\displaystyle 2\alpha'(q^{2})(\nabla h)^{t}\left(\partial_{t}\nabla h+({\bf u}^{t}\nabla)\nabla h\right)} \\
    \end{array}
\]
By substituting the value of $\partial_{t}\nabla h$ given by \eqref{eq:nablah}
into the former equation, one finds
\begin{equation}
    \partial_{t}\alpha(q^{2})+\boldsymbol{u}^{t}\nabla\alpha(q^{2})=-2\alpha'(q^{2})\left((\nabla h)^{t}{\rm div}(h\nabla\boldsymbol{u}^{t})+q^{2}{\rm div}\boldsymbol{u}\right).\label{q}
\end{equation}
Finally, by using the fact that $h{\bf v}=\alpha(q^{2})G(h){\bf a}$,
one finds that the advective term ${\rm div}(h{\bf v}\otimes{\bf u})$
is given by

\[
    {\displaystyle {\rm div}\left(\alpha(q^{2}) G(h)\,{\bf a}\otimes\boldsymbol{u}\right)=%
    \alpha(q^{2}) G(h)\left(({\bf u}^{t}\nabla){\bf a}+{\rm div}({\bf u}){\bf a}\right)+(({\bf u}^{t}\nabla)(\alpha(q^{2})\Margo{G}(h))){\bf a}.}
\]

We can now calculate the equation satisfied by $\boldsymbol{v}$ using
formula \eqref{eqa}--\eqref{q}. More precisely we have {\setlength{\arraycolsep}{1pt}
        \begin{eqnarray*}
            \partial_{t}(h\boldsymbol{v})+{\rm div}(h\boldsymbol{v}\otimes\boldsymbol{u}) & = & \partial_{t}\left(\alpha(q^{2}) G(h){\bf a}\right)+{\rm div}\left(\alpha(q^{2}) G(h){\bf a}\otimes{\bf u}\right)\\[2mm]
            & = & \alpha(q^{2})G(h)\left((\partial_{t}+{\bf u}^{t}\nabla){\bf a}+{\rm div}({\bf u}){\bf a}\right)+\left((\partial_{t}+{\bf u}^{t}\nabla)(\alpha(q^{2})G(h))\right){\bf a}\\[2mm]
            & = & \alpha(q^{2}) G(h)\left(({\bf u}^{t}\nabla){\bf a}+{\rm div}({\bf u}){\bf a}-\nabla({\bf u}^{t}{\bf a}+{\rm div}({\bf u))}\right)\\
            &  & -\left(\alpha(q^{2})h G'(h){\rm div}({\bf u})+2 G(h)\alpha'(q^{2})((\nabla h)^{t}{\rm div}(h\nabla{\bf u}^{t})+{\rm div}({\bf u})q^{2})\right){\bf a}\\[2mm]
            & = & \alpha(q^{2}) G(h)\left(({\bf u}^{t}\nabla){\bf a}-\nabla({\bf u}^{t}{\bf a}+{\rm div}({\bf u))}\right)\\
            &  & -\left((\frac{h G'(h)}{G(h)}-1){\rm div}({\bf u})+\frac{2\alpha'(q^{2})}{\alpha(q^{2})}((\nabla h)^{t}{\rm div}(h\nabla{\bf u}^{t})+{\rm div}({\bf u})q^{2})\right)h{\bf v}.
        \end{eqnarray*}
    } Note that we have the relation
\[
    {\displaystyle \alpha(q^{2}){G}(h)\left(({\bf u}^{t}\nabla){\bf a}-\nabla({\bf u}^{t}{\bf a}+{\rm div}({\bf u))}\right)=-\frac{\alpha(q^{2}){G}(h)}{h}{\rm{div}}(h\nabla{\bf u}^{t})}
\]
and therefore, by using the relation $h \, {\bf v}=\alpha(q^{2}) \, {G}(h) \, {\bf a}$,
one finds {\setlength{\arraycolsep}{1pt}
        \begin{eqnarray*}
            \partial_{t}(h\boldsymbol{v})+{\rm div}(h\boldsymbol{v}\otimes\boldsymbol{u}) & = & -{\rm div}({\bf u})\left(\frac{h{G}'(h)}{{G}(h)}-1+\frac{2\alpha'(q^{2})q^{2}}{\alpha(q^{2})}\right)h{\bf v}\\
            &  & -\frac{2\alpha'(q^{2})}{\alpha(q^{2})}\left(\frac{h^{2}{\bf v}}{\alpha(q^{2}){G}(h)}\right)^{t}{\rm div}(h\nabla{\bf u}^{t})h{\bf v}-\frac{\alpha(q^{2}){G}(h)}{h}{\rm div}(h\nabla{\bf u}^{t}).\\[2mm]
            & = & -{\rm div}({\bf u})\left(\frac{h{G}'(h)}{{G}(h)}-1+\frac{2\alpha'(q^{2})q^{2}}{\alpha(q^{2})}\right)h{\bf v}\\
            &  & -\left(\frac{2\alpha'(q^{2})}{\alpha(q^{2})^{2}}\frac{h^{3}}{{G}(h)}{\bf v}\otimes{\bf v}+\frac{\alpha(q^{2}){G}(h)}{h}I_{d}\right){\rm div}(h\nabla{\bf u}^{t}).
        \end{eqnarray*}
    }
This yields the conclusion on the evolution of $h\boldsymbol{v}$.

\noindent \bigskip{}

\noindent \textit{Equation satisfied by $\boldsymbol{u}$.} Let us
first note that
\[
    \boldsymbol{p}=\frac{\boldsymbol{v}}{\alpha(q^{2})\sqrt{\dfrac{\sigma(h)}{h}}}
\]
and therefore
\[
    \nabla_{\mathbf{p}}E=\sigma(h)\left(\alpha(q^{2})^{2}+2\alpha(q^{2})\alpha'(q^{2})\|\boldsymbol{p}\|^{2}\right)\boldsymbol{p}=\sqrt{\sigma(h)}\sqrt{h}\left(\alpha(q^{2})+2\alpha'(q^{2})q^{2}\right)\boldsymbol{v}
\]
Next, we expand ${f} (h,\boldsymbol{v})\boldsymbol{v}$ and $g(h,\boldsymbol{v})\cdot \boldsymbol{v}$.
First, one has
\[
    {f}(h,\boldsymbol{v})\boldsymbol{v}=\sqrt{\sigma(h)}\sqrt{h}\left(2\frac{\alpha'(q^{2})h}{\alpha(q^{2})^{2}\sigma(h)}\|\boldsymbol{v}\|^{2}+\alpha(q^{2})\right)\boldsymbol{v}=\sqrt{\sigma(h)}\sqrt{h}\left(2\alpha'q^{2}+\alpha\right)\boldsymbol{v}=\nabla_{\mathbf{p}}E.
\]
Now we observe that
\[
    \boldsymbol{p}\cdot \nabla_{\mathbf{p}}E=\left(2\alpha'(q^{2})q^{2}+\alpha(q^{2})\right)\alpha(q^{2})q^{2}\sigma(h).
\]
This yields
\begin{eqnarray*}
    \boldsymbol {g} \left(h,\boldsymbol{v}\right) \cdot \boldsymbol{v} &  & =\left(\left(\frac{\sigma'(h)h}{2\sigma(h)}+\frac{1}{2}\right)\,+2\frac{\alpha'(q^{2})}{\alpha(q^{2})}q^{2}\right)h\|\boldsymbol{v}\|^{2}\\
    &  & =\left(\left(\frac{\sigma'(h)h}{2\sigma(h)}+\frac{1}{2}\right)\alpha(q^{2})\,+2\alpha'(q^{2})q^{2}\right)\alpha(q^{2})q^{2}\sigma(h)\\
    &  & =\left(2\alpha'(q^{2})q^{2}+\alpha(q^{2})\right)\alpha(q^{2})q^{2}\sigma(h)-\left(1-\frac{\sigma'(h)h}{\sigma(h)}\right)\frac{1}{2}(\alpha(q^{2}))^{2}q^{2}\sigma(h)
\end{eqnarray*}
and thus
\[
    \boldsymbol {g}  \left(h,\boldsymbol{v}\right) \cdot \boldsymbol{v}=\boldsymbol{p}^{t}\nabla_{\mathbf{p}}E-\left(\sigma(h)-h\sigma'(h)\right)\mathcal{E_{\mathrm{cap}}}\left(q\right).
\]
Note that the momentum conservation equation of \eqref{eq:Sys_SW_gen}
can be written as: %\[
%\partial_{t}\left(h\mathbf{u}\right)+div\left(h\mathbf{u}\otimes\mathbf{u}\right)+\nabla P=-div\left(\nabla h\otimes\nabla_{\mathbf{p}}E\right)+\nabla\left(hdiv\left(\nabla_{\mathbf{p}}E\right)\right)
%\]
%or
\[
    \partial_{t}\left(h\mathbf{u}\right)+{\rm div}\left(h\mathbf{u}\otimes\mathbf{u}\right)+\nabla\pi=-{\rm div}\left(\nabla h\otimes\nabla_{\mathbf{p}}E\right)+\nabla\left(h{\rm div}\left(\nabla_{\mathbf{p}}E\right)\right)+\nabla\left(\left(\sigma(h)-h\sigma'(h)\right)\mathcal{E_{\mathrm{cap}}}\left(\|\nabla h\|\right)\right).
\]
We now remark that
\[
    \operatorname{div}\left(h\nabla({f}  (h,\boldsymbol{v})\boldsymbol{v})^{t}\right)-\nabla(\boldsymbol{g}(h,\boldsymbol{v})\cdot\boldsymbol{v})=\operatorname{div}\left(h\nabla(\nabla_{\mathbf{p}}E)^{t}\right)-\nabla\left(\boldsymbol{p}^{t}\nabla_{\mathbf{p}}E-\left(\sigma(h)-h\sigma'(h)\right)\mathcal{E_{\mathrm{cap}}}\left(q\right)\right).
\]
Then, by taking $\boldsymbol{p}=\nabla h$, we obtain
\[
    \begin{array}{c}
        {\rm div}\left(h\nabla(\nabla_{\mathbf{p}}E)^{t}\right)-\nabla\left(\boldsymbol{p}^{t}\nabla_{\mathbf{p}}E\right)={\rm div}\left(h\nabla(\nabla_{\mathbf{p}}E)^{t}\right)-\nabla\left((\nabla h)^{t}\nabla_{\mathbf{p}}E\right) \\
        =-{\rm div}\left(\nabla h\otimes\nabla_{\mathbf{p}}E\right)+\nabla\left(h{\rm div}\left(\nabla_{\mathbf{p}}E\right)\right)
    \end{array}
\]
and consequently  the right-hand side of the
momentum equation in the augmented system is :
\[
    \operatorname{div}\left(h\nabla(f (h,\boldsymbol{v})\boldsymbol{v})^{t}\right)-\nabla(\boldsymbol {g} (h,\boldsymbol{v})^{t}\boldsymbol{v})=-{\rm div}\left(\nabla h\otimes\nabla_{\mathbf{p}}E\right)+\nabla\left(h{\rm div}\left(\nabla_{\mathbf{p}}E\right)\right)+\nabla\left(\left(\sigma(h)-h\sigma'(h)\right)\mathcal{E_{\mathrm{cap}}}\left(\|\nabla h\|\right)\right)
\]
and the momentum equation in the original system is satisfied, which
gives the conclusion on $\mathbf{u}$ for i).

\smallskip{}
Note that if $(h,\boldsymbol{u})$ is regular enough and the initial
velocity $\boldsymbol{v}_{0}$ satisfies
\[
    \boldsymbol{v}_{0}=\alpha(\left\Vert \nabla h_{0}\right\Vert ^{2})\sqrt{\frac{\sigma(h_{0})}{h_{0}}}\nabla h_{0}
\]
then $\boldsymbol{v}$ satisfies also \eqref{v} and $(h,\boldsymbol{u})$
solves the original system.

\noindent \medskip{}

\noindent \textbf{Part ii).} Recall that
\[
    E_{tot}(\boldsymbol{U})=\frac{1}{2h}\left(\left\Vert h\boldsymbol{u}\right\Vert ^{2}+\left\Vert h\boldsymbol{v}\right\Vert ^{2}\right)+\Phi\left(h\right)
\]
where $\boldsymbol{U}$ is given by \eqref{def_U_F} and
\[
    \left(\frac{\Phi}{h}\right)'=\frac{\pi}{h^{2}}.
\]

\noindent Let us consider the augmented system written as
\begin{equation}
    \partial_{t}\boldsymbol{U}+{\rm div}\left(\mathcal{F}\left(\boldsymbol{U}\right)\right)
    =\mathbf{M}\label{eq:Augmented_FHS}
\end{equation}
with the first order part given by
\[
    \boldsymbol{U}=\left(\begin{array}{c}
            h               \\
            h\boldsymbol{u} \\
            h\boldsymbol{v}
        \end{array}\right),\qquad \mathcal{F} \left(\boldsymbol{U}\right)=\left(\begin{array}{c}
            h\boldsymbol{u}                                             \\
            h\boldsymbol{u}\otimes\boldsymbol{u}+\pi\left(h\right)I_{d} \\
            h\boldsymbol{v}\otimes\boldsymbol{u}
        \end{array}\right)
\]
whereas the capillary term on the right hand side of \eqref{eq:Augmented_FHS}
is given by
\[
    \mathbf{M}=\left(\begin{array}{l}
            0                                                                                                                                          \\
            \operatorname{div}\left(h\nabla({f}(h,\boldsymbol{v})\boldsymbol{v})^{t}\right)-\nabla(\boldsymbol{g}(h,\boldsymbol{v})^{t}\boldsymbol{v}) \\
            - f(h,\boldsymbol{v})\operatorname{div}\left(h\nabla\boldsymbol{u}^{t}\right)\quad-\boldsymbol{g}(h,\boldsymbol{v})\operatorname{div}\boldsymbol{u}
        \end{array}\right).
\]
Note that the left-hand side of \eqref{eq:Augmented_FHS} is conservative
and hyperbolic and the right-hand site is skew-symmetric for the $L^2$
scalar product. This properties is real important to allow an appropriate
splitting method which preserve the energy conservation at the discrete level.

The entropy variable $\boldsymbol{V}$ for the first order part of \eqref{eq:Augmented_FHS}
is  given by
\[
    \boldsymbol{V}^t =\left(\nabla_{U}E_{tot}\right)^{t}=\left(-\frac{1}{2}\left(\left\Vert \boldsymbol{u}\right\Vert ^{2}+\left\Vert \boldsymbol{v}\right\Vert ^{2}\right)+\Phi'(h),\boldsymbol{u}^{t},\boldsymbol{v}^{t}\right).
\]
The energy equation is thus {\setlength{\arraycolsep}{1pt}
        \begin{eqnarray*}
            \partial_{t}E_{tot}+{\rm div}(\boldsymbol{u}\left(E_{tot}+\pi)\right)= &  & \left(\nabla_{U}E_{tot}\right)^{t}\mathcal{M}\\
            = &  & \boldsymbol{u}^{t}{\rm div}\left(h\nabla(f(h,\boldsymbol{v})\boldsymbol{v})^{t}\right)-\boldsymbol{u}^{t}\nabla(\boldsymbol{g}(h,\boldsymbol{v})^{t}\boldsymbol{v})\\
            &  & -\boldsymbol{v}^{t}f(h,\boldsymbol{v}){\rm div}\left(h\nabla\boldsymbol{u}^{t}\right)-\boldsymbol{v}^{t}
            \boldsymbol{g}(h,\boldsymbol{v}){\rm div}\left(\boldsymbol{u}\right)\\
            %%%%%%%%%%%%%%%%%%%%
            = &  & \boldsymbol{u}^{t}{\rm div}\left(h\nabla(f(h,\boldsymbol{v})\boldsymbol{v})^{t}\right)-(f(h,\boldsymbol{v}){\bf v})^{t}{\rm div}\left(h\nabla\boldsymbol{u}^{t}\right)\\
            &  & -{\rm div}\left(\boldsymbol{u}\boldsymbol{g}(h,\boldsymbol{v})^{t}\boldsymbol{v}\right)\\
            %-\boldsymbol{u}^{t}\nabla(g(h,\boldsymbol{v})^{t}\boldsymbol{v})-\boldsymbol{v}^{t}g(h,\boldsymbol{v})div\left(\boldsymbol{u}\right)
            %=
            %
            %
            = &  & {\rm div}\left(h(\boldsymbol{u}^{t}\nabla)(f(h,\boldsymbol{v})\boldsymbol{v})\right)-{\rm div}\left(h\left((f(h,\boldsymbol{v}){\bf v})^{t}\nabla\right)\boldsymbol{u}\right)-{\rm div}\left(\boldsymbol{u}\boldsymbol{g}(h,\boldsymbol{v})^{t}\boldsymbol{v}\right).
        \end{eqnarray*}
    } %Finally
%\[
%\begin{array}{c}
%\partial_{t}\left(E_{tot}\right)+div(\boldsymbol{u}\left(E_{tot}+\pi)\right)=\\
%div\left(h\boldsymbol{u}^{t}\nabla^{t}(f(h,\boldsymbol{v})\boldsymbol{v})\right)-div\left(h\boldsymbol{v}^{t}f(h,\boldsymbol{v})\nabla^{t}\boldsymbol{u}\right)-div\left(\boldsymbol{u}(g(h,\boldsymbol{v})^{t}\boldsymbol{v})\right)
%\end{array}
%\]
Recall that $f(h,\boldsymbol{v}) \cdot\boldsymbol{v}=\nabla_{\mathbf{p}}E_{tot}$
and $\boldsymbol{g}\left(h,\boldsymbol{v}\right)\cdot \boldsymbol{v}=\boldsymbol{p}^{t}\nabla_{\mathbf{p}}E_{tot}-\left(\sigma-h\sigma'\right)\mathcal{E_{\mathrm{cap}}}\left(q\right)$
and therefore
\[
    \begin{array}{c}
        \partial_{t}\left(E_{tot}\right)+{\rm div}\left(\boldsymbol{u}\left(E_{tot}+\pi\right)\right)=\left({\rm div}\left(h(\boldsymbol{u}^{t}\nabla)(\nabla_{\mathbf{p}}E_{tot})\right)-{\rm div}(h(\nabla_{\mathbf{p}}E_{tot}^{t}\nabla){\bf u})\right) \\
        -{\rm div}\left(\boldsymbol{u}(\boldsymbol{p}^{t}\nabla_{\mathbf{p}}E_{tot}-\left(\sigma-h\sigma'\right)\mathcal{E_{\mathrm{cap}}})\right).
    \end{array}
\]
By chosing $\boldsymbol{p}=\nabla h$, we easily verify that
\[
    \begin{array}{c}
        \left({\rm div}\left(h(\boldsymbol{u}^{t}\nabla)(\nabla_{\mathbf{p}}E_{tot})\right)-{\rm div}(h(\nabla_{\mathbf{p}}E_{tot}^{t}\nabla){\bf u})\right)-{\rm div}\left(\boldsymbol{u}\boldsymbol{p}^{t}\nabla_{\mathbf{p}}E_{tot}\right) \\
        =\left({\rm div}\left(h(\boldsymbol{u}^{t}\nabla)(\nabla_{\mathbf{p}}E_{tot})\right)-{\rm div}(h(\nabla_{\mathbf{p}}E_{tot}^{t}\nabla){\bf u})\right)-{\rm div}\left(\boldsymbol{u}(\nabla h)^{t}\nabla_{\mathbf{p}}E_{tot}\right)    \\
        ={\rm div}\left(h{\rm div}\left(\nabla_{\mathbf{p}}E_{tot}\right)\boldsymbol{u}\right)-{\rm div}\left({\rm div}\left(h\boldsymbol{u}\right)\nabla_{\mathbf{p}}E_{tot}\right).
    \end{array}
\]
Then we get
\[
    \begin{array}{c}
        \partial_{t}\left(E_{tot}\right)+{\rm div}\left(\boldsymbol{u}\left(E_{tot}+\pi-\left(\sigma-h\sigma'\right)\mathcal{E_{\mathrm{cap}}}\right)\right)={\rm div}\left(h{\rm div}\left(\nabla_{\mathbf{p}}E_{tot}\right)\boldsymbol{u}\right)-{\rm div}\left({\rm div}\left(h\boldsymbol{u}\right)\nabla_{\mathbf{p}}E_{tot}\right)\end{array}
\]
which is exactly the formulation (\ref{eq:Energie_K}) of the Energy
balance of the original system.

\bigskip{}
\bigskip{}

% \newpage
\bigskip
\bigskip
\noindent \textbf{Acknowledgments.} D. Bresch, M. Gisclon, N. Cellier and C. Ruyer-Quil have been supported by the Fraise project, grant ANR-16-CE06-0011 of the French National Research Agency (ANR) and with G.--L.~Richard by the project Optiwind through Horizon 2020/Clean Sky2 (call H2020-CS2-CFP06-2017-01) with Saint-Gobain.
This work was granted access to the HPC resources of CALMIP
supercomputing center under the allocation 2019-P1234.

\bigskip
\bigskip
\noindent \textbf{Authors coordinates:}
\noindent
\sc D. Bresch, M. Gisclon, G.-L. Richard,
\rm LAMA UMR5127 CNRS, Universit\'e Savoie Mont-Blanc,
73376 Le Bourget du Lac, France.

Emails: didier.bresch@univ-smb.fr, marguerite.gisclon@univ-smb.fr, gael.loic.richard@orange.fr

\smallskip

\noindent
\sc N. Cellier and C. Ruyer-Quil,
\rm LOCIE UMR5271 CNRS, Universit\'e Savoie Mont-Blanc,
73376 Le Bourget du Lac, France.

Emails: contact@nicolas-cellier.net, christian.ruyer-quil@univ-smb.fr

\smallskip

\noindent
\sc F. Couderc, P. Noble and J.P. Vila.
\rm Institut de Math\'ematiques de Toulouse,
UMR5219 CNRS, INSA de Toulouse,
31077 Toulouse cedex 4, France.

Emails: couderc@math.univ-toulouse.fr, pascal.noble@math.univ-toulouse.fr,  vila@insa-toulouse.fr

\end{document}